\documentclass[11pt,oneside,a4paper,mathscr]{amsart}
\usepackage{amsmath,amssymb}
\usepackage{tikz-cd}
\usepackage{adjustbox}
\usepackage{array}
\usepackage{amscd,amsmath,amssymb,euscript}
\usepackage[frame,cmtip,curve,arrow,matrix,line,graph]{xy}
\usepackage{tikz-cd}
\usepackage{hyperref}
\usepackage{bm}
\usepackage{here}
\usepackage{enumitem}
\hypersetup{hidelinks}
\usepackage{amsmath,amssymb}
\usepackage{tikz-cd}
\usepackage{xcolor}

\tikzset{
  mybox/.style={
    draw,
    rounded corners=3pt,
    fill=gray!8,
    inner sep=5pt,
    align=center
  },
  resultbox/.style={
    draw,
    rounded corners=3pt,
    fill=gray!14,
    inner sep=5pt,
    align=center
  }
}
\title[Semiorthogonal decompositions of limit categories]{Semiorthogonal decompositions of limit categories for symmetric stacks}
\author{Yukinobu Toda}

\newtheorem{thm}{Theorem}[section]
\newtheorem{claim}[thm]{Claim}

\newtheorem{cor}[thm]{Corollary}
\newtheorem{prop}[thm]{Proposition}

\newtheorem{lemma}[thm]{Lemma}

\theoremstyle{definition}
\newtheorem{defn}[thm]{Definition}

\newtheorem{thm*}[thm]{Theorem$^*$}
\newtheorem{remark}[thm]{Remark}

\newtheorem{example}[thm]{Example}

\newtheorem{assum}[thm]{Assumption}

\newcommand{\comment}[1]{}

\renewcommand{\leq}{\leqslant}
\renewcommand{\geq}{\geqslant}

\newcommand{\Y}{\mathcal{Y}}
\newcommand{\cS}{\mathcal{S}}

\newcommand{\LL}{\operatorname{L}}
\newcommand{\Pic}{\operatorname{Pic}}

\newcommand{\QCoh}{\operatorname{QCoh}}

\newcommand{\X}{\mathcal{X}}

\newcommand{\cH}{\mathcal{H}}

\newcommand{\cL}{\mathcal{L}}

\newcommand{\cC}{\mathcal{C}}

\newcommand{\cP}{\mathcal{P}}
\newcommand{\cE}{\mathcal{E}}

\newcommand{\cD}{\mathcal{D}}
\newcommand{\ev}{\mathrm{ev}}
\newcommand{\bgm}{B\mathbb{G}_m}

\newcommand{\diasquare}{\ar@{}[rd]|\square}

\newcommand{\fU}{\mathfrak{U}}

\newcommand{\Coh}{\operatorname{Coh}}

\newcommand{\id}{\operatorname{id}}

\newcommand{\Ind}{\operatorname{Ind}}
\newcommand{\IndCoh}{\operatorname{IndCoh}}
\newcommand{\IndL}{\operatorname{IndL}}
\newcommand{\Hom}{\operatorname{Hom}}

\newcommand{\Spec}{\operatorname{Spec}}

\newcommand{\GL}{\operatorname{GL}}

\newcommand{\ch}{\operatorname{ch}}

\newcommand{\inclusion}{\ar@<-0.3ex>@{^{(}->}[r]}
\newcommand{\iinclusion}{\ar@<-0.3ex>@{^{(}->}[rr]}
\newcommand{\linclusion}{\ar@<0.3ex>@{_{(}->}[l]}
\newcommand{\dinclusion}{\ar@<-0.3ex>@{^{(}->}[d]}
\newcommand{\uinclusion}{\ar@<-0.3ex>@{^{(}->}[u]}

\newcommand{\ssslash}{/\!\!/}

\newcommand{\wt}{\mathrm{wt}}
\newcommand{\Hig}{\mathrm{Higgs}}

\makeatletter
\newcommand{\colim@}[2]{%
  \vtop{\m@th\ialign{##\cr
    \hfil$#1\operator@font colim$\hfil\cr
    \noalign{\nointerlineskip\kern1.5\ex@}#2\cr
    \noalign{\nointerlineskip\kern-\ex@}\cr}}%
}
\newcommand{\colim}{%
  \mathop{\mathpalette\colim@{}}\nmlimits@
}
\makeatother

\usetikzlibrary{positioning,fit,calc}
\tikzstyle{block}=[draw=black, width=1cm, minimum height=2cm, align=center] 
\tikzstyle{block2}=[draw=black, text width=2cm, minimum height=1cm, align=center] 
\tikzstyle{block3}=[draw=black, text width=2cm, minimum height=1cm, align=center] 

\makeatletter

\@addtoreset{equation}{section}	
\makeatother

\begin{document}
\begin{abstract}
In our earlier work with Tudor P{\u a}durariu, we introduced the notion of limit categories for derived stacks with self-dual cotangent complexes and used them to formulate the Dolbeault geometric Langlands conjecture. In this paper, we generalize the notion of limit categories to a broader class of derived stacks, which we call symmetric stacks.

One of the motivations for this generalization is to formulate Dolbeault geometric Langlands-type correspondences for a broader class of moduli stacks. We establish several fundamental properties of limit categories for symmetric stacks, including compact generation, the existence of semiorthogonal decompositions, functoriality, and a Künneth formula.

In a companion paper, we apply the results of this paper to formulate a Dolbeault geometric Langlands-type conjecture for moduli stacks of sheaves on smooth surfaces with nef anticanonical divisors.
\end{abstract}
\maketitle
\setcounter{tocdepth}{2}
\tableofcontents

\section{Introduction}
\subsection{Overview and main results}

In our earlier paper with P{\u a}durariu~\cite{PTlim},
we introduced the notion of \textit{limit categories}
for derived stacks with self-dual cotangent complexes,
and formulated a version of the \textit{Dolbeault geometric Langlands conjecture}
for moduli stacks of Higgs bundles in terms of limit categories.
In the present paper, we extend the theory of limit categories
to a broader class of derived stacks, which we call \textit{symmetric stacks}.
We establish general results on the compact generation of these categories
and the existence of semiorthogonal decompositions associated with
$\Theta$-stratifications.
We also prove functoriality under suitable correspondences,
including Lagrangian correspondences between zero-shifted symplectic stacks,
as well as a K\"unneth formula.

One of the main motivations for this generalization is to formulate
\textit{Dolbeault geometric Langlands-type correspondences}
for a broader class of moduli stacks.
In~\cite{TodaSurfaceLimits}, we will apply this theory to moduli stacks
of pure one-dimensional sheaves on smooth surfaces with nef anticanonical divisors,
formulate a Dolbeault geometric Langlands-type conjecture in this setting,
and construct Hecke actions on limit and quasi-BPS categories.
This framework gives rise to new and interesting problems,
many of which are more tractable than the original Dolbeault geometric Langlands correspondence.
Furthermore, it will be applied in~\cite{PTtwist} to prove
the twisted Dolbeault geometric Langlands correspondence
for $\GL_2$ over the entire Hitchin base.

In~\cite{BPT}, we established the existence of semiorthogonal decompositions
for coherent sheaves on stacks admitting good moduli spaces, based on the combinatorics
of weight polytopes.
In Section~\ref{subsec:toy}, we examine a toy model
to explain the relationship between the semiorthogonal decompositions
arising from the construction in this paper and those of the dg-category
of coherent sheaves in~\cite{BPT}.

\subsection{Background}

Let $\mathfrak{M}$ be a QCA (quasi-compact with affine geometric stabilizers, see Subsection~\ref{subsec:notation0}) stack locally finitely presented over an algebraically closed field $k$ of characteristic zero.
In~\cite[Theorem~0.4.5]{MR3037900}, Drinfeld--Gaitsgory proved that the dg-category
\begin{align*}
\IndCoh(\mathfrak{M})
\end{align*}
is compactly generated, with compact objects given by $\Coh(\mathfrak{M})$.
However, such compact generation fails in general for non-quasi-compact stacks.

In our previous paper with Tudor P{\u a}durariu~\cite{PTlim},
for stacks $\mathfrak{M}$ whose cotangent complexes are self-dual,
we introduced the notion of a \textit{limit category}
\begin{align}\label{intro:dglim}
\IndL(\mathfrak{M})_{\delta} \subset \IndCoh(\mathfrak{M})
\end{align}
for $\delta \in \mathrm{Pic}(\mathfrak{M})_{\mathbb{Q}}$ and developed its general theory.
When $\mathfrak{M}$ is of cotangent type $\Omega_{\X}$ for a smooth stack $\X$,
the limit category is expected to provide an effective version of the classical limit of the dg-category of D-modules on $\X$.

We then proposed a formulation of the \textit{Dolbeault geometric Langlands conjecture}~\cite{DPlecture} in terms of limit categories of moduli stacks of Higgs bundles on curves.
It predicts an equivalence
\begin{align}\label{intro:DGLC}
\IndCoh(\Hig_{^L G}(w)^{\mathrm{ss}})_{-\chi} \simeq \IndL(\Hig_G(\chi))_w
\end{align}
between categories associated with the moduli stacks of $G$-Higgs bundles and ${}^{L}G$-Higgs bundles on a curve $C$; see~\cite[Conjecture~1.1]{PTlim} for details of the notation.
In the case of Higgs bundles, we proved in~\cite{PTlim} that the dg-category~(\ref{intro:dglim}) is compactly generated, even though 
$\Hig_G(\chi)$ is in general non-quasi-compact.
A conjectural equivalence (\ref{intro:DGLC}) is expected to be the
classical limit of the geometric Langlands correspondence, proposed and proved in~\cite{BD0, AG, GLC1, GLC2, GLC3, GLC4, GLC5}.

Moreover, its subcategory of compact objects admits a semiorthogonal decomposition into \textit{quasi-BPS categories}~\cite{PThiggs}, which categorify BPS invariants in Donaldson--Thomas theory~\cite{DM}.
In this way, the formulation~(\ref{intro:DGLC}) connects \textit{categorical DT theory}~\cite{T,HHR1,HHR2} with the geometric Langlands correspondence.

\subsection{Semiorthogonal decompositions of the limit category}
We extend the definition of the category~(\ref{intro:dglim}) to 
\textit{symmetric stacks} (see Definition~\ref{def:sym}), and generalize the results in~\cite{PTlim}. 
Let $\mathfrak{M}$ be a quasi-smooth, locally QCA and symmetric derived stack over $k$. 
For $\delta \in \mathrm{Pic}(\mathfrak{M})_{\mathbb{Q}}$ (or
more precisely, for the corresponding element of the dual of the component lattice), we 
define the subcategory
\begin{align}\label{intro:defL}
\IndL(\mathfrak{M})_{\delta} \subset \IndCoh(\mathfrak{M}). 
\end{align}

For example, 
a perfect complex $\cE$ on $\mathfrak{M}$ lies in this 
subcategory if and only if for any $\nu \colon \bgm \to \mathfrak{M}$ 
we have the weight condition 
\begin{align*}
    \wt(\nu^* \cE) \subset \left[\frac{1}{2}c_1(\nu^* \mathbb{L}_{\mathfrak{M}}^{<0}), \frac{1}{2}c_1(\nu^* \mathbb{L}_{\mathfrak{M}}^{>0})  \right]+c_1(\nu^* \delta). 
\end{align*}
We refer to Subsection~\ref{subsec:limit-definition} for details of the definition 
of the subcategory (\ref{intro:defL}). 
The subcategory 
\begin{align*}
    \LL(\mathfrak{M})_{\delta} \subset \IndL(\mathfrak{M})_{\delta}
\end{align*}
is defined to be the subcategory of compact objects. 

Suppose that the stack $\mathfrak{M}$ admits a (possibly infinite) locally finite $\Theta$-stratification~\cite{Halpinstab}, indexed by $I\subset\mathbb{R}_{\geq0}$ as in Definition~\ref{def:ndata}, with $I_{>0}=I\setminus\{0\}$,
\begin{align}\label{intro:thetastrata}
    \mathfrak{M}=\bigsqcup_{\mu \in I_{>0}} \mathfrak{S}_{\mu} \sqcup \mathfrak{M}^{\ell\text{-ss}}
\end{align}
with centers $\mathfrak{S}_{\mu} \to \mathfrak{Z}_{\mu}$, 
such that each $\Theta$-stratum 
\begin{align*}p_{\mu} \colon \mathfrak{S}_{\mu} \hookrightarrow \mathfrak{M}_{\preceq \mu}=\bigsqcup_{\mu'\in I\cap(0,\mu]} \mathfrak{S}_{\mu'} \sqcup \mathfrak{M}^{\ell\text{-ss}}
\end{align*}
satisfies Assumption~\ref{assum:dag}; it consists of the following four conditions:
\begin{enumerate}
    \item the stack $\mathfrak{M}_{\preceq \mu}$ is quasi-smooth, QCA, has a separated diagonal, and is symmetric; 
    \item the $\Theta$-stratification is defined by a linear form
    $\ell\in\mathrm{CL}(\mathfrak{M})_{\mathbb{Q}}^{\vee}$ and a positive definite rational quadratic norm $b$ on the component lattice, as in Definition~\ref{def:ndata};
    \item the dual obstruction spaces at the centers of the $\Theta$-strata have non-negative $\mathbb{G}_m$-weights; 
    \item the center $\mathfrak{Z}_{\mu}$ admits a good moduli space. 
\end{enumerate}

Under the above assumptions, 
we prove the following result, which generalizes~\cite[Theorem~1.2]{PTlim} (we refer to Subsection~\ref{subsec:component} for the notation $\delta_P$):
\begin{thm}\emph{(Theorem~\ref{thm:cgen})}\label{intro:thm1}
The dg-category $\IndL(\mathfrak{M})_{\delta}$ is compactly generated, 
and its subcategory of compact objects $\LL(\mathfrak{M})_{\delta}$ admits 
a semiorthogonal decomposition 
\begin{align}\label{intro:sod1}
    \LL(\mathfrak{M})_{\delta}=\left\langle \Upsilon_{\mu}\LL(\mathfrak{Z}_{\mu})_{\delta_{\mu}}, j_{!}\LL(\mathfrak{M}^{\ell\text{-ss}})_{\delta} : 
    \mu \in I_{>0} \right\rangle. 
\end{align}
Here 
$\delta_{\mu}=\delta|_{\mathfrak{Z}_{\mu}}+\delta_P(\det \mathbb{L}_{p_{\mu}}^{1/2})|_{\mathfrak{Z}_{\mu}}$, and 
$\Upsilon_{\mu}$ and $j_!$ are fully faithful functors 
\begin{align*}
    \Upsilon_{\mu} \colon \LL(\mathfrak{Z}_{\mu})_{\delta_{\mu}} \hookrightarrow \LL(\mathfrak{M})_{\delta}, \ j_{!} \colon 
    \LL(\mathfrak{M}^{\ell\text{-ss}})_{\delta}  \hookrightarrow 
    \LL(\mathfrak{M})_{\delta}. 
\end{align*}
The order of the semiorthogonal decomposition (\ref{intro:sod1}) is given by the decreasing order of $\mu$.
\end{thm}
The functor $j_{!}$ is the left adjoint functor of $j^* \colon \LL(\mathfrak{M})_{\delta} \to  \LL(\mathfrak{M}^{\ell\text{-ss}})_{\delta}$.
Note that the left adjoint functor of the pull-back of open immersions does not exist in the usual setting in algebraic geometry, e.g. 
$\Coh$, $\QCoh$, or $\IndCoh$. The existence of $j_{!}$ is one of the key 
properties of this theory. 

In our previous paper with Chenjing Bu and Tudor P{\u a}durariu~\cite{BPT}, 
we proved the existence of semiorthogonal decompositions of dg-categories of coherent sheaves on 
QCA stacks admitting good moduli spaces. In contrast, here we assume neither 
the QCA condition nor the existence of a good moduli space for $\mathfrak{M}$. 

The proof of Theorem~\ref{intro:thm1} is based on a local-to-global argument. For a single $\Theta$-stratum, we reduce the assertion, using \'etale slices, to a local model given by a derived zero locus of an equivariant section. After passing to a symmetric enlargement, the desired semiorthogonal decomposition follows from Koszul duality and the local result established in~\cite{PTlim}. We then iterate the one-stratum decomposition and use a limit-colimit argument to treat possibly infinitely many strata and prove compact generation.

We also use Theorem~\ref{intro:thm1} to prove the
K\"unneth formula for limit categories, see~Theorem~\ref{thm:kunneth}.

\subsection{Comparison of two semiorthogonal decompositions: a toy example}
In~\cite{BPT}, we proved that if $\mathfrak{M}$ satisfies condition~(i) above 
(quasi-smooth, QCA, with separated diagonal, and symmetric) and admits 
a \textit{good moduli space}
\begin{align*}
    \mathfrak{M} \to M
\end{align*}
then $\Coh(\mathfrak{M})$ admits a semiorthogonal decomposition 
of the form 
\begin{align}\label{intro:sod2}
     \Coh(\mathfrak{M})=\left\langle \ast_{\lambda}\LL(\mathfrak{M}^{\lambda})_{\delta_{\lambda}} : \lambda \in \lvert\mathrm{CL}(\mathfrak{M})_{\mathbb{Q}}\rvert\right\rangle.
 \end{align}
 We refer to Theorem~\ref{thm:BPT} for the details of the notation. 

\begin{figure}[htbp]
\centering
\[
\begin{tikzcd}[
ampersand replacement=\&,
column sep=0.7em,
row sep=2.7em,
cells={nodes={inner sep=2pt}}
]
\&
|[draw, rounded corners, inner sep=4pt]|
{\begin{array}{c}
\mathfrak{M}\\[-1pt]
\text{\scriptsize stack with additional}\\[-2pt]
\text{\scriptsize data and properties}
\end{array}}
\arrow[dl, dashed]
\arrow[dr, dashed]
\&
\\[0.3em]
|[draw, rounded corners, inner sep=4pt]|
{\begin{array}{c}
\mathfrak{M}^{\mathrm{ss}} \subset \mathfrak{M}\\[1pt]
\text{\scriptsize semistable locus}
\end{array}}
\arrow[d]
\&
\&
|[draw, rounded corners, inner sep=4pt]|
{\begin{array}{c}
\LL(\mathfrak{M})_{\delta}\\[1pt]
\text{\scriptsize limit category}
\end{array}}
\arrow[d]
\\
|[draw, rounded corners, inner sep=4pt]|
{\begin{array}{c}
\text{\scriptsize semiorthogonal decomposition of}\\
\Coh(\mathfrak{M}^{\mathrm{ss}})\\[-1pt]
\text{\scriptsize (combinatorics of}\\[-2pt]
\text{\scriptsize the weight polytope)}
\end{array}}
\arrow[rr, leftrightarrow, "\text{Langlands duality}" description]
\&
\&
|[draw, rounded corners, inner sep=4pt]|
{\begin{array}{c}
\text{\scriptsize semiorthogonal decomposition of}\\
\LL(\mathfrak{M})_{\delta}\\[-1pt]
\text{\scriptsize (geometry of}\\[-2pt]
\text{\scriptsize HN stratification)}
\end{array}}
\end{tikzcd}
\]
\caption{Conceptual overview}\label{fig:intro}
\end{figure}
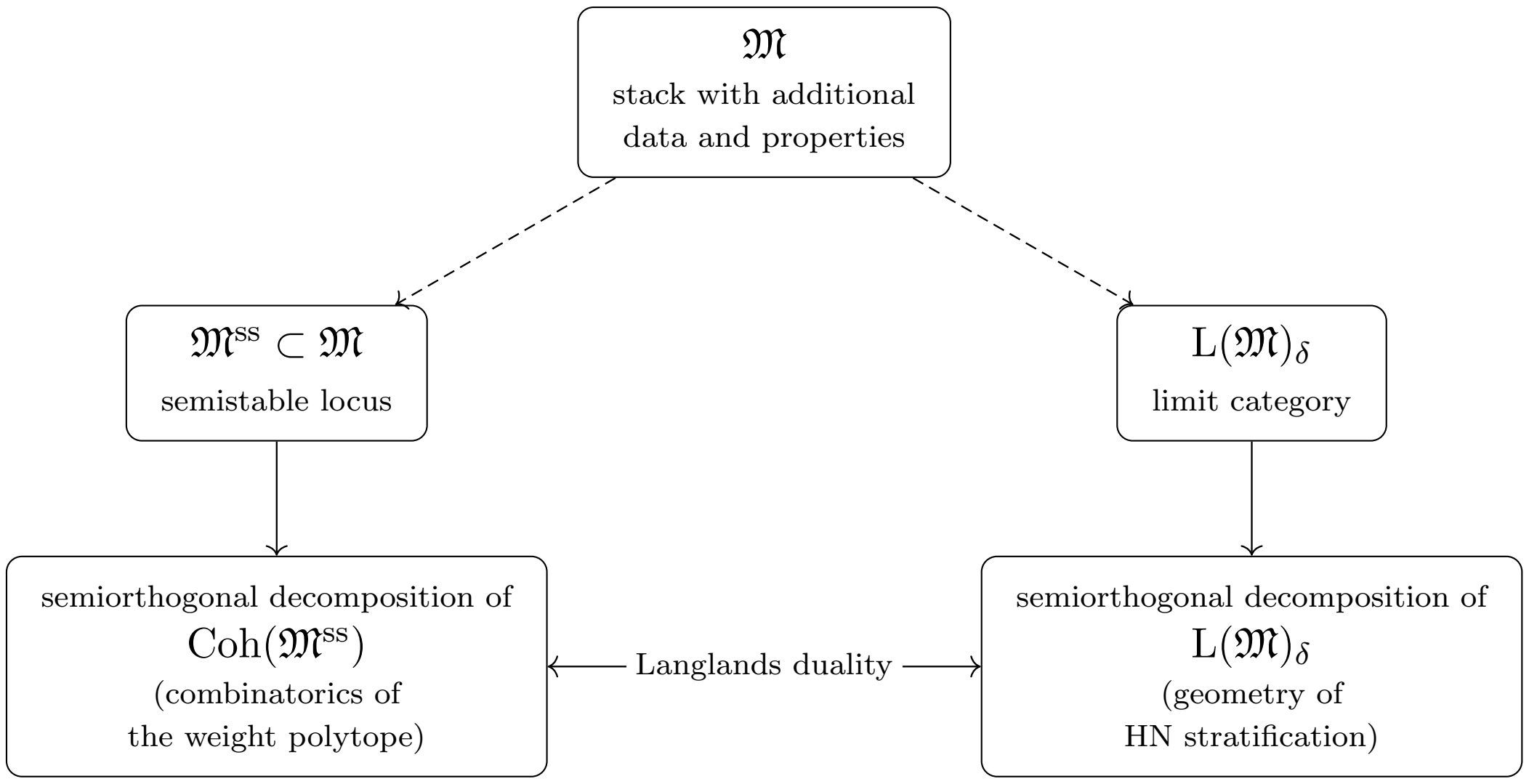

As mentioned above, 
the semiorthogonal decompositions~(\ref{intro:sod1}) and~(\ref{intro:sod2}) have different origins: the decomposition~(\ref{intro:sod2}) arises from the 
combinatorics of weight polytopes as in~\cite{SvdB,hls}, whereas~(\ref{intro:sod1})
arises from the geometry of the $\Theta$-stratification~(\ref{intro:thetastrata}). 
Our point of view is that these similar semiorthogonal decompositions, despite arising from different sources, 
are related by a Langlands-type duality;
see Figure~\ref{fig:intro}.

The following example illustrates this expectation.
On one side, consider the following stack with a 
good moduli space 
\begin{align*}
   \mathcal{Y}= \mathbb{A}_{xy}^2/\mathbb{G}_m
   \stackrel{\pi_{\mathcal{Y}}}{\to} \mathbb{A}^1
\end{align*}
where $\mathbb{G}_m$ acts on $\mathbb{A}_{xy}^2$ with weights $(1, -1)$. 
The semiorthogonal decomposition (\ref{intro:sod2}) is
\begin{align}\label{sod:introY}
    \Coh(\mathcal{Y})=\left\langle \mathcal{O}_{L^+}(-i), \mathcal{O}_{L^-}(i) : 
    i\geq 1, \pi_{\mathcal{Y}}^* \Coh(\mathbb{A}^1) \right\rangle. 
\end{align}
Here $L^+=\{y=0\}$ and $L^-=\{x=0\}$. 

On the other side, we construct a non-quasi-compact stack over $\mathbb{A}^1$
\begin{align}\label{intro:X}
    \X=\bigcup_{i\in 1/2+\mathbb{Z}} \X_i \stackrel{\pi_{\X}}{\to} \mathbb{A}^1.
\end{align}
Here each $\X_i$ is a copy of $\mathcal{Y}$, and 
$\X_i$ and $\X_{i+1}$ are glued along common open substacks 
isomorphic to $\mathbb{A}^1$. 
For $\delta=\delta_P(\mathcal{O}_{\mathcal{X}}(-1/2))$, the semiorthogonal decomposition~(\ref{intro:sod1}) of its limit category 
is 
\begin{align}\label{sod:introX}
    \LL(\X)_{\delta}=\left\langle A_i, B_i : i\in 1/2+\mathbb{Z}_{\geq 0}, s_{!}\Coh(\mathbb{A}^1) \right\rangle.
\end{align}
Here $s$ is a (non-closed) section of $\pi_{\mathcal{X}}$, and 
 \(A_i\) and \(B_i\) are certain \(!\)-pushforwards of (twisted) structure sheaves of \(L_i^+\subset\mathcal X_i\) and \(L_{-i}^-\subset\mathcal X_{-i}\), respectively.
The details of these computations will be given in Section~\ref{sec:toy}.

One can check that there is a rank-one Cohen--Macaulay sheaf 
$\cP$ on $\mathcal{Y} \times_{\mathbb{A}^1}\X$ 
whose associated Fourier--Mukai transform gives an equivalence; 
see Proposition~\ref{prop:equivX}:
\begin{align*}
    \Phi_{\mathcal{P}} \colon \Coh(\Y) \stackrel{\sim}{\to} 
    \LL(\X)_{\delta}.
\end{align*}
Moreover, under this equivalence, the semiorthogonal decompositions~(\ref{sod:introY}) and~(\ref{sod:introX}) are compatible.
Indeed, this is an example of a more general conjectural phenomenon, 
namely the \textit{Dolbeault geometric Langlands-type duality for 
surfaces} formulated and studied in the companion papers~\cite{TodaSurfaceLimits, PTtwist}.

\subsection{Further examples and $\infty$-flop functor}
We extend the example in the previous subsection to higher dimensions
by gluing copies of the quotient stack $(V\oplus V^{\vee})/\mathbb{G}_m$:
\begin{align*}
\X=\bigcup_{i\in \frac{1}{2}+\mathbb{Z}}\X_i,
\quad
\X_i\cong(V\oplus V^{\vee})/\mathbb{G}_m,
\end{align*}
where $V$ is an $n$-dimensional $\mathbb{G}_m$-representation of weight
one.
This gives a higher-dimensional generalization of~(\ref{intro:X}).
We describe its $\Theta$-stratification and the associated semiorthogonal
decomposition of $\LL(\X)_{\delta}$ given by
Theorem~\ref{intro:thm1}.
An interesting feature of this example is that the $!$-pushforward functor
\begin{align*}
j_! \colon \Coh(\X^{\mathrm{ss}}) \hookrightarrow \LL(\X)_{\delta}
\end{align*}
can be described by successively applying Bondal--Orlov derived
equivalences for standard flops~\cite{B-O2}. We call this the \textit{$\infty$-flop functor}.

Similarly, we consider the derived fiber product 
\begin{align*}
\mathfrak{X}_0:=\X\times_{\mathbb{A}^1}\{0\}
\end{align*}
of the moment map $h\colon\X\to\mathbb{A}^1$, and apply
Theorem~\ref{intro:thm1}.
The resulting functor
\begin{align*}
j_! \colon \Coh(\mathfrak{X}_0^{\mathrm{ss}})
\hookrightarrow \LL(\mathfrak{X}_0)_{\delta}
\end{align*}
is described analogously using successive Kawamata--Namikawa derived
equivalences for Mukai flops~\cite{Na,MR1949787}.
We refer to Subsections~\ref{subsec:toy-higher}
and~\ref{subsec:toy-moment-zero} for details.
\subsection{Organization of the paper}
In Section~\ref{sec:prelim}, we recall weights, filtered and graded objects, component lattices, numerical data, and the window theorem. In Section~\ref{sec:symmetric}, we define symmetric stacks and their limit categories, recall the semiorthogonal decomposition of~\cite{BPT}, and establish functoriality, including the functors associated with a $\Theta$-stratum. In Section~\ref{sec:decomposition}, we construct the semiorthogonal decomposition of the limit category and prove compact generation and the K\"unneth formula. In Section~\ref{sec:toy}, we study the two-dimensional toy model and its Fourier--Mukai equivalence, and then treat higher-dimensional toy models and their derived zero loci of moment maps via standard and Mukai flops.

\subsection{Acknowledgements}
The author would like to thank Tudor P{\u a}durariu, Tatsuki Kinjo,
Francesco Sala, Duiliu Diaconescu, Mauro Porta, Philasang Yoo,
and David Fang for helpful discussions.

The final stages of this work were carried out while the author was visiting the Morningside Center of Mathematics, Chinese Academy of Sciences, as a Distinguished Visiting Professor from October 18 to November 18, 2026. The author thanks the Morningside Center of Mathematics for its hospitality.

The author is supported by the World Premier International Research
Center Initiative (WPI Initiative), MEXT, Japan, the Inamori Research
Institute for Science, and JSPS KAKENHI Grant Number JP24H00180.

\subsection{AI disclosure}
ChatGPT-6 Astra was used to explore a proof of Lemma~\ref{lem:lift}, to assist with the exposition of Subsections~\ref{subsec:toy-higher} and~\ref{subsec:toy-moment-zero}, to check the proofs, and for editorial assistance. The author takes full responsibility for the correctness and content of the paper.

 \subsection{Notation and conventions}\label{subsec:notation0}
Figure~\ref{fig:notation} collects the main notation used in the paper,
together with the locations of its definitions.

    \subsubsection*{Conventions on stacks}
    Throughout the paper, $k$ is an algebraically closed field of characteristic zero.
All (derived) stacks are locally finitely presented over $k$, and all morphisms are over $k$.
For a stack $\X$, we write $\X(k)$ for the set of isomorphism classes of $k$-valued points of $\X$.
We denote its classical truncation by $\X^{\mathrm{cl}}$. 
A stack $\X$ is called \textit{QCA} if it is quasi-compact with affine geometric stabilizer groups and its classical inertia stack
is of finite presentation over $\X^{\mathrm{cl}}$, 
see~\cite[Definition~1.1.8]{MR3037900}. We say that $\X$ is \textit{locally QCA} if 
$\X$ admits a Zariski open cover by QCA open substacks. 

For a morphism of derived stacks $f\colon \X \to \Y$, let $\mathbb{L}_{f}$ and $\mathbb{T}_{f}$ denote the $f$-relative cotangent and tangent complexes. The morphism $f$ is called \textit{quasi-smooth} if $\mathbb{L}_{f}$ is perfect 
such that $\mathbb{L}_f|_{x}$ is 
of cohomological amplitude $[-1, 1]$ for any $x\in \mathcal{X}(k)$. In this case, we write $\omega_{f} :=\det \mathbb{L}_f$ and $\dim f:=\mathrm{rank}(\mathbb{L}_f)$.
All fiber products are derived fiber products, and we use the notation $\square$ for a 
Cartesian square. 

For a vector space $V$, we use the notation 
\begin{align*}
    V[-1]:=\Spec \mathrm{Sym}^*(V^{\vee}[1]).
\end{align*}

\subsubsection*{Notation for (ind-)coherent and quasi-coherent sheaves}\label{subsec:notation}
For a dg-category $\cC$, we denote by $\cC^{\mathrm{cp}} \subset \cC$ the subcategory of compact objects. The dg-category $\cC$ is called \textit{compactly generated} if $\cC=\Ind(\cC^{\mathrm{cp}})$, where $\Ind(-)$ is the ind-completion. 
A dg-functor $F \colon \cC_1 \to \cC_2$ is called \textit{continuous} 
if it commutes with direct sums. 

For compactly generated dg-categories $\cC_1, \cC_2$, their 
tensor product is denoted by $\cC_1 \otimes \cC_2$. 
It is given by the ind-completion of 
$\cC_1^{\mathrm{cp}}\boxtimes \cC_2^{\mathrm{cp}}$. 
The tensor product of the subcategories of compact objects
is defined by $\cC_1^{\mathrm{cp}} \otimes \cC_{2}^{\mathrm{cp}}:=(\cC_1 \otimes \cC_2)^{\mathrm{cp}}$.
It is split-generated by the objects $c_1 \boxtimes c_2$
for $c_i \in \cC_i^{\mathrm{cp}}$. 

We follow the notation of~\cite{MR3701352} for quasi-coherent and ind-coherent sheaves. 
For a derived stack $\X$, we denote by 
$\Coh(\X)$ the dg-category of coherent sheaves,
by $\QCoh(\X)$ the dg-category of quasi-coherent sheaves,
and by $\IndCoh(\X)$ the dg-category of ind-coherent sheaves. 
The heart of the standard t-structure on $\Coh(\X)$ is denoted by 
$\Coh^{\heartsuit}(\X)$. We will implicitly use the fact from~\cite{MR3037900} that $\IndCoh(\X)$ for a QCA stack $\X$
is compactly generated with compact objects $\Coh(\X)$. 
We write $\mathrm{Vect}:=\Coh(\Spec k)$, the dg-category of bounded 
complexes of 
finite-dimensional $k$-vector spaces. All functors in this paper
are derived; e.g. we write $f_{*}$ instead of $Rf_{*}$ for 
$f\colon \X \to \Y$. 

When $\mathfrak{M}$ is quasi-smooth, we denote by $\mathbb{D}_{\mathfrak{M}}$ its derived dual 
\begin{align*}
    \mathbb{D}_{\mathfrak{M}}(-):=\cH om(-, \mathcal{O}_{\mathfrak{M}}) \colon 
\Coh(\mathfrak{M})^{\mathrm{op}} \stackrel{\sim}{\to} \Coh(\mathfrak{M}).
\end{align*}

\subsubsection*{Semiorthogonal decomposition}
Let $\cD$ be a pretriangulated dg-category, let $I$ be a totally ordered set, and let $\cD_i \subset \cD$, for $i\in I$, be a family of 
pretriangulated full subcategories of $\cD$. 
We say that the subcategories $\cD_i$ form a semiorthogonal decomposition of $\cD$, denoted by 
\begin{align*}
    \cD=\langle \cD_i : i\in I\rangle
\end{align*}
if $\Hom(\cD_i, \cD_j)=0$ for $i>j$ and 
the subcategories $\cD_i$
generate $\cD$. 

\begin{figure}[p]
\centering
\footnotesize
\setlength{\tabcolsep}{4pt}
\renewcommand{\arraystretch}{1.08}
\begin{tabular}{|>{\raggedright\arraybackslash}p{.29\textwidth}|
                 >{\raggedright\arraybackslash}p{.43\textwidth}|
                 >{\raggedright\arraybackslash}p{.20\textwidth}|}
\hline
\textbf{Notation} & \textbf{Description} & \textbf{Location defined} \\\hline
$\X^{\mathrm{cl}}$ & Classical truncation of a derived stack &
Subsection~\ref{subsec:notation0} \\\hline
$V[-1]$ & Derived affine scheme $\Spec\mathrm{Sym}^*(V^\vee[1])$ &
Subsection~\ref{subsec:notation0} \\\hline
$\mathbb L_f,\ \mathbb T_f,\ \omega_f$ & Relative cotangent and tangent complexes; determinant of $\mathbb L_f$ &
Subsection~\ref{subsec:notation0} \\\hline
$\Coh,\ \QCoh,\ \IndCoh$ & Coherent, quasi-coherent and ind-coherent dg-categories &
Subsection~\ref{subsec:notation0} \\\hline
$\cC^{\mathrm{cp}},\ \Ind(\cC)$ & Compact objects; ind-completion &
Subsection~\ref{subsec:notation0} \\\hline
$\mathbb D_{\mathfrak M}$ & Derived dual with respect to $\mathcal O_{\mathfrak M}$ &
Subsection~\ref{subsec:notation0} \\\hline
$\Theta=\mathbb A^1/\mathbb G_m$ & Stack parametrizing filtrations &
Subsection~\ref{subsec:filtered} \\\hline
$\mathrm{Filt}(\mathfrak M),\ \mathrm{Grad}(\mathfrak M)$ & Stacks of filtered and graded objects &
Subsection~\ref{subsec:filtered} \\\hline
$\ev_0,\ \ev_1$ & Associated graded and evaluation at $1$ &
Equation~\eqref{dia:ev} \\\hline
$\mathcal Y^\lambda,\ \mathcal Y^{\lambda\geq0}$ & Fixed and attracting stacks for a cocharacter $\lambda$ &
Subsection~\ref{subsec:filtered} \\\hline
$\Coh(\X)_w,\ \wt(A)$ & Weight-$w$ subcategory; set of weights of $A$ &
Subsection~\ref{subsec:Gmwt} \\\hline
$A^{>0},\ \wt^{\mathrm{max}}(A),\ \wt^{\mathrm{min}}(A)$ & Positive-weight part; largest and smallest weights &
Subsection~\ref{subsec:Gmwt} \\\hline
$\mathrm{CL}(\mathfrak M)_{\mathbb Q}$ & Rational component lattice &
Subsection~\ref{subsec:component} \\\hline
$\delta_P,\ \delta$ & Character map from rational line bundles; linear parameter for limit categories &
Subsection~\ref{subsec:component}; Definition~\ref{defn:lim1} \\\hline
$\mathrm{Pic}(\mathfrak M)_{\mathbb Q}^{\mathrm{CL}}$ & Image of rational line bundles in the dual component lattice &
Subsection~\ref{subsec:component} \\\hline
$\nu_V^{\mathrm{reg}}$ & Quasi-smooth regularization of a graded point $\nu$ &
Subsection~\ref{subsec:limit-definition} \\\hline
$\IndL(\mathfrak M)_\delta$ & Ind-limit category on a locally QCA stack &
Definition~\ref{defn:lim2} \\\hline
$\LL(\mathfrak M)_\delta$ & Limit category; compact objects of $\IndL(\mathfrak M)_\delta$ &
Definitions~\ref{defn:lim1}, \ref{defn:lim2} \\\hline
$\widetilde{\LL}(\mathfrak M)_\delta$ & Limit of the categories $\LL(\mathfrak U)_\delta$ over QCA opens &
Equation~\eqref{def:Ltilde} \\\hline
$\ell,\ b,\ \mu(f)$ & Linear form, quadratic norm and slope of a filtration &
Subsection~\ref{subsec:numerical-data}; Equation~\eqref{def:slope:lb} \\\hline
$\mathfrak S_\mu,\ \mathfrak Z_\mu$ & $\Theta$-stratum and its center &
Subsection~\ref{subsec:theta-notation} \\\hline
$\mathfrak M_{\preceq\mu},\ \mathfrak M^{\ell\text{-ss}}$ & Open substack of the stratification; semistable locus &
Subsection~\ref{subsec:theta-notation}; Definition~\ref{def:ndata} \\\hline
$j_!,\ \Upsilon_\mu$ & Fully faithful extension and induction functors in the semiorthogonal decomposition of the limit category &
Theorem~\ref{thm:cgen} \\\hline
$\delta_\mu$ & The shifted parameter of the limit category on the center $\mathfrak Z_\mu$ &
Theorem~\ref{intro:thm1} \\\hline
$\X,\ \X^{\mathrm{ss}}$ & Higher-dimensional toy model and its semistable locus &
Subsection~\ref{subsec:toy-higher} \\\hline
$\mathfrak X_0,\ \mathfrak X_0^{\mathrm{ss}}$ & Derived zero locus of moment maps and their semistable locus &
Subsection~\ref{subsec:toy-moment-zero} \\\hline
$\mathcal A_{\mu,0},\ \mathcal B_{\mu,0}$ & Components of the semiorthogonal decomposition &
Proposition~\ref{prop:toy-moment-zero}\\\hline
$\mathrm{BO}_{i,t},\ \mathrm{BO}'_{i,t}$ & Bondal--Orlov standard-flop equivalences &
Remark~\ref{rmk:toy-BO} \\\hline
$\mathrm{KN}_{i,t},\ \mathrm{KN}'_{i,t}$ & Kawamata--Namikawa Mukai-flop equivalences &
Subsection~\ref{subsec:toy-moment-zero}  \\\hline
\end{tabular}
\caption{Main notation used in the paper}
\label{fig:notation}
\end{figure}

\section{Filtered and graded objects and the window theorem}\label{sec:prelim}
This section collects the notation for weights and filtered and graded objects,
recalls component lattices and numerical data, and states the window theorem.

\subsection{\texorpdfstring{Notation for $\mathbb{G}_m$-weights}{Notation for Gm-weights}}\label{subsec:Gmwt}
For a $\mathbb{G}_m$-gerbe $\X \to X$ of a quasi-compact stack $X$, there is an orthogonal weight 
decomposition 
\begin{align*}
    \Coh(\X)=\bigoplus_{w\in \mathbb{Z}}\Coh(\X)_w,
\end{align*}
where $\Coh(\X)_w$ is the subcategory of objects of $\mathbb{G}_m$-weight $w$. 
We have an equivalence 
\begin{align*}
    \Coh(\X)_{w} \simeq \Coh(X, \alpha^w)
\end{align*}
where $\alpha$ is the Brauer class of $X$ corresponding to the 
$\mathbb{G}_m$-gerbe $\X \to X$, and $\Coh(X, \alpha^w)$ is the dg-category of $\alpha^w$-twisted 
coherent sheaves~\cite{MR2309155}. 

For a non-zero object $A \in \Coh(\bgm)$,
let $A=\oplus_{w\in \mathbb{Z}} A_w$ 
be its weight decomposition. We write $A^{>0}=\oplus_{w>0} A_w$. 
We denote by $\wt(A)\subset \mathbb{Z}$ the set of integers $w$ such that $A_w\neq 0$, 
and by $\wt^{\mathrm{max}}(A)\in \mathbb{Z}$
(resp. $\wt^{\mathrm{min}}(A) \in \mathbb{Z}$) the maximal
(resp. minimal) element of $\wt(A)$. 
We also set $\wt^{\mathrm{min}}(0)=\infty$
and $\wt^{\mathrm{max}}(0)=-\infty$.
We will use the following map 
\begin{align*}
c_1 \colon K(\bgm)_{\mathbb{Q}} \to H^2(\bgm, \mathbb{Q})=\mathbb{Q}, 
\end{align*}
which is given by $c_1(A)=\sum_w w \cdot \dim A_w$, where $\dim A_w$ denotes the Euler characteristic of the weight-$w$ complex.

\subsection{Filtered and graded objects}\label{subsec:filtered}
We recall the notion of stacks of \textit{filtered and graded objects}, following~\cite{Halpinstab,HalpK32}.

Let $\Theta:=\mathbb{A}^1/\mathbb{G}_m$, where $\mathbb{G}_m$ acts on $\mathbb{A}^1$ by weight one. For a derived stack $\mathfrak{M}$, 
the stacks of filtered and graded objects are defined by 
\begin{align*}
    \mathrm{Filt}(\mathfrak{M}):=\mathrm{Map}(\Theta, \mathfrak{M}), \ 
    \mathrm{Grad}(\mathfrak{M}):=\mathrm{Map}(\bgm, \mathfrak{M}).
\end{align*}
The natural maps 
\begin{align*}
   \xymatrix{
   \bgm  \ar@<0.5ex>[r] & \Theta  \ar@<0.5ex>[l] & \mathrm{pt} \linclusion 
   }
\end{align*}
induce maps 
\begin{align}\label{dia:ev}
  \xymatrix{
   \mathrm{Grad}(\mathfrak{M})  \ar@<-0.5ex>[r]&  \ar@<-0.5ex>[l]_-{\ev_0} \mathrm{Filt}(\mathfrak{M}) \ar[r]^-{\ev_1} & \mathfrak{M}.
   }
\end{align}
A map $\nu \colon \bgm \to \mathfrak{M}$ is called a \textit{graded point}, and $f \colon \Theta \to \mathfrak{M}$ is called a \textit{filtration} of $f(1)$. For a filtration $f$, the map  $f_0:=\ev_0(f) \colon \bgm \to \mathfrak{M}$ is called its \textit{associated graded}. 

Let $G$ be a reductive group and let $Y$ be a representation of $G$. 
For a cocharacter $\lambda \colon \mathbb{G}_m \to G$, 
we denote by 
\begin{align*}
    Y^{\lambda} \subset Y^{\lambda \geq 0} \subset Y
\end{align*}
the $\lambda$-fixed locus and the subspace of non-negative $\lambda$-weights, 
respectively. 
We also denote by $G^{\lambda} \subset G^{\lambda \geq 0} \subset G$
the associated Levi and parabolic subgroups, respectively. 
For the quotient stack $\mathcal{Y}=Y/G$, we use the following notation 
\begin{align*}
    \mathcal{Y}^{\lambda}:=Y^{\lambda}/G^{\lambda}, \ 
     \mathcal{Y}^{\lambda \geq 0}:=Y^{\lambda \geq 0}/G^{\lambda \geq 0}.
\end{align*}
We have the diagram of natural maps
  \begin{align*}
 \xymatrix{
\mathcal{Y}^{\lambda} \ar@<-0.5ex>[r]&  \ar@<-0.5ex>[l] \mathcal{Y}^{\lambda \geq 0} \ar[r] & \mathcal{Y}
   }  
\end{align*}  
which is a special case of the diagram (\ref{dia:ev}) for $\mathfrak{M}=Y/G$ (and restricted to components of $\mathrm{Filt}(\mathfrak{M})$ and $\mathrm{Grad}(\mathfrak{M})$).

\subsection{Component lattice}\label{subsec:component}
Let $\mathfrak{M}$ be a derived stack.
Here we recall the notion of the \textit{component lattice} of $\mathfrak{M}$ introduced in~\cite{BHIK}. 
We refer to~\cite[Section~3]{BHIK} for more details. 

For a commutative ring $R$, a \textit{formal lattice} is a functor 
\begin{align*}
 F\colon    \mathrm{Lat}(R)^{\mathrm{op}} \to \mathrm{Set}
\end{align*}
where $\mathrm{Lat}(R)$ is the category of finitely generated free $R$-modules. 
Its underlying set is defined by $\lvert F \rvert :=F(R)$. 
The category of formal lattices is denoted by $\mathrm{fLat}(R)$. 
Any $R$-module $M$ is regarded as a formal lattice by the Yoneda embedding. 
For $F\in \mathrm{fLat}(R)$, we set  
\begin{align*}
    F^{\vee}:=\mathrm{Map}_{\mathrm{fLat}(R)}(F, R). 
\end{align*}

The \textit{component lattice} is defined as a formal lattice 
\begin{align*}
    \mathrm{CL}(\mathfrak{M})\in \mathrm{fLat}(\mathbb{Z}), \ 
    \Lambda \mapsto \pi_0(\mathrm{Map}(BT_{\Lambda}, \mathfrak{M})).
\end{align*}
Here 
$\pi_0$ denotes taking the set of connected components, $\mathrm{Map}$ denotes the derived mapping 
stack, and $T_{\Lambda}$ is the algebraic torus whose cocharacter lattice is $\Lambda$. 
There is also a rational version 
\begin{align*}
    \mathrm{CL}(\mathfrak{M})_{\mathbb{Q}} \in \mathrm{fLat}(\mathbb{Q})
\end{align*}
with its underlying set 
$\lvert \mathrm{CL}(\mathfrak{M})_{\mathbb{Q}}\rvert$. 

We have the natural map 
\begin{align*}
\delta_P \colon 
    \mathrm{Pic}(\mathfrak{M}) \to \mathrm{CL}(\mathfrak{M})^{\vee}, \ 
    \mathcal{L} \mapsto \left((\nu \colon BT_{\Lambda} \to \mathfrak{M})\mapsto 
  \nu^* \cL \in \Lambda^{\vee} \right).
\end{align*}
The above map naturally extends to 
\begin{align*}
\delta_P \colon 
     \mathrm{Pic}(\mathfrak{M})_{\mathbb{Q}} \to \mathrm{CL}(\mathfrak{M})_{\mathbb{Q}}^{\vee}.
\end{align*}
Its image is denoted by 
\begin{align*}
    \mathrm{Pic}(\mathfrak{M})^{\mathrm{CL}}_{\mathbb{Q}}\subset \mathrm{CL}(\mathfrak{M})_{\mathbb{Q}}^{\vee}.
\end{align*}
Below we often regard elements of $\mathrm{Pic}(\mathfrak{M})_{\mathbb{Q}}$
as elements of $\mathrm{CL}(\mathfrak{M})_{\mathbb{Q}}^{\vee}$ via the above map. 

For a map $f \colon \mathfrak{M}_1 \to \mathfrak{M}_2$, 
there is a pull-back 
\begin{align*}
    f^{*} \colon \mathrm{CL}(\mathfrak{M}_2)_{\mathbb{Q}}^{\vee} \to 
    \mathrm{CL}(\mathfrak{M}_1)_{\mathbb{Q}}^{\vee}
\end{align*}
induced by $\mathrm{Map}(BT_{\Lambda}, \mathfrak{M}_1) \to \mathrm{Map}(BT_{\Lambda}, \mathfrak{M}_2)$ by composing $BT_{\Lambda} \to \mathfrak{M}_1$ with $f$. 
It restricts to the map 
\begin{align*}
    f^* \colon   \mathrm{Pic}(\mathfrak{M}_2)^{\mathrm{CL}}_{\mathbb{Q}} \to   \mathrm{Pic}(\mathfrak{M}_1)^{\mathrm{CL}}_{\mathbb{Q}}.
\end{align*}

\subsection{$\Theta$-stratification}\label{subsec:theta-notation}
Recall from Subsection~\ref{subsec:filtered} the stacks of filtered and graded objects on  a given stack $\mathfrak{M}$.
Here we recall $\Theta$-strata and $\Theta$-stratifications, 
following~\cite{Halpinstab,HalpK32}. 

By definition, a \textit{$\Theta$-stratum} of a stack $\mathfrak{M}$ is an open and closed 
substack $\mathfrak{S} \subset \mathrm{Filt}(\mathfrak{M})$ such
that the map
\begin{align*}\ev_1|_{\mathfrak{S}} \colon \mathfrak{S} \to \mathfrak{M}
\end{align*}
is a closed 
immersion. 
For a $\Theta$-stratum $\mathfrak{S}$, there is an associated 
open and closed substack $\mathfrak{Z} \subset \mathrm{Grad}(\mathfrak{M})$, 
called the \textit{center} of $\mathfrak{S}$, such that 
the diagram (\ref{dia:ev}) restricts to the diagram 
\begin{align}\label{dia:ev2}
  \xymatrix{
   \mathfrak{Z}  \ar@<-0.5ex>[r]&  \ar@<-0.5ex>[l]_-{q} \mathfrak{S} \ar[r]^-{p} & \mathfrak{M}.
   }
\end{align}

Suppose that there is a well-ordered set $(I, \preceq)$
with minimal element $o\in I$,
and open substacks 
\begin{align*}
    \mathfrak{M}_{\preceq \mu} \subset \mathfrak{M}
\end{align*}
for each $\mu \in I$
such that the following conditions are satisfied:
\begin{itemize}[leftmargin=1.5em]
    \item For $\mu \preceq \mu'$, we have 
    $\mathfrak{M}_{\preceq \mu} \subset \mathfrak{M}_{\preceq \mu'}$. 
    \item There is a $\Theta$-stratum 
    $\mathfrak{S}_{\mu} \hookrightarrow \mathfrak{M}_{\preceq \mu}$
    with center $\mathfrak{S}_{\mu} \to \mathfrak{Z}_{\mu}$
such that 
\begin{align*}
    \mathfrak{M}_{\prec \mu}:=\bigcup_{\mu'\prec \mu} \mathfrak{M}_{\preceq \mu'}=\mathfrak{M}_{\preceq \mu} \setminus \mathfrak{S}_{\mu}.
\end{align*}
\item We have the open cover of $\mathfrak{M}$
\begin{align}\label{open:Mu}
    \mathfrak{M}=\bigcup_{\mu\in I}\mathfrak{M}_{\preceq \mu}.
\end{align}
\end{itemize}
Then we have the stratification (with a possibly infinite number of strata)
\begin{align}\label{Theta:M}
    \mathfrak{M}=\bigsqcup_{\mu\in I} \mathfrak{S}_{\mu}
\end{align}
where for simplicity we have set $\mathfrak{S}_{o}=\mathfrak{M}_{\preceq o}$.
The above stratification is called a \textit{$\Theta$-stratification}
as in~\cite[Definition~2.1.2]{Halpinstab}.

\subsection{$\Theta$-stratifications with respect to numerical data}\label{subsec:numerical-data}
We keep the setting of the previous subsection. 
A \textit{linear form} is an element
\begin{align*}
    \ell\in\mathrm{CL}(\mathfrak{M})_{\mathbb{Q}}^{\vee}.
\end{align*}
A \textit{rational quadratic form} $b$ on $\mathrm{CL}(\mathfrak{M})$
assigns to every $\gamma\in\mathrm{CL}(\mathfrak{M})(\Lambda)$ a quadratic
form $b_\gamma$ on $\Lambda_{\mathbb{Q}}$, compatibly with pull-back along
homomorphisms of lattices: for $a\colon\Lambda'\to\Lambda$,
\begin{align*}
    b_{a^*\gamma}(v)=b_\gamma(a(v)).
\end{align*}
 It extends uniquely to rational cocharacters.
For $\nu\colon\bgm\to\mathfrak{M}$, we abbreviate
\begin{align*}
    \ell(\nu):=\ell([\nu]),\quad
   \lVert\nu\rVert:= b(\nu)(:=b_{[\nu]}(1)).
\end{align*}
We call $b$ a \textit{positive definite quadratic norm} if
$b(\nu)>0$ for every non-degenerate graded point $\nu$.
Here non-degenerate means that the homomorphism on stabilizers is
non-trivial. 
We call the pair $(\ell, b)$ \textit{numerical data}. 
They pull back along representable 
morphisms of stacks.
\begin{remark}
For a $\mathbb{Q}$-line bundle $\mathcal{L}$ on $\mathfrak{M}$
and a class $B\in H^4(\mathfrak{M},\mathbb Q)$, 
there is an associated linear form $\ell$ and a rational quadratic form $b$ by 
\begin{align*}
    \ell=\delta_P(\mathcal{L}),\quad
    b(\nu)=\frac{\nu^*B}{u^2},\quad
    H^*(\bgm,\mathbb Q)=\mathbb Q[u].
\end{align*}
The pair $(\ell, b)$ gives numerical data if $b$ is positive definite. 
However there are numerical data which do not come from the above construction, 
see Remark~\ref{rmk:nonglobal}.
\end{remark}

For a map $f\colon\Theta\to\mathfrak{M}$ with non-degenerate
$f_0\colon\bgm\to\mathfrak{M}$, set
\begin{align}\label{def:slope:lb}
    \mu(f):=-\frac{\ell(f_0)}{\sqrt{b(f_0)}}.
\end{align}
The slope $\mu(f)$ is unchanged by replacing $f$
with the composition of the $r$-th power map $(-)^r\colon\Theta\to\Theta$ for $r>0$. 
\begin{defn}\label{def:stability}
A point $p\in\mathfrak{M}(k)$ 
is called $\ell$-semistable if for every map $g\colon\Theta\to\mathfrak{M}$ with
$g(1)\cong p$ and non-degenerate $g_0$, we have $\mu(g)\leq 0$.
\end{defn}

\begin{defn}\label{def:ndata}
Let $(\ell,b)$ be numerical data as above, with $b$ a positive
definite rational quadratic norm, and let $\mu$ be the slope
defined in~(\ref{def:slope:lb}).
We say that $(\ell,b)$ defines a $\Theta$-stratification of
$\mathfrak{M}$ if there is a $\Theta$-stratification
\begin{align*}
    \mathfrak{M}
    =\mathfrak{M}_{\preceq 0}
      \sqcup\bigsqcup_{c\in I\setminus\{0\}}\mathfrak{S}_c
\end{align*}
as in Subsection~\ref{subsec:theta-notation}, indexed by a set
$I\subset\mathbb{R}_{\geq 0}$ containing $0$, such that:
\begin{enumerate}
    \item For every $R\geq 0$, the set $I\cap[0,R]$ is finite.
    \item The open substack $\mathfrak{M}_{\preceq 0}$ is
    precisely the locus of $\ell$-semistable points.
    We denote it by $\mathfrak{M}^{\ell\text{-ss}}$.
    \item For every 
    $c\in I\setminus\{0\}$ and point 
    $f\in\mathfrak{S}_c(k)$ regarded as a map
    $f\colon\Theta\to\mathfrak{M}$, the graded point $f_0$
    is non-degenerate and $\mu(f)=c$.
    Moreover, for every map $g\colon\Theta\to\mathfrak{M}$
    with non-degenerate $g_0$ and $g(1)\cong f(1)$, we have
    $\mu(g)\leq c$. If $\mu(g)=c$, there exist positive integers
    $m,n$ such that
    \begin{align*}
        f\circ(-)^m \cong g\circ(-)^n.
    \end{align*}
\end{enumerate}
\end{defn}

\subsection{The window theorem}\label{subsec:window}
For a derived stack $\mathfrak{M}$, let $\mathfrak{S} \subset \mathrm{Filt}(\mathfrak{M})$ be a $\Theta$-stratum 
with center $\mathfrak{Z}$.
We have the diagram 
\begin{align}\label{theta:str}
    \xymatrix{
\mathfrak{S} \inclusion^-{p} \ar[d]^-{q} & \mathfrak{M} & \linclusion_-{j} \mathfrak{M}^{\circ}:=\mathfrak{M} \setminus \mathfrak{S} \\
\mathfrak{Z}. \ar@/^10pt/[u]^-{\sigma} \ar[ru]_-{\sigma'}& &
    }
\end{align}
If $\mathfrak{Z}$ is quasi-compact, 
there is a direct sum decomposition with respect to the canonical $\bgm$-action on $\mathfrak{Z}$
\begin{align*}
    \Coh(\mathfrak{Z})=\bigoplus_{w\in \mathbb{Z}} \Coh(\mathfrak{Z})_w.
\end{align*}

For $\delta \in \mathrm{CL}(\mathfrak{M})_{\mathbb{Q}}^{\vee}$, by the definition of the component lattice, its restriction to 
$\mathfrak{Z}$ gives a map 
\begin{align*}
    c_{\delta}\colon \pi_0(\mathfrak{Z}) \to \mathbb{Q}, \quad 
    \nu \mapsto \delta(\nu).
\end{align*}
When $\mathfrak{Z}$ is connected, the image of the above map in $\mathbb{Q}$ is called \textit{central weight}, 
and denoted by 
\begin{align*}\wt^{\mathrm{cent}}(\delta|_{\mathfrak{Z}}) \in \mathbb{Q}.
\end{align*}
If 
$\delta=\delta_P(\mathcal{L})$ for $\mathcal{L} \in \mathrm{Pic}(\mathfrak{M})$, the central weight of $\delta|_{\mathfrak{Z}}$ 
equals $w$ if and only if $\cL|_{\mathfrak{Z}} \in \Coh(\mathfrak{Z})_w$.

Suppose that $\mathfrak{M}$ is quasi-smooth. Then 
the stacks $\mathfrak{S}, \mathfrak{Z}$ are also quasi-smooth and
the map $q$ is quasi-smooth by~\cite[Lemma~3.1.4]{HalpK32}. 
Moreover $p$ is a closed immersion (in particular proper) by the definition of $\Theta$-stratum. 
Therefore the diagram (\ref{theta:str}) induces the functor 
\begin{align}\label{funct:up}
    \Upsilon \colon \Coh(\mathfrak{Z}) \stackrel{q^*}{\to} \Coh(\mathfrak{S}) \stackrel{p_{*}}{\to} \Coh(\mathfrak{M}).
\end{align}

We recall the window theorem from~\cite{HalpK32}:
\begin{thm}\emph{(\cite[Theorem~3.3.1]{HalpK32})}\label{thm:window}
Suppose that $\mathfrak{M}$ is quasi-smooth and QCA with separated diagonal,
and that $\mathfrak{Z}$ is connected. Assume that, for every
$z\in\mathfrak{Z}(k)$ with associated map $\nu\colon\bgm\to\mathfrak{M}$,
the weights of $\cH^{-1}(\nu^*\mathbb{L}_{\mathfrak{M}})$ are non-negative.
Then the functor (\ref{funct:up}) is fully faithful on each $\Coh(\mathfrak{Z})_i$, and 
there is a semiorthogonal decomposition for each $m\in \mathbb{R}$
\begin{align}\label{sod:dag}
&\Coh(\mathfrak{M})= \\
\notag&\langle \ldots, \Upsilon \Coh(\mathfrak{Z})_{\lfloor m \rfloor-1}, \Upsilon \Coh(\mathfrak{Z})_{\lfloor m \rfloor}, 
\mathcal{W}(m),\Upsilon \Coh(\mathfrak{Z})_{\lfloor m \rfloor+1}, \Upsilon \Coh(\mathfrak{Z})_{\lfloor m \rfloor+2}, \ldots \rangle.
\end{align}
Here $\lfloor m \rfloor\in  \mathbb{Z}$ is the round-down of $m$.
Moreover the following composition functor is an equivalence 
\begin{align*}
    \mathcal{W}(m) \subset \Coh(\mathfrak{M}) \stackrel{j^*}{\to} \Coh(\mathfrak{M}^{\circ}).
\end{align*}
\end{thm}
\begin{remark}\label{rmk:window}
The non-negativity hypothesis of the weights of 
$\cH^{-1}(\nu^* \mathbb{L}_{\mathfrak{M}})$ is the condition $(\dag)$ in~\cite{HalpK32}.
The theorem does not require numerical data or a good moduli space for the center. 

If we do not impose this condition, there is a
similar semiorthogonal decomposition for the singular support quotient of
$\Coh(\mathfrak{M})$ with respect to a $\Theta$-stratum of the
$(-1)$-shifted cotangent $\Omega_{\mathfrak{M}}[-1]$, see~\cite{Totheta}.
\end{remark}

The subcategory $\mathcal{W}(m)$ is characterized as follows. 
Here we use the dual $\mathbb{D}_{\mathfrak{M}}=\mathcal{H}om(-,\mathcal{O}_{\mathfrak{M}})$
of Subsection~\ref{subsec:notation0}. 
The following lemma follows from~\cite[Proposition~3.2.6]{HalpK32}
(also see
the proof of~\cite[Proposition~4.3.9]{HalpK32}):
\begin{lemma}\label{lem:wtmin}
In the setting of Theorem~\ref{thm:window}, an object $\cE \in \Coh(\mathfrak{M})$ lies in $\mathcal{W}(m)$ if and only if $\cE|_{\mathfrak{Z}}, \mathbb{D}_{\mathfrak{M}}(\cE)|_{\mathfrak{Z}} \in \mathrm{QCoh}(\mathfrak{Z})$ satisfy the inequalities 
\begin{align*}
     \wt^{\mathrm{min}}(\mathcal{E}|_{\mathfrak{Z}})> m, \
    \wt^{\mathrm{min}}(\mathbb{D}_{\mathfrak{M}}(\cE)|_{\mathfrak{Z}}) \geq -c_1(\mathbb{L}_{\mathfrak{M}}^{>0}|_{\mathfrak{Z}})-m.
\end{align*}
Here $\wt^{\mathrm{min}}$ is the minimal $\mathbb{G}_m$-weight 
with respect to the canonical $\bgm$-action on $\mathfrak{Z}$, 
and $(-)|_{\mathfrak{Z}}:=(\sigma')^*(-)$ for the map $\sigma'$ in (\ref{theta:str}). 

Equivalently for any 
$z\in \mathfrak{Z}(k)$ and the corresponding map 
$\nu \colon \bgm \to \mathfrak{M}$, we have 
\begin{align}\label{wtmin}
    \wt^{\mathrm{min}}(\nu^* \mathcal{E})> m, \
    \wt^{\mathrm{min}}(\nu^* \mathbb{D}_{\mathfrak{M}}(\cE)) \geq -c_1(\nu^* \mathbb{L}_{\mathfrak{M}}^{>0})-m.
\end{align}
\end{lemma}
 
\section{Limit categories for symmetric stacks}\label{sec:symmetric}
In this section, we introduce the notion of symmetric stacks and limit categories on them. 
We also establish functoriality of limit categories under
some class of correspondences (such as Lagrangian correspondences between zero-shifted symplectic stacks), including those associated with a $\Theta$-stratum.

\subsection{Definition of symmetric stacks}
The notion of symmetric stacks is introduced in~\cite[Definition~4.2.2]{BDIKP}
when the stacks admit good moduli spaces. 
Here we give a generalized (and simpler) notion of symmetric stacks, including those without good moduli spaces:
\begin{defn}\label{def:sym}
A quasi-smooth derived stack $\mathfrak{M}$ over $k$ is \textit{symmetric} if 
for any point $x\in \mathfrak{M}(k)$, the virtual $G_x:=\mathrm{Aut}(x)$-representation $\mathbb{T}_{\mathfrak{M}}|_{x}$ is self-dual, i.e. 
\begin{align*}
    [\mathbb{T}_{\mathfrak{M}}|_{x}]=[\mathbb{T}_{\mathfrak{M}}|_{x}^{\vee}] 
    \in K(BG_x).
\end{align*}
\end{defn}

When $\mathfrak{M}$ admits a good moduli space
and is either smooth, zero-shifted symplectic, or $(-1)$-shifted symplectic, 
the following lemma together with~\cite[Lemma~4.2.4]{BDIKP} 
shows that 
the above definition 
coincides with the one in~\cite[Definition~4.2.2]{BDIKP}:
\begin{lemma}\label{lem:sym1}
Suppose that $\mathfrak{M}$ is quasi-smooth, admits a good moduli space $\pi \colon \mathfrak{M}\to M$ and has an affine diagonal. 
Then $\mathfrak{M}$ is symmetric if and only if for any closed point $x\in \mathfrak{M}(k)$, the virtual representation \begin{align}\label{Tx}
[T^0_{\mathfrak{M}}|_x]
-
[T^1_{\mathfrak{M}}|_x]\in K(BG_x)
\end{align}
of the reductive group $G_x=\mathrm{Aut}(x)$ is self-dual. Here $T_{\mathfrak{M}}^i|_{x}:=\cH^i(\mathbb{T}_{\mathfrak{M}}|_{x})$.
\end{lemma}
\begin{proof}
Suppose first that $\mathfrak{M}$ is symmetric. Let
$x\in \mathfrak{M}(k)$ be a closed point. Then
$$
[\mathbb{T}_{\mathfrak{M}}|_x]
=
[T^0_{\mathfrak{M}}|_x]
-
[T^1_{\mathfrak{M}}|_x]
-
[\mathfrak{g}_x]
$$
is self-dual as a virtual $G_x$-representation. 
Since $G_x$ is reductive,
its adjoint representation $\mathfrak{g}_x$ is self-dual. It follows that
(\ref{Tx}) 
is also self-dual.

Conversely, suppose that (\ref{Tx})
is self-dual for every closed point
$x\in \mathfrak{M}(k)$. We show that $\mathfrak{M}$ is symmetric.
Take an arbitrary point
$y\in \mathfrak{M}(k)$, and let
$x\in \pi^{-1}(\pi(y))$ be the unique closed point in the fiber of the
good moduli space morphism
$
\pi\colon \mathfrak{M}\to M
$.
By the étale local description in~\cite[Lemma 4.2.3]{BPT}, étale locally
over $M$ at $\pi(y)=\pi(x)$, we may write
$$
\mathfrak{M}=s^{-1}(0)/G_x,
$$
where $s\colon U\to V$ is a $G_x$-equivariant map, $U$ is a
$G_x$-invariant open neighborhood of the origin in a
$G_x$-representation $Y$, and $V$ is a $G_x$-representation.
The point $y$ is represented by a point of $s^{-1}(0)$, and hence there is
an inclusion
$
i\colon G_y\hookrightarrow G_x$. 
Moreover, in $K(BG_y)$ we have
$$
[\mathbb{T}_{\mathfrak{M}}|_y]
=
[Y|_{G_y}]-[V|_{G_y}]-[\mathfrak{g}_x|_{G_y}]
=
i^*[\mathbb{T}_{\mathfrak{M}}|_x].
$$

On the other hand, by the assumption at the closed point $x$ and the
self-duality of $\mathfrak{g}_x$, the virtual $G_x$-representation
$$
[\mathbb{T}_{\mathfrak{M}}|_x]
=
[T^0_{\mathfrak{M}}|_x]
-
[T^1_{\mathfrak{M}}|_x]
-
[\mathfrak{g}_x]
$$
is self-dual. Its pull-back to $G_y$ is therefore self-dual, and hence
$[\mathbb{T}_{\mathfrak{M}}|_y]$ is self-dual. Since $y$ was arbitrary,
$\mathfrak{M}$ is symmetric.
\end{proof}

\begin{example}
    (1) If $\mathbb{L}_{\mathfrak{M}}$ is self-dual, then $\mathfrak{M}$ is symmetric. 
    It includes the case that $\mathfrak{M}$ is zero-shifted symplectic, and the cotangent 
    stacks $\Omega_{\mathcal{X}}$. For example, the moduli stacks of coherent sheaves on K3 or abelian surfaces and the moduli stacks of $G$-Higgs bundles for a reductive group $G$ on a smooth projective curve are symmetric. 

    (2) For any quasi-projective smooth surface $S$, the derived moduli stack of compactly supported one-dimensional coherent sheaves on $S$ is symmetric, see \cite[Proposition~3.1]{TodaSurfaceLimits}. Note that in this case, the cotangent complex of the moduli stack is not necessarily self-dual. 
\end{example}
\subsection{Limit categories for symmetric stacks}\label{subsec:limit-definition}
Let $\mathfrak{M}$ be a quasi-smooth symmetric stack over $k$. 
Consider a map
\begin{align*}
    \nu \colon \bgm \to \mathfrak{M}.
\end{align*}
A pair consisting of a finite-dimensional $\mathbb{G}_m$-representation $V$ 
and a map 
\begin{align*}
    \nu_V^{\mathrm{reg}} \colon V[-1]/\mathbb{G}_m \to \mathfrak{M}
\end{align*}
is called a \textit{regularization} of $\nu$ if $\nu_V^{\mathrm{reg}}$ is quasi-smooth 
and its composition 
\begin{align*}\bgm \to V[-1]/\mathbb{G}_m \stackrel{\nu_V^{\mathrm{reg}}}{\to} \mathfrak{M}
\end{align*}
is isomorphic to $\nu$. 
Such a regularization always exists; one may take a representation $V$ admitting a surjection
$V \twoheadrightarrow \cH^1(\nu^{*}\mathbb{T}_{\mathfrak{M}})$, see the proof of~\cite[Lemma~4.1.2]{BPT}.
We define 
the limit categories as follows: 
\begin{defn}\label{defn:lim1}
Let $\mathfrak{M}$ be a quasi-smooth, QCA and symmetric derived stack over $k$. For $\delta \in \mathrm{CL}(\mathfrak{M})_{\mathbb{Q}}^{\vee}$, we define the subcategory 
\begin{align}\label{def:limit:qca}
    \LL(\mathfrak{M})_{\delta} \subset \Coh(\mathfrak{M})
\end{align}
to consist of objects $\mathcal{E}$ such that for every map 
$\nu \colon \bgm \to \mathfrak{M}$ 
and every diagram 
\begin{align*}\bgm \stackrel{\iota}{\leftarrow} V[-1]/\mathbb{G}_m \stackrel{\nu_V^{\mathrm{reg}}}{\to} \mathfrak{M}
\end{align*}
where $\nu_V^{\mathrm{reg}}$ is a regularization of $\nu$ and $\iota$ is the projection, 
we have the condition 
\begin{align}\label{wt:condL}
   &\wt( \iota_{*}(\nu_V^{\mathrm{reg}})^*(\cE))
   \subset \\ 
 \notag  &\left[\frac{1}{2}c_1(\nu^* \mathbb{L}_{\mathfrak{M}}^{<0})
   +c_1((V^{\vee})^{<0}), \frac{1}{2}c_1(\nu^* \mathbb{L}_{\mathfrak{M}}^{>0})
   +c_1((V^{\vee})^{>0})
   \right]+\delta(\nu).
\end{align}
\end{defn}

\begin{remark}\label{lem:independent}
For a given $\nu \colon \bgm \to \mathfrak{M}$, the condition (\ref{wt:condL}) is independent of a choice of a regularization 
$\nu_V^{\mathrm{reg}}$, see~\cite[Lemma~4.1.2]{BPT}. 
In particular if $\mathbb{L}_{\mathfrak{M}}$ is self-dual, the 
definition of the subcategory (\ref{def:limit:qca}) coincides with 
the one in~\cite[Definition~3.3]{PTlim}. 
\end{remark}

\begin{remark}\label{rmk:closed-points}
If $\mathfrak{M}$ admits a good moduli space, it is enough in Definition~\ref{defn:lim1} to consider maps 
$\nu \colon \bgm \to \mathfrak{M}$ whose image is a closed point in $\mathfrak{M}(k)$. 
\end{remark}

\begin{remark}\label{lem:cond:equiv}
In the setting of Definition~\ref{defn:lim1},
if $\wt(\nu^* \cE)$ for $\nu^* \cE \in \QCoh(\bgm)$ is a finite set for any map $\nu \colon \bgm \to \mathfrak{M}$ (e.g. 
if $\cE$ is perfect), then the condition (\ref{wt:condL}) is equivalent to the condition that
for any such map $\nu$ we have 
    \begin{align}\label{wt:condL2}
    \wt(\nu^* \cE) \subset 
    \left[\frac{1}{2}c_1(\nu^* \mathbb{L}_{\mathfrak{M}}^{<0}), \frac{1}{2}c_1(\nu^* \mathbb{L}_{\mathfrak{M}}^{>0})
   \right]+\delta(\nu).
    \end{align}
    The proof is identical to~\cite[Lemma~6.1]{TodaGL2}.

    However in general $\wt(\nu^* \cE)$ is not a finite set, and the condition (\ref{wt:condL})
    is weaker than (\ref{wt:condL2}). 
\end{remark}

In the case of locally QCA stacks, which are not necessarily quasi-compact,  
we define the dg-category 
$\IndL(\mathfrak{M})_{\delta}$ 
which generalizes~\cite[Subsection~3.7]{PTlim} as follows: 
\begin{defn}\label{defn:lim2}
Let $\mathfrak{M}$ be a quasi-smooth, locally QCA and symmetric derived stack over $k$. For $\delta \in \mathrm{CL}(\mathfrak{M})_{\mathbb{Q}}^{\vee}$, we define 
\begin{align*}
    \IndL(\mathfrak{M})_{\delta}:=\lim_{\mathfrak{U} \subset \mathfrak{M}}\Ind(\LL(\mathfrak{U})_{\delta})
\end{align*}
where the limit is over QCA open substacks $\mathfrak{U}\subset \mathfrak{M}$, 
and with respect to the pull-backs for open immersions $j \colon \mathfrak{U}_1 \subset \mathfrak{U}_2$
    \begin{align*}
    j^* \colon 
        \Ind(\LL(\mathfrak{U}_2)_{\delta})
        \to \Ind(\LL(\mathfrak{U}_1)_{\delta}).
    \end{align*}
  
    The subcategory 
    \begin{align*}
\LL(\mathfrak{M})_{\delta} \subset \IndL(\mathfrak{M})_{\delta}
    \end{align*}
    is defined to be the subcategory of compact objects. 
    We also define 
    \begin{align}\label{def:Ltilde}
        \widetilde{\LL}(\mathfrak{M})_{\delta}:=\lim_{\mathfrak{U} \subset \mathfrak{M}}\LL(\mathfrak{U})_{\delta}
        \subset \IndL(\mathfrak{M})_{\delta}.
    \end{align}
\end{defn}
\begin{remark}\label{rmk:Ltilde}
    The subcategory (\ref{def:Ltilde}) is also defined as a subcategory of $\Coh(\mathfrak{M})$ by the weight condition (\ref{wt:condL}). In general, an object of $\widetilde{\LL}(\mathfrak{M})_{\delta}$ is not 
    compact in $\IndL(\mathfrak{M})_{\delta}$, see Corollary~\ref{cor:Ltilde}. 
\end{remark}

We now express the window condition of Lemma~\ref{lem:wtmin}
in terms of regularizations.
\begin{lemma}\label{lem:wtW}
In the setting of Lemma~\ref{lem:wtmin}, suppose in addition that $\mathfrak{M}$ is symmetric, and let 
 \begin{align*}\nu_V^{\mathrm{reg}} \colon V[-1]/\mathbb{G}_m\to \mathfrak{M}
 \end{align*}
 be a regularization of $\nu$.
The condition (\ref{wtmin}) holds if and only if 
\begin{align*}
    \wt(\iota_{*}(\nu_V^{\mathrm{reg}})^* \cE) \subset \left(m+c_1((V^{\vee})^{<0}), m+c_1((V^{\vee})^{>0})+c_1(\nu^* \mathbb{L}_{\mathfrak{M}}^{>0})  \right].
\end{align*}
\end{lemma}
\begin{proof}
    The equivalence of the lower bound of the weight follows from~\cite[Lemma~4.3.11]{HalpK32}, 
which says that 
\begin{align*}
    \wt^{\mathrm{min}}(\nu^*\cE)=\wt^{\mathrm{min}}(\iota_{*}(\nu_V^{\mathrm{reg}})^* \cE)-c_1((V^{\vee})^{<0}).
\end{align*}
    The upper bound is also obtained by the above 
    formula applied to $\mathbb{D}_{\mathfrak{M}}(\cE)$
    together with the Grothendieck duality for $\iota \colon V[-1]/\mathbb{G}_m \to \bgm$, 
    see the last part of the proof of~\cite[Lemma~4.3.10]{HalpK32}.
Here the symmetric condition gives $c_1(\nu^*\mathbb{L}_{\mathfrak{M}})=0$, so that
$-c_1(\nu^*\mathbb{L}_{\mathfrak{M}}^{<0})=c_1(\nu^*\mathbb{L}_{\mathfrak{M}}^{>0})$.
\end{proof}

\subsection{Semiorthogonal decomposition for stacks with good moduli spaces}
Here we recall the semiorthogonal decompositions of $\Coh(\mathfrak{M})$ 
into limit categories constructed in~\cite{BPT}, in the presence of 
good moduli spaces. 

Let $\mathfrak{M}$ be a quasi-smooth and QCA derived stack over $k$, 
which is symmetric, has a separated diagonal and 
a good moduli space
\begin{align*}
    \mathfrak{M} \to M.
\end{align*}
For each $\lambda \in \lvert\mathrm{CL}(\mathfrak{M})_{\mathbb{Q}}\rvert$, 
let 
\begin{align*}
    \mathfrak{M}^{\lambda} \subset \mathrm{Grad}(\mathfrak{M}), \ 
    \mathfrak{M}^{\lambda \geq 0} \subset \mathrm{Filt}(\mathfrak{M})
\end{align*}
be the corresponding components, such that the 
diagram of attracting loci (\ref{dia:ev}) restricts to the diagram 
\begin{align*}
    \mathfrak{M}^{\lambda} \stackrel{\mathrm{gr}}{\leftarrow} \mathfrak{M}^{\lambda \geq 0} 
    \stackrel{\ev}{\to} \mathfrak{M}.
\end{align*}
The map $\mathrm{gr}$ is quasi-smooth and, 
under the existence of a good moduli space, the map 
$\ev$ is proper. Therefore we have the functor 
(called \textit{Hall induction functor} in~\cite{BPT})
\begin{align*}
    \ast_{\lambda}=\ev_{*}\mathrm{gr}^* \colon \Coh(\mathfrak{M}^{\lambda}) \to \Coh(\mathfrak{M}).
\end{align*}
We recall the following result from~\cite{BPT}:
\begin{thm}\emph{(\cite[Theorem~4.1.3]{BPT})}\label{thm:BPT}
 For a given $\delta \in \mathrm{CL}(\mathfrak{M})_{\mathbb{Q}}^{\vee}$
 and an auxiliary quadratic norm $b_{\mathrm{aux}}$ on $\mathrm{CL}(\mathfrak{M})_{\mathbb{Q}}$, 
 there is a semiorthogonal decomposition 
 \begin{align*}
     \Coh(\mathfrak{M})=\left\langle \ast_{\lambda}\LL(\mathfrak{M}^{\lambda})_{\delta_{\lambda}} : \lambda \in \lvert \mathrm{CL}(\mathfrak{M})_{\mathbb{Q}}\rvert\right\rangle.
 \end{align*}
 Here each $*_{\lambda} \colon \LL(\mathfrak{M}^{\lambda})_{\delta_{\lambda}} \to \Coh(\mathfrak{M})$
 is fully faithful, and $\delta_{\lambda}\in \mathrm{CL}(\mathfrak{M}^{\lambda})_{\mathbb{Q}}^{\vee}$ is defined by
 \begin{align*}
     \delta_{\lambda}(\mu)=\delta(\mu)+\frac{1}{2}c_1(\mu^* \det \mathbb{L}_{\ev}|_{\mathfrak{M}^{\lambda}})-b_{\mathrm{aux}}(\lambda, \mu)
 \end{align*}
 for $\mu \in \lvert \mathrm{CL}(\mathfrak{M}^{\lambda})_{\mathbb{Q}}\rvert$.
 Here $b_{\mathrm{aux}}(\lambda, \mu)$ is evaluated on the two-dimensional face determined by 
 $\lambda$ and $\mu$, as in~\cite[Theorem~3.2.3]{BPT}.
\end{thm}

The following example will be compared with the semiorthogonal decomposition of the limit category in Section~\ref{sec:toy}:
\begin{example}\label{exam:1}
Let $\mathcal{Y}=\mathbb{A}^2/\mathbb{G}_m$ where 
the $\mathbb{G}_m$-action is given by weight $(1, -1)$. 
For $\lambda=\id \colon \mathbb{G}_m \to \mathbb{G}_m$, let 
\begin{align*}
    L^{+}=\mathcal{Y}^{\lambda \geq 0}, \ L^-=\mathcal{Y}^{\lambda^{-1} \geq 0}. 
\end{align*}
    Let $\mathcal{O}(1) \in \mathrm{Pic}(\mathcal{Y})$ be the pull-back of 
    $\mathcal{O}_{\bgm}(1)$ by $\mathbb{A}^2/\mathbb{G}_m \to \bgm$. 
    For each $i\in \mathbb{Z}_{>0}$, let 
    \begin{align*}
        \cC(i):=\langle \mathcal{O}_{L^+}(-i), \mathcal{O}_{L^-}(i) \rangle
        \subset \Coh(\mathcal{Y})
    \end{align*}
    where $\mathcal{O}_{L^+}(-i)$ and $\mathcal{O}_{L^-}(i)$ are 
    orthogonal exceptional objects, so it is equivalent to $\mathrm{Vect} \oplus \mathrm{Vect}$. 

We apply Theorem~\ref{thm:BPT} with $\delta=0$ and any standard quadratic norm $b_{\mathrm{aux}}$; it gives
    \begin{align}\label{sod:A2gm}
        \Coh(\mathcal{Y})=
        \left\langle \ldots, \cC(2), \cC(1), \LL(\mathcal{Y})_0  \right\rangle. 
    \end{align}
    Moreover $\LL(\mathcal{Y})_0$ is generated by $\mathcal{O}_{\mathcal{Y}}$. 

    By direct calculation, the Hom spaces in the other directions
    are also computed as follows (here $i, j>0$): 
    \begin{align}\label{calc1}
       & \Hom(\mathcal{O}_{L^+}(-i), \mathcal{O}_{\mathcal{Y}})
        = \Hom(\mathcal{O}_{L^-}(i), \mathcal{O}_{\mathcal{Y}})
        =k[-1], \\
       \notag &\Hom(\mathcal{O}_{L^+}(-i), \mathcal{O}_{L^+}(-j))
        =\Hom(\mathcal{O}_{L^-}(i), \mathcal{O}_{L^-}(j))=k\oplus k[-1], \ (i>j), \\
     \notag   &\Hom(\mathcal{O}_{L^+}(-i), \mathcal{O}_{L^-}(j))
        =\Hom(\mathcal{O}_{L^-}(j), \mathcal{O}_{L^+}(-i))=0.
    \end{align}
\end{example}

\subsection{Functors given by correspondences}\label{subsec:functors}\label{sec:functor}
We establish functoriality for limit categories, generalizing the results
for cotangent stacks in~\cite[Section~5]{PTlim}.
We consider a diagram of derived stacks over $k$
\begin{align}\label{dia:M12}
    \mathfrak{M}_1 \stackrel{q}{\leftarrow} \mathfrak{N} \stackrel{p}{\to} \mathfrak{M}_2
\end{align}
such that 
\begin{itemize}[leftmargin=1.5em]
    \item $\mathfrak{M}_i$ are quasi-smooth, QCA and symmetric stacks, 
and $\mathfrak{N}$ is quasi-smooth;
\item $q$ is quasi-smooth, and $p$ is projective; 
\item for any map $g \colon \bgm \to \mathfrak{N}$, we have 
\begin{align}\label{c1zero}
    c_1(g^*\mathbb{L}_p^{>0})= c_1(g^*\mathbb{L}_q^{<0}), \ c_1(g^*\mathbb{L}_p^{<0})= c_1(g^*\mathbb{L}_q^{>0}).
\end{align}
\end{itemize}
By the assumption on $p$ and $q$, we have the functors
\begin{align}\label{funct:coh}
p_{*}q^* \colon \Coh(\mathfrak{M}_1) \stackrel{q^*}{\to} \Coh(\mathfrak{N}) \stackrel{p_{*}}{\to} \Coh(\mathfrak{M}_2).
\end{align}

Under the above setting, we 
have the following proposition: 
\begin{prop}\label{prop:functL}
Suppose that there exist 
$\delta_i \in \mathrm{CL}(\mathfrak{M}_i)_{\mathbb{Q}}^{\vee}$
such that 
\begin{align}\label{pqdelta}
    q^* \delta_1=p^* \delta_2+\frac{1}{2}\delta_P(\det \mathbb{L}_p) \in \mathrm{CL}(\mathfrak{N})_{\mathbb{Q}}^{\vee}.
\end{align}
Then the functor (\ref{funct:coh})
restricts to the functor 
\begin{align*}
 p_{*}q^* \colon   \LL(\mathfrak{M}_1)_{\delta_1} \to \LL(\mathfrak{M}_2)_{\delta_2}.
\end{align*}
\end{prop}
\begin{proof}
  We take $A \in \LL(\mathfrak{M}_1)_{\delta_1}$, 
  and a finite-dimensional $\mathbb{G}_m$-representation $V$. 
  Let $\nu_V^{\mathrm{reg}}$ be a quasi-smooth map 
  \begin{align*}
      \nu_V^{\mathrm{reg}} \colon V[-1]/\mathbb{G}_m \to \mathfrak{M}_2, 
  \end{align*}
  which is a regularization of $\nu \colon \bgm \to \mathfrak{M}_2$.
  We need to show that 
  \begin{align}\label{wt:condf}
    & \wt(\iota_{*}(\nu_V^{\mathrm{reg}})^* p_{*}q^* A) \subset \\
   & \notag  \left[\frac{1}{2}c_1(\nu^* \mathbb{L}_{\mathfrak{M}_2}^{<0})
   +c_1((V^{\vee})^{<0}), \frac{1}{2}c_1(\nu^* \mathbb{L}_{\mathfrak{M}_2}^{>0})
   +c_1((V^{\vee})^{>0})
   \right] 
   +\delta_2(\nu).
  \end{align}
  
  We have the following commutative diagram 
  \begin{align}\label{dia:N}
      \xymatrix{
\bgm \ar@{=}[d] & N/\mathbb{G}_m  \ar@/^15pt/[rr]^-{q'} \ar[r]_-{r} \ar[d]_-{p'} \ar[l]_-{\iota'} \diasquare & \mathfrak{N} \ar[d]_-{p} 
\ar[r]_-{q} & \mathfrak{M}_1 \\
\bgm & V[-1]/\mathbb{G}_m \ar[l]_-{\iota} \ar[r]^-{\nu_V^{\mathrm{reg}}} & \mathfrak{M}_2. 
      }
  \end{align}
  Here $N/\mathbb{G}_m$ is defined by the above Cartesian square. Note
  that $N$ is a quasi-smooth derived scheme over $k$
  (as $\nu_V^{\mathrm{reg}}$ is quasi-smooth, $\mathfrak{N}$ is also quasi-smooth and $p'$ is projective)
  with $\mathbb{G}_m$-action. 
  By the above diagram and the base change, we have 
  \begin{align}\label{isom:pqA}
      \iota_{*}(\nu_V^{\mathrm{reg}})^* p_{*}q^* A &\cong \iota_{*}p'_{*}r^*q^*A 
      \cong \iota'_{*}(q')^*A.
  \end{align}
  Since $p'$ is projective,
  there is a commutative diagram 
  \begin{align*}
      \xymatrix{
N/\mathbb{G}_m \ar[rd]^-{\iota'} \inclusion^-{i} \ar[d]_-{p'} & \mathbb{P}^n/\mathbb{G}_m \ar[d]^-{p''} \\
V[-1]/\mathbb{G}_m \ar[r]^-{\iota} & \bgm.
      }
  \end{align*}
  Here $\mathbb{G}_m$ acts on $\mathbb{P}^n$ linearly and $i$ is a closed immersion. Note that $i$ is quasi-smooth, as $N$ is quasi-smooth and $\mathbb{P}^n$ is smooth.  
By the above diagram, we have 
\begin{align}\label{isom:pprime}
   \iota'_{*}(q')^*A \cong p''_{*} i_{*}(q')^*A.  
\end{align}

Let $x\in N\subset \mathbb{P}^n$ be a $\mathbb{G}_m$-fixed closed point, 
and take the corresponding map $\lambda \colon \bgm \to \mathbb{P}^n/\mathbb{G}_m$. 
We have the following diagram 
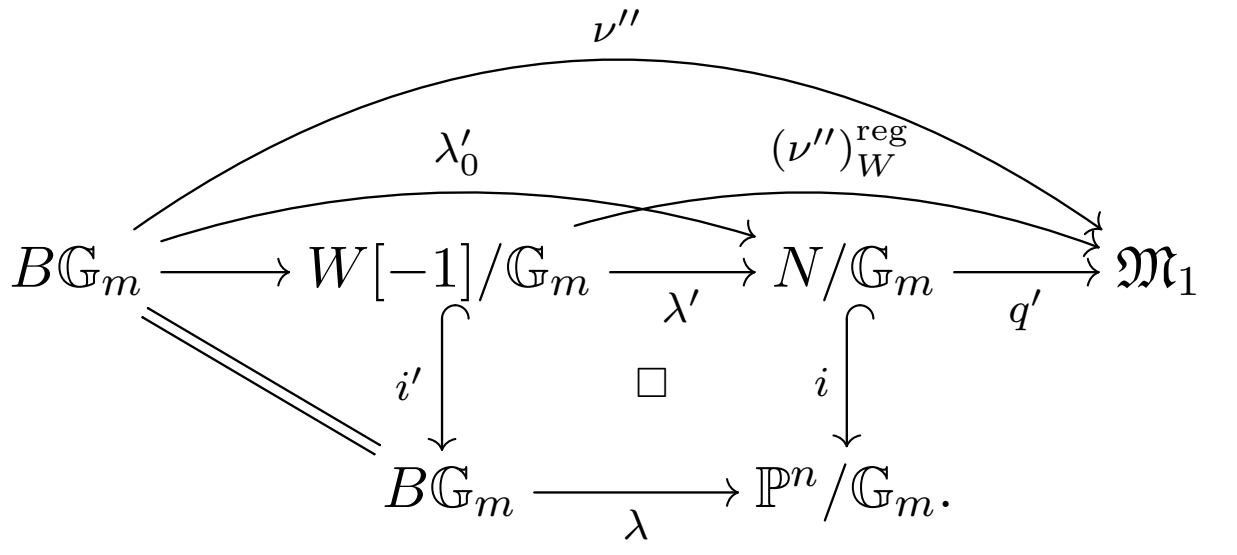
\begin{align*}
\xymatrix{
 \bgm \ar[r] \ar@{=}[rd] \ar@/^15pt/[rr]^-{\lambda_0'} \ar@/^40pt/[rrr]^-{\nu''}& W[-1]/\mathbb{G}_m \ar@/^15pt/[rr]^-{(\nu'')_W^{\mathrm{reg}}}\diasquare \dinclusion_-{i'} \ar[r]_-{\lambda'} & N/\mathbb{G}_m \dinclusion_-{i} \ar[r]_-{q'} & \mathfrak{M}_1 \\
 & \bgm \ar[r]_-{\lambda} & \mathbb{P}^n/\mathbb{G}_m.  &
 }
\end{align*}
Here $W=\cH^1(\mathbb{T}_i|_{x})$ is regarded as a $\mathbb{G}_m$-representation. We have the isomorphism 
\begin{align}\label{id:lambda}
\lambda^* i_{*}(q')^* A \cong i'_{*}(\lambda')^*(q')^*A \cong i'_{*}((\nu'')_W^{\mathrm{reg}})^*A.
\end{align}
The map $(\nu'')_W^{\mathrm{reg}}$ is quasi-smooth, since both $\lambda'$ and $q'$ are quasi-smooth. Therefore $(\nu'')_W^{\mathrm{reg}}$ is a regularization of $\nu''$. 
Since $A\in \LL(\mathfrak{M}_1)_{\delta_1}$, we have 
\begin{align}\label{wt:iprime}
  &  \wt(i'_{*}((\nu'')_W^{\mathrm{reg}})^*A) \subset \\
 \notag & \left[\frac{1}{2}c_1((\nu'')^* \mathbb{L}_{\mathfrak{M}_1}^{<0})
   +c_1((W^{\vee})^{<0}), \frac{1}{2}c_1((\nu'')^* \mathbb{L}_{\mathfrak{M}_1}^{>0})
   +c_1((W^{\vee})^{>0})
   \right]+\delta_1(\nu'').
\end{align}
We also have 
\begin{align}\label{equal:delta}
\delta_1(\nu'')=\delta_2(\nu)+\frac{1}{2}c_1(\det \mathbb{L}_p|_{x})
\end{align}
by the condition (\ref{pqdelta}). 

By the distinguished triangle 
\begin{align*}
    \Omega_{\mathbb{P}^n}|_{x} \to \mathbb{L}_{N}|_{x} \to \mathbb{L}_i|_{x}=W^{\vee}[1]
\end{align*}
we have the identity 
\begin{align*}
    c_1(\Omega_{\mathbb{P}^n}|_{x}^{<0})-c_1((W^{\vee})^{<0})=c_1(\mathbb{L}_{N}|_{x}^{<0}).
\end{align*}
By the distinguished triangle 
\begin{align*}
    (p')^* \mathbb{L}_{V[-1]}|_{x} \to \mathbb{L}_{N}|_{x} \to \mathbb{L}_p|_{x}
\end{align*}
we have 
\begin{align*}
    c_1(\mathbb{L}_N|_{x}^{<0})=-c_1((V^{\vee})^{<0})+c_1(\mathbb{L}_{p}|_{x}^{<0}).
\end{align*}
It follows that 
\begin{align*}
  &\frac{1}{2}c_1((\nu'')^* \mathbb{L}_{\mathfrak{M}_1}^{<0})
   +c_1((W^{\vee})^{<0}) +\delta_1(\nu'') \\
   &=\frac{1}{2}c_1((\lambda_0')^*r^*\mathbb{L}_{\mathfrak{N}}^{<0})-\frac{1}{2}c_1(\mathbb{L}_q|_{x}^{<0})+c_1(\Omega_{\mathbb{P}^n}|_{x}^{<0})-c_1(\mathbb{L}_{N}|_{x}^{<0})+\delta_1(\nu'') \\
   &=\frac{1}{2}c_1(\nu^* \mathbb{L}_{\mathfrak{M}_2}^{<0})+\frac{1}{2}c_1(\mathbb{L}_p|_{x}^{<0})-\frac{1}{2}c_1(\mathbb{L}_q|_{x}^{<0})+c_1(\Omega_{\mathbb{P}^n}|_{x}^{<0})+c_1((V^{\vee})^{<0})-c_1(\mathbb{L}_p|_{x}^{<0})+\delta_1(\nu'') \\
   &=\frac{1}{2}c_1(\nu^* \mathbb{L}_{\mathfrak{M}_2}^{<0})+c_1(\Omega_{\mathbb{P}^n}|_{x}^{<0})+c_1((V^{\vee})^{<0})-\frac{1}{2}c_1(\mathbb{L}_p|_{x}^{<0})-\frac{1}{2}c_1(\mathbb{L}_q|_{x}^{<0})+\delta_1(\nu'') \\
   &=\frac{1}{2}c_1(\nu^* \mathbb{L}_{\mathfrak{M}_2}^{<0})+c_1(\Omega_{\mathbb{P}^n}|_{x}^{<0})+c_1((V^{\vee})^{<0})+\delta_2(\nu).
\end{align*}
Here
we have used (\ref{c1zero}), (\ref{equal:delta}) for the last identity. 
Similarly, we have 
\begin{align*}
  &\frac{1}{2}c_1((\nu'')^* \mathbb{L}_{\mathfrak{M}_1}^{>0})
   +c_1((W^{\vee})^{>0}) +\delta_1(\nu'') \\
   &=\frac{1}{2}c_1(\nu^* \mathbb{L}_{\mathfrak{M}_2}^{>0})+c_1(\Omega_{\mathbb{P}^n}|_{x}^{>0})+c_1((V^{\vee})^{>0})+\delta_2(\nu).
\end{align*}
By combining (\ref{id:lambda}), (\ref{wt:iprime}), 
  it follows that 
  \begin{align}\label{wt:condf2}
     &\wt(\lambda^* i_{*}(q')^* A ) \subset \left[c_1(\Omega_{\mathbb{P}^n}|_{x}^{<0}), 
     c_1(\Omega_{\mathbb{P}^n}|_{x}^{>0}) 
     \right] \\
    \notag &+\left[\frac{1}{2}c_1(\nu^* \mathbb{L}_{\mathfrak{M}_2}^{<0})
   +c_1((V^{\vee})^{<0}), \frac{1}{2}c_1(\nu^* \mathbb{L}_{\mathfrak{M}_2}^{>0})
   +c_1((V^{\vee})^{>0})
   \right]+\delta_2(\nu).
  \end{align}
  Since we have 
  \begin{align*}
      \iota_{*}(\nu_V^{\mathrm{reg}})^*p_{*}q^* A \cong p''_{*}i_{*}(q')^*A
  \end{align*}
  by (\ref{isom:pqA}), (\ref{isom:pprime}), the condition (\ref{wt:condf})
  follows from (\ref{wt:condf2}) and applying Lemma~\ref{lem:pspace} below. 
\end{proof}
We have used the following lemma from~\cite{PTlim}:
\begin{lemma}\emph{(\cite[Lemma~5.11]{PTlim})}\label{lem:pspace}
Let $\mathbb{G}_m$ act on $\mathbb{P}^n$ linearly. 
For $\mathcal{E} \in \Coh(\mathbb{P}^n/\mathbb{G}_m)$, 
suppose that for any $\mathbb{G}_m$-fixed closed point 
$x \in \mathbb{P}^n$ we have 
\begin{align*}
    \wt(\mathcal{E}|_{x}) \subset \left[c_1 (\Omega_{\mathbb{P}^n}|_{x}^{<0}), 
    c_1 (\Omega_{\mathbb{P}^n}|_{x}^{>0})
    \right]+I
\end{align*}
for an interval $I\subset \mathbb{R}$. 
Then we have 
\begin{align*}
    \wt(\Gamma(\mathbb{P}^n, \mathcal{E}))
    \subset I. 
\end{align*}
\end{lemma}

We also have the following variant of Proposition~\ref{prop:functL}:
\begin{prop}\label{prop:functL2}
Suppose that there exist 
$\delta_i \in \mathrm{CL}(\mathfrak{M}_i)_{\mathbb{Q}}^{\vee}$
such that 
\begin{align}\label{pqdelta2}
    q^* \delta_1=p^* \delta_2-\frac{1}{2}\delta_P(\det \mathbb{L}_p) \in \mathrm{CL}(\mathfrak{N})_{\mathbb{Q}}^{\vee}.
\end{align}
Then the functor
$p_{*}(q^*(-)\otimes \omega_q)$ restricts to the functor 
\begin{align*}
 p_{*}(q^*(-)\otimes \omega_q) \colon   \LL(\mathfrak{M}_1)_{\delta_1} \to \LL(\mathfrak{M}_2)_{\delta_2}.
\end{align*}
\end{prop}
\begin{proof}
    By noting that 
    \begin{align*}\delta_P(\omega_q)=\delta_P(\det \mathbb{L}_q)
    =\delta_P(\det \mathbb{L}_p)
    \in\Pic(\mathfrak{N})_{\mathbb{Q}}^{\mathrm{CL}},
    \end{align*}
   from (\ref{c1zero}), the 
    result holds as in the proof of Proposition~\ref{prop:functL}. 
\end{proof}
\begin{cor}\label{cor:right}
In the setting of Proposition~\ref{prop:functL}, suppose furthermore that both \(p\) and \(q\) are quasi-smooth and projective.
    Then the functor (\ref{funct:coh}) admits a right adjoint
    \begin{align*}
        q_{*}p^! \colon \LL(\mathfrak{M}_2)_{\delta_2} \to \LL(\mathfrak{M}_1)_{\delta_1}.
    \end{align*}
    Here $p^!(-)=p^*(-)\otimes \omega_p[\dim p]$.
\end{cor}
\begin{proof}
    Since both $p$ and $q$ are quasi-smooth and projective, the functor 
    \begin{align*}
      q_{*}p^! \colon  \Coh(\mathfrak{M}_2) \stackrel{p^!}{\to} \Coh(\mathfrak{N}) \stackrel{q_{*}}{\to} \Coh(\mathfrak{M}_1)
    \end{align*}
    gives a right adjoint of (\ref{funct:coh}). By Proposition~\ref{prop:functL2}, 
    the above functor sends $\LL(\mathfrak{M}_2)_{\delta_2}$
    to $\LL(\mathfrak{M}_1)_{\delta_1}$, therefore the corollary holds. 
\end{proof}

We also have the following variant: 
\begin{prop}\label{prop:functL0.5}
Suppose that 
$q^* \delta_1=p^* \delta_2$ and 
$\det \mathbb{L}_p^{1/2}\in \mathrm{Pic}(\mathfrak{N})$ exists. Then 
we have the functor 
\begin{align*}
    p_{*}(q^*(-)\otimes \det \mathbb{L}_p^{1/2}) \colon 
    \LL(\mathfrak{M}_1)_{\delta_1} \to \LL(\mathfrak{M}_2)_{\delta_2}.
\end{align*}
\end{prop}
    \begin{proof}
        The same proof as for Proposition~\ref{prop:functL} applies. 
    \end{proof}

\subsection{Lagrangian correspondences}
As a corollary of the result in the previous subsection, we have the functoriality with respect to Lagrangian correspondences for zero-shifted symplectic stacks. 
Here we refer to~\cite{PTVV} for shifted symplectic stacks and 
their Lagrangians. 
This fact will not be used in this paper, but it is of independent interest. 

\begin{cor}\label{cor:Lag}
Let $\mathfrak{M}_1, \mathfrak{M}_2$ be quasi-smooth, QCA, zero-shifted symplectic stacks over $k$. For a Lagrangian
correspondence 
\begin{align*}
\mathfrak{M}_2 \stackrel{q}{\leftarrow} \mathfrak{L} \stackrel{p}{\to}
\mathfrak{M}_1
\end{align*}
suppose that $q$ is quasi-smooth, $p$ is projective and $\det \mathbb{L}_{\mathfrak{L}}^{1/2} \in \mathrm{Pic}(\mathfrak{L})$ exists. 
Let $\delta_i \in \mathrm{CL}(\mathfrak{M}_i)_{\mathbb{Q}}^{\vee}$ for $i=1, 2$ satisfy $p^* \delta_1=q^* \delta_2$. 
Then we have 
the induced functor 
\begin{align*}
    p_{*}(q^*(-)\otimes \det \mathbb{L}_{\mathfrak{L}}^{1/2}) \colon \LL(\mathfrak{M}_2)_{\delta_2} \to \LL(\mathfrak{M}_1)_{\delta_1}.
\end{align*}
\end{cor}
\begin{proof}
Since $(p,q) \colon \mathfrak{L} \to \mathfrak{M}_1\times \mathfrak{M}_2$
is Lagrangian, 
we have an equivalence 
    \begin{align*}
    \mathbb{T}_p \stackrel{\sim}{\to} \mathbb{L}_q[-1].
    \end{align*}
    Therefore the condition (\ref{c1zero}) is satisfied. 
    Also $\det \mathbb{L}_p=\det \mathbb{L}_{\mathfrak{L}}$
    as $\det \mathbb{L}_{\mathfrak{M}_1}=\mathcal{O}_{\mathfrak{M}_1}$. 
    The corollary follows from Proposition~\ref{prop:functL0.5}.
\end{proof}

By applying the above corollary to $\mathfrak{M}_2=\mathrm{pt}$, we have the following corollary, which implies that $\LL(\mathfrak{M})_{0}$ may be regarded 
as an algebro-geometric category of $A$-branes for zero-shifted symplectic stacks: 
\begin{cor}\label{cor:Abrane}
Let $\mathfrak{M}$ be a quasi-smooth, QCA, zero-shifted symplectic 
stack over $k$. Let $p \colon \mathfrak{L} \to \mathfrak{M}$ be a Lagrangian 
such that $p$ is projective, $\mathfrak{L}$ is quasi-smooth
and $\det \mathbb{L}_{\mathfrak{L}}^{1/2} \in \mathrm{Pic}(\mathfrak{L})$ exists. Then we have 
\begin{align*}
    p_{*}\det \mathbb{L}_{\mathfrak{L}}^{1/2} \in \LL(\mathfrak{M})_{0}.
\end{align*}
\end{cor}

\subsection{The case of closed immersions}
We return to the setting of Subsection~\ref{subsec:functors}. 
We furthermore assume that 
\begin{align*}p \colon \mathfrak{N} \to \mathfrak{M}_2
\end{align*}
is a closed immersion. In this subsection, we prove a kind of converse statement of Proposition~\ref{prop:functL}.
In this case, we have the following lemma: 
\begin{lemma}\label{lem:pU}
If $p$ is a closed immersion, then the derived scheme $N$ in the diagram (\ref{dia:N})
is either $\emptyset$, or equivalent to $U[-1]$
    for a finite-dimensional $\mathbb{G}_m$-representation $U$. 
\end{lemma}
\begin{proof}
We may assume that $N\neq \emptyset$. 
Since $p$ is a closed immersion, the map
$N \to V[-1]$ in the diagram (\ref{dia:N}) is a closed immersion, 
therefore $N^{\mathrm{cl}}=\Spec k$. Moreover $N$ is quasi-smooth 
as $N\to \mathfrak{N}$ is quasi-smooth. 
Since $N$ is quasi-smooth and $N^{\mathrm{cl}}=\Spec k$, it is equivalent to $U[-1]$ 
for a finite-dimensional $\mathbb{G}_m$-representation $U$. 
\end{proof} 

The diagram (\ref{dia:N}) is then written as 
 
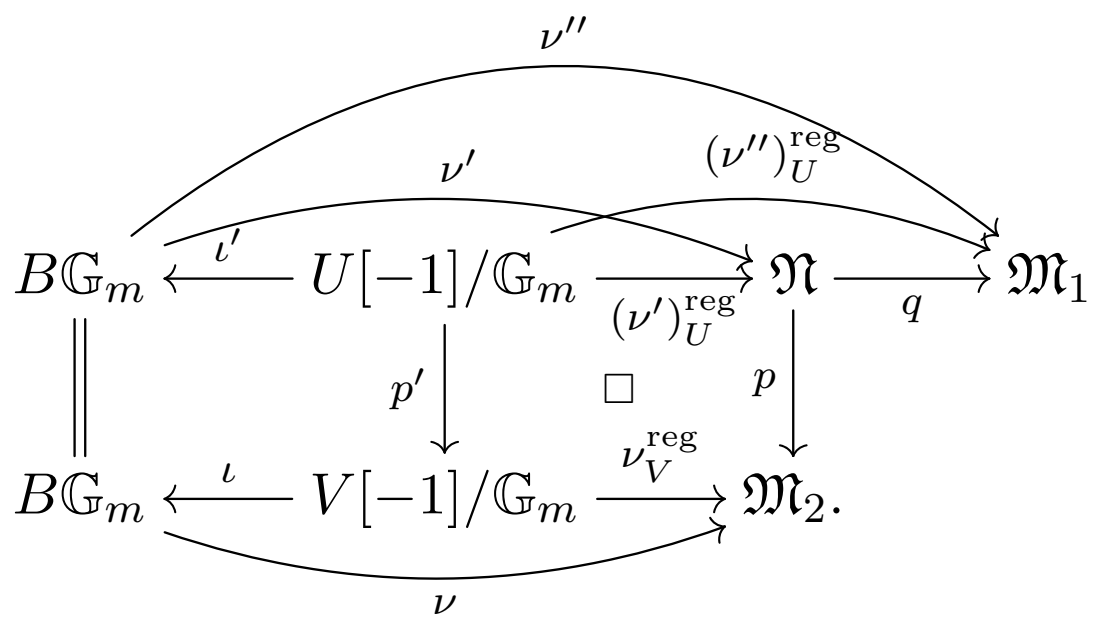
\begin{align}\label{dia:N2}
      \xymatrix{
\bgm \ar@{=}[d]  \ar@/^15pt/[rr]^-{\nu'}  \ar@/^40pt/[rrr]^-{\nu''}& U[-1]/\mathbb{G}_m  \ar@/^15pt/[rr]^-{(\nu'')_U^{\mathrm{reg}}} \ar[r]_-{(\nu')_U^{\mathrm{reg}}} \ar[d]_-{p'} \ar[l]_-{\iota'} \diasquare & \mathfrak{N} \ar[d]_-{p} 
\ar[r]_-{q} & \mathfrak{M}_1 \\
\bgm \ar@/_15pt/[rr]_-{\nu} & V[-1]/\mathbb{G}_m \ar[l]_-{\iota} \ar[r]^-{\nu_V^{\mathrm{reg}}} & \mathfrak{M}_2. 
      }
  \end{align}
Since $(\nu')_U^{\mathrm{reg}}$, $(\nu'')_U^{\mathrm{reg}}$ are quasi-smooth, 
they are regularization maps for $\nu'$, $\nu''$. 
Here we have replaced the notation $r, q'$ with  $(\nu')_U^{\mathrm{reg}}$, $(\nu'')_U^{\mathrm{reg}}$ in order to indicate 
that they are regularization maps. 
We will use the following lemma: 
\begin{lemma}\label{lem:nuineq}
Suppose that $\delta_1, \delta_2$ satisfy (\ref{pqdelta}). 
For an object $A\in \Coh(\mathfrak{M}_1)$ and an interval $I\subset \mathbb{R}$, we have 
\begin{align*}
\wt(\iota'_{*} ((\nu'')_U^{\mathrm{reg}})^* A) \subset 
&\left[\frac{1}{2}c_1((\nu'')^* \mathbb{L}_{\mathfrak{M}_1}^{<0})
   +c_1((U^{\vee})^{<0}), \frac{1}{2}c_1((\nu'')^* \mathbb{L}_{\mathfrak{M}_1}^{>0})
   +c_1((U^{\vee})^{>0})  \right] \\
   &+\delta_1(\nu'') +I
\end{align*}
if and only if 
    \begin{align*}
        \wt(\iota_{*}(\nu_V^{\mathrm{reg}})^*p_{*}q^*A)\subset & \left[ \frac{1}{2}c_1(\nu^* \mathbb{L}_{\mathfrak{M}_2}^{<0})+c_1((V^{\vee})^{<0}),  \frac{1}{2}c_1(\nu^* \mathbb{L}_{\mathfrak{M}_2}^{>0})+c_1((V^{\vee})^{>0}) \right] \\
        &+\delta_2(\nu)+I.
    \end{align*}
\end{lemma}
\begin{proof}
  We have 
    \begin{align}\label{wt:identity}
        \wt(\iota_{*}(\nu_V^{\mathrm{reg}})^* p_{*}q^*A)=\wt(\iota'_{*}((\nu'')_U^{\mathrm{reg}})^*A)
    \end{align}
    by (\ref{isom:pqA}).    
    As in the computations in Proposition~\ref{prop:functL}, we 
 have the equality 
      \begin{align}\label{c1identity}
&\frac{1}{2}c_1((\nu'')^* \mathbb{L}_{\mathfrak{M}_1}^{<0})
   +c_1((U^{\vee})^{<0}) +\delta_1(\nu'') \\
  \notag &=\frac{1}{2}c_1(\nu^* \mathbb{L}_{\mathfrak{M}_2}^{<0})+c_1((V^{\vee})^{<0})+\delta_2(\nu).
   \end{align}
   The analogous equality holds for the positive-weight parts.
   Therefore the lemma holds. 
\end{proof}

By applying the above lemma, we also have the following proposition: 
\begin{prop}\label{prop:closed2}
Let $\mathfrak{M}_1$ be a quasi-smooth, QCA and symmetric stack, 
and $\delta_1 \in \mathrm{CL}(\mathfrak{M}_1)_{\mathbb{Q}}^{\vee}$. 
For an object $A\in \Coh(\mathfrak{M}_1)$, suppose that the following holds: 
for any $\nu'' \colon \bgm \to \mathfrak{M}_1$, 
there is a diagram (\ref{dia:M12}) with closed immersion $p$
such that 
$\nu''$ lifts to a map $\nu' \colon \bgm \to \mathfrak{N}$, 
and that 
\begin{align*}
      p_{*}q^*A \in \LL(\mathfrak{M}_2)_{\delta_2}
\end{align*}
for $\delta_2 \in \mathrm{CL}(\mathfrak{M}_2)_{\mathbb{Q}}^{\vee}$ satisfying (\ref{pqdelta}). 
Then we have $A\in \LL(\mathfrak{M}_1)_{\delta_1}$.
\end{prop}
\begin{proof}
For maps $\nu''$ and $\nu'$ as above, consider the composition
\begin{align*}\nu\colon \bgm \stackrel{\nu'}{\to} \mathfrak{N} \stackrel{p}{\to} \mathfrak{M}_2.
\end{align*}
    Let $\nu_V^{\mathrm{reg}} \colon V[-1]/\mathbb{G}_m \to \mathfrak{M}_2$
    be a regularization map of $\nu$. The proposition follows by applying Lemma~\ref{lem:nuineq} 
    to the above map $\nu_V^{\mathrm{reg}}$. 
\end{proof}
\begin{remark}
In Proposition~\ref{prop:closed2}, if $\mathfrak{M}_1$ admits a good moduli space, it is enough 
to assume that the maps $\nu'' \colon \bgm \to \mathfrak{M}_1$ whose image
is a closed point in $\mathfrak{M}_1(k)$ satisfy the condition of Proposition~\ref{prop:closed2}. See Remark~\ref{rmk:closed-points}.
\end{remark}

\begin{remark}\label{rmk:locQCA}
The same statements 
of Propositions~\ref{prop:functL}, \ref{prop:functL2}, \ref{prop:closed2} 
hold for locally QCA stacks $\mathfrak{M}_i$ with 
$\LL(\mathfrak{M}_i)_{\delta_i}$ replaced by $\widetilde{\LL}(\mathfrak{M}_i)_{\delta_i}$, by applying the above propositions on QCA open substacks and passing to the limit.
\end{remark}

\subsection{Functors associated with a $\Theta$-stratum}\label{subsec:theta-functors}
Let $\mathfrak{M}$ be a quasi-smooth, QCA and symmetric stack with separated
diagonal, and let $\mathfrak{S}\hookrightarrow\mathfrak{M}$ be a $\Theta$-stratum
with center $\mathfrak{Z}$. We use the maps $p,q,\sigma$ of~(\ref{theta:str}).
The following is a quasi-smooth analogue of~\cite[Lemma~6.2]{PTlim}:
\begin{lemma}\label{lem:Zsym}
The center $\mathfrak{Z}$ is also symmetric. Moreover for each $x\in \mathfrak{S}(k)$, we have
\begin{align}\label{Lpq}
    [\mathbb{L}_p|_{x}]=-[\mathbb{L}_q^{\vee}|_{x}] \in K(BG_x).
\end{align}
\end{lemma}
\begin{proof}
By the computations of cotangent complexes in
\cite[Lemmas~1.3.2 and~1.5.5]{HalpK32}, we have 
\[
\mathbb{L}_{\mathfrak{S}}|_{\mathfrak{Z}}\cong\mathbb{L}_{\mathfrak{M}}|_{\mathfrak{Z}}^{\leq 0},
\quad
\mathbb{L}_{\mathfrak{Z}}\cong\mathbb{L}_{\mathfrak{M}}|_{\mathfrak{Z}}^{0}.
\]
Therefore $\mathfrak{Z}$ is symmetric from the symmetric condition of $\mathfrak{M}$, and 
\begin{align}\label{tri:NL}
\mathbb{L}_p|_{\mathfrak{Z}}
\cong \mathbb{L}_{\mathfrak{M}}|_{\mathfrak{Z}}^{>0}[1],
\quad 
\mathbb{L}_q|_{\mathfrak{Z}}
\cong \mathbb{L}_{\mathfrak{M}}|_{\mathfrak{Z}}^{<0}.
\end{align}
Therefore (\ref{Lpq}) holds if $x$ lies in the image of 
$\sigma \colon \mathfrak{Z}(k) \to \mathfrak{S}(k)$.

Let $y=q(x) \in \mathfrak{Z}(k)$, and let $\rho \colon G_x \to G_y$
be the induced map. 
Then for $F\in \mathrm{Perf}(\mathfrak{S})$, we have 
$[F|_x]=\rho^*[(\sigma^*F)|_y]$ in $K(BG_x)$. 
This follows from the semiorthogonal decomposition 
\begin{align*}
    \mathrm{Perf}(\mathfrak{S})=\langle \cdots, \mathrm{Perf}(\mathfrak{S})_{i-1}, 
    \mathrm{Perf}(\mathfrak{S})_i, \mathrm{Perf}(\mathfrak{S})_{i+1}, \ldots \rangle
\end{align*}
such that we have the equivalence $q^* \colon \mathrm{Perf}(\mathfrak{Z})_i \stackrel{\sim}{\to} \mathrm{Perf}(\mathfrak{S})_i$, 
see~\cite[Proposition~1.1.2(3), Lemmas~1.5.4 and~1.5.6]{HalpK32}. 
Therefore the lemma holds. 
\end{proof}

Using Propositions~\ref{prop:functL} and~\ref{prop:closed2}, we have the following proposition: 
\begin{prop}\label{prop:restL}
For $\delta \in \mathrm{CL}(\mathfrak{M})_{\mathbb{Q}}^{\vee}$, 
the functor (\ref{funct:up}) restricts to the functor 
\begin{align*}
    \Upsilon=p_{*}q^* \colon \LL(\mathfrak{Z})_{\delta|_{\mathfrak{Z}}+\delta_P(\det \mathbb{L}_p^{1/2}|_{\mathfrak{Z}})} \to \LL(\mathfrak{M})_{\delta}.
\end{align*}
Moreover for an object $\cE \in \Coh(\mathfrak{Z})$, we have:
\begin{align*}
\mbox{ if }\Upsilon(\cE) \in \LL(\mathfrak{M})_{\delta}, \mbox{ then } 
\cE \in \LL(\mathfrak{Z})_{\delta|_{\mathfrak{Z}}+\delta_P(\det \mathbb{L}_p^{1/2}|_{\mathfrak{Z}})}.
\end{align*}
\end{prop}
\begin{proof}
    By Lemma~\ref{lem:Zsym}, the condition of Proposition~\ref{prop:functL} is satisfied for $(\mathfrak{M}_1, \mathfrak{N}, \mathfrak{M}_2)=(\mathfrak{Z}, \mathfrak{S}, \mathfrak{M})$, 
    therefore the first claim follows. 
    The second claim follows from Proposition~\ref{prop:closed2}
    applied to $(\mathfrak{M}_1, \mathfrak{N}, \mathfrak{M}_2)=(\mathfrak{Z}, \mathfrak{S}, \mathfrak{M})$. Note that the assumption of
    Proposition~\ref{prop:closed2}
    is satisfied since $q \colon \mathfrak{S} \to \mathfrak{Z}$
    admits a section $\sigma$ in the diagram (\ref{theta:str}). 
\end{proof}

For $z\in\mathfrak{Z}(k)$, let $\nu\colon\bgm\to\mathfrak{M}$ be the
associated map, and choose a regularization $\nu_V^{\mathrm{reg}}$.
The weight identities in the proof of Lemma~\ref{lem:nuineq} give the following lemma, which will be used in the next section.
\begin{lemma}\label{lem:closed}
For $A\in \Coh(\mathfrak{Z})_w$ with $z \in \mathrm{Supp}(A)$, 
we have 
\begin{align*}
   & \wt^{\mathrm{max}}(\iota_{*}(\nu_V^{\mathrm{reg}})^* p_{*}q^* A)=w+\frac{1}{2}c_1(\nu^* \mathbb{L}_{\mathfrak{M}}^{>0})+c_1((V^{\vee})^{>0})-\frac{1}{2}c_1( \det \nu^*\mathbb{L}_p), \\
    &\wt^{\mathrm{min}}(\iota_{*}(\nu_V^{\mathrm{reg}})^* p_{*}q^* A)=w+\frac{1}{2}c_1(\nu^* \mathbb{L}_{\mathfrak{M}}^{<0})+c_1((V^{\vee})^{<0})-\frac{1}{2}c_1( \det \nu^*\mathbb{L}_p).
\end{align*}
    
\end{lemma}
\begin{proof}
We apply diagram~(\ref{dia:N2}) with
$(\mathfrak{M}_1,\mathfrak{N},\mathfrak{M}_2)=(\mathfrak{Z},\mathfrak{S},\mathfrak{M})$.
The canonical $\mathbb{G}_m$-action on $\mathbb{L}_{\mathfrak{Z}}$ has weight zero,
whereas $A$ has pure weight $w$ and $A|_{z}\neq 0$.
The calculations used in Lemma~\ref{lem:wtW} give 
\begin{align*}
\wt^{\mathrm{min}}\bigl(\iota'_*((\nu'')_U^{\mathrm{reg}})^*A\bigr)
 &=w+c_1((U^\vee)^{<0}),\\
\wt^{\mathrm{max}}\bigl(\iota'_*((\nu'')_U^{\mathrm{reg}})^*A\bigr)
 &=w+c_1((U^\vee)^{>0}).
\end{align*}
By Lemma~\ref{lem:Zsym}, we may apply the identities
(\ref{wt:identity}) and~(\ref{c1identity}) in the proof of Lemma~\ref{lem:nuineq},
with $\delta_2=0$ and
$\delta_1=\delta_P(\det\mathbb{L}_p^{1/2}|_{\mathfrak{Z}})$.
Since $(\nu'')^*\mathbb{L}_{\mathfrak{Z}}$ has only weight zero, these identities
and the analogue of~(\ref{c1identity}) for positive weights give the two asserted equalities.
\end{proof}

\section{Semiorthogonal decomposition of the limit category}\label{sec:decomposition}
In this section, we construct a semiorthogonal decomposition 
which is of different type from the one in Theorem~\ref{thm:BPT}. 
Instead of $\Coh(\mathfrak{M})$ for a stack $\mathfrak{M}$ with good moduli space, we construct a semiorthogonal 
decomposition of $\LL(\mathfrak{M})_{\delta}$ for a stack $\mathfrak{M}$, 
which does not necessarily admit a good moduli space but is equipped with a $\Theta$-stratification. 
Our point of view is that the semiorthogonal decomposition constructed in this section may be viewed as `Langlands dual' to the one in Theorem~\ref{thm:BPT}. 

We impose the following condition for stacks with $\Theta$-stratum:
\begin{assum}\label{assum:dag}
The derived stack $\mathfrak{M}$ over $k$ and its $\Theta$-stratum $\mathfrak{S} \hookrightarrow \mathfrak{M}$ satisfy the following conditions: 
\begin{enumerate}
\item 
The derived stack $\mathfrak{M}$ is quasi-smooth, QCA, has a separated diagonal, and is symmetric; 
\item The $\Theta$-stratum $\mathfrak{S} \hookrightarrow \mathfrak{M}$ is 
defined by numerical data $(\ell,b)$ as in
Definition~\ref{def:ndata};
\item 
For any $z\in \mathfrak{Z}(k)$ corresponding to 
$\nu \colon \bgm \to \mathfrak{M}$, 
the $\mathbb{G}_m$-weights of $\cH^{-1}(\nu^* \mathbb{L}_{\mathfrak{M}})$ are non-negative;
\item The center $\mathfrak{Z}$ of the $\Theta$-stratum $\mathfrak{S}$ admits 
a good moduli space $\mathfrak{Z} \to Z$.
\end{enumerate}
\end{assum}

\begin{remark}\label{rmk:connected}
In this section, until Theorem~\ref{thm:sodlim}, we further impose the condition that $\mathfrak{Z}$ is connected. 
This is just for simplifying the notation, and later in Theorem~\ref{thm:sodlim} we will summarize the 
result without the connectedness of $\mathfrak{Z}$. 
\end{remark}
\subsection{Semiorthogonal decompositions of limit categories}
The purpose of this subsection is to formulate Claim~\ref{claim:1}, 
which will be proved in the rest of this section. 

We consider a derived stack and its $\Theta$-stratum 
$(\mathfrak{M}, \mathfrak{S}, \mathfrak{Z})$
satisfying Assumption~\ref{assum:dag}. 
For $\delta \in \mathrm{CL}(\mathfrak{M})_{\mathbb{Q}}^{\vee}$, 
we take $m_{\delta}\in \mathbb{Q}$ to be the following central weight
\begin{align}\label{cent:m}
   m_{\delta}:= \wt^{\mathrm{cent}}(\delta_P(\det \mathbb{L}_p^{1/2}|_{\mathfrak{Z}})+\delta|_{\mathfrak{Z}}) \in \mathbb{Q}.
\end{align}
The above choice of $m_{\delta}$ is relevant to the limit category by 
the following lemma: 
\begin{lemma}\label{lem:choicem}
For $z \in \mathfrak{Z}(k)$, let $\nu \colon \bgm \to \mathfrak{M}$ be the corresponding map.
For the choice of $m_{\delta}$ as in (\ref{cent:m}),
we have 
\begin{align*}
m_{\delta}=\frac{1}{2}c_1(\nu^* \mathbb{L}_{\mathfrak{M}}^{<0})+\delta(\nu).
\end{align*}  
\end{lemma}
\begin{proof}
    The lemma follows from the fact 
    that the central weight of $\delta_P(\det \mathbb{L}_p^{1/2}|_{\mathfrak{Z}})$
    is equal to 
    \begin{align*}-\frac{1}{2}c_1(\nu^* \mathbb{L}_{\mathfrak{M}}^{>0})
    =\frac{1}{2}c_1(\nu^* \mathbb{L}_{\mathfrak{M}}^{<0}),
    \end{align*}
    which follows from (\ref{tri:NL})
    and the symmetric condition for $\mathfrak{M}$. 
\end{proof}

We have the following lemma: 
\begin{lemma}\label{lem:sodL1}
For $\cE \in \LL(\mathfrak{M})_{\delta}$, there is a distinguished triangle 
in $\Coh(\mathfrak{M})$
\begin{align}\label{sod:l1}
\cE_1 \to \cE \to \cE_2 
\end{align}
such that 
 $\cE_1 \in \mathcal{W}(m_{\delta})$ and 
$\cE_2=\Upsilon(A)$
for $A\in \Coh(\mathfrak{Z})_{m_{\delta}}$. 
In particular if $m_{\delta}\notin \mathbb{Z}$, we have 
$\cE_2=0$ and $\cE \in \mathcal{W}(m_{\delta})$. 
\end{lemma}
\begin{proof}
Recall the semiorthogonal decomposition (\ref{sod:dag}) with $m=m_{\delta}$.
Let 
\begin{align*}
    \Upsilon(A_{w^{\mathrm{min}}}), \ldots, \Upsilon(A_{\lfloor m_{\delta}\rfloor}), 
    P, \Upsilon(A_{\lfloor m_{\delta}\rfloor+1}), \ldots, \Upsilon(A_{w^{\mathrm{max}}})
\end{align*}
be the non-zero semiorthogonal factors of $\cE$, where 
$A_i \in \Coh(\mathfrak{Z})_i$ and $P\in \mathcal{W}(m_{\delta})$. 
It is enough to show that 
$\Upsilon(A_i)=0$ except $i=m_{\delta}$. 

Suppose that $\Upsilon(A_i)\neq 0$ for some $i>m_{\delta}$. 
We take $z\in \mathfrak{Z}(k)$ which lies in the support of $A_{w^{\mathrm{max}}}$
and the associated map $\nu \colon \bgm \to \mathfrak{M}$. 
Let $\nu_V^{\mathrm{reg}} \colon V[-1]/\mathbb{G}_m \to \mathfrak{M}$ be a regularization of $\nu$. 
Then by Lemma~\ref{lem:closed} and Lemma~\ref{lem:wtW},
we have 
\begin{align*}
    \wt^{\mathrm{max}}(\iota_{*}(\nu_V^{\mathrm{reg}})^* \cE)
    &=\wt^{\mathrm{max}}(\iota_{*}(\nu_V^{\mathrm{reg}})^* \Upsilon(A_{w^{\mathrm{max}}})) \\
    &=
    w^{\mathrm{max}}+\frac{1}{2}c_1(\nu^* \mathbb{L}_{\mathfrak{M}}^{>0})+c_1((V^{\vee})^{>0})-\frac{1}{2}c_1(\det \nu^*\mathbb{L}_p) \\
    &>\frac{1}{2}c_1(\nu^* \mathbb{L}_{\mathfrak{M}}^{>0})+c_1((V^{\vee})^{>0})+\delta(\nu).
\end{align*}
Here the last inequality follows from 
\begin{align*}
    w^{\mathrm{max}}>m_{\delta}= \delta(\nu)+\frac{1}{2}c_1(\det \nu^*\mathbb{L}_p).
\end{align*}
This contradicts $\cE \in \LL(\mathfrak{M})_{\delta}$, 
therefore $\Upsilon(A_i)=0$ for $i>m_{\delta}$. 
The same argument also shows that $\Upsilon(A_i)=0$ for $i<m_{\delta}$. 
\end{proof}

We also have the following lemma: 
\begin{lemma}\label{lem:E1}
In Lemma~\ref{lem:sodL1}, suppose that $\cE_1 \in \LL(\mathfrak{M})_{\delta}$ for all $\cE \in \LL(\mathfrak{M})_{\delta}$. 
Then we have 
\begin{align*}\cE_2 \in \Upsilon(\LL(\mathfrak{Z})_{\delta|_{\mathfrak{Z}}+\delta_P(\det \mathbb{L}_p^{1/2}|_{\mathfrak{Z}})})
\end{align*}
    and we have the semiorthogonal decomposition 
    \begin{align}\label{sod:l}
        \LL(\mathfrak{M})_{\delta}=\left\langle \Upsilon(\LL(\mathfrak{Z})_{\delta|_{\mathfrak{Z}}+\delta_P(\det \mathbb{L}_p^{1/2}|_{\mathfrak{Z}})}), \mathcal{W}(m_{\delta}) \cap \LL(\mathfrak{M})_{\delta} \right\rangle.
    \end{align}
\end{lemma}
\begin{proof}
Take $\cE \in \LL(\mathfrak{M})_{\delta}$
and the distinguished triangle (\ref{sod:l1}) in $\Coh(\mathfrak{M})$.
    If $\cE_1 \in \LL(\mathfrak{M})_{\delta}$, then 
    $\cE_2=\Upsilon(A) \in \LL(\mathfrak{M})_{\delta}$. 
    It follows that 
    \begin{align*}A \in \LL(\mathfrak{Z})_{\delta|_{\mathfrak{Z}}+\delta_P(\det \mathbb{L}_p^{1/2}|_{\mathfrak{Z}})}
    \end{align*}
    by Proposition~\ref{prop:restL}.
    Then the distinguished triangle (\ref{sod:l1}) implies that the 
    left-hand side of (\ref{sod:l}) is contained in the right-hand side. 

    On the other hand, the right-hand side of (\ref{sod:l}) is semiorthogonal by Theorem~\ref{thm:window}, and is 
contained in the left-hand side by Proposition~\ref{prop:restL}. 
Therefore the lemma holds. 
\end{proof}

Below we will consider the following claim:
\begin{claim}\label{claim:1}
For $\cE \in \LL(\mathfrak{M})_{\delta}$, take the 
distinguished triangle (\ref{sod:l1}). Then we 
    have $\cE_1 \in \LL(\mathfrak{M})_{\delta}$ and the restriction functor 
    gives an equivalence 
    \begin{align}\label{equiv:WK}
   j^* \colon \mathcal{W}(m_{\delta}) \cap \LL(\mathfrak{M})_{\delta} \stackrel{\sim}{\to}
    \LL(\mathfrak{M}^{\circ})_{\delta}. 
    \end{align}
    Here $j$ is the open immersion in (\ref{theta:str}). 
\end{claim}

By Lemma~\ref{lem:E1}, if Claim~\ref{claim:1} is true, then 
we have the semiorthogonal decomposition 
\begin{align*}
      \LL(\mathfrak{M})_{\delta}=\left\langle \Upsilon(\LL(\mathfrak{Z})_{\delta|_{\mathfrak{Z}}+\delta_P(\det \mathbb{L}_p^{1/2}|_{\mathfrak{Z}})}), j_{!}\LL(\mathfrak{M}^{\circ})_{\delta} \right\rangle
\end{align*}
where $j_{!}$ is the fully faithful left adjoint of 
$j^* \colon \LL(\mathfrak{M})_{\delta} \to \LL(\mathfrak{M}^{\circ})_{\delta}$. 

\subsection{The local model}
In this subsection, we prove Claim~\ref{claim:1} for a local model. The results of this subsection 
are essentially proved in~\cite{PTlim}, but here we prove them in a more general setting. We will later reduce the result in a global setting to 
the local model. 

Let $G$ be a reductive group and $Y$ a finite-dimensional $G$-representation. 
We consider the quotient stack 
\begin{align*}\mathcal{Y}=Y/G. 
\end{align*}
Note that $\Y$ is symmetric if $Y$ is a symmetric $G$-representation, i.e. $Y \cong Y^{\vee}$ as 
a $G$-representation. 

For a rational character $\chi\in\mathrm{Pic}(BG)_{\mathbb Q}$,
write $\ell=\delta_P(\chi)$ and let
$\mathcal{Y}^{\ell\text{-ss}}\subset\mathcal{Y}$ be its semistable locus. By fixing a Weyl-invariant norm 
on $N(T)_{\mathbb{Q}}$, where $N(T)$ is the cocharacter lattice of
a maximal torus $T$ of $G$,
there is an associated \textit{Kempf--Ness
stratification}, see~\cite[Section~2.1]{halp}
\begin{align}\label{strata:Y}
    \mathcal{Y}=\mathcal{S}_1 \sqcup \cdots \sqcup \mathcal{S}_k \sqcup \mathcal{Y}^{\ell\text{-ss}}
\end{align}
with centers $\mathcal{S}_i \to \mathcal{Z}_i$
and associated one-parameter subgroups $\lambda_i \colon \mathbb{G}_m \to 
G$. 
Then 
\begin{align*}
    \mathcal{S}_i \subset \mathcal{Y}_i:=\mathcal{Y}\setminus (\cS_1 \sqcup \cdots \sqcup \mathcal{S}_{i-1})
\end{align*}
is a $\Theta$-stratum. 

Note that Assumption~\ref{assum:dag} is automatic for $\mathcal{S}_i \subset \mathcal{Y}_i$ if $Y$ is a symmetric $G$-representation. The $\Theta$-stratification (\ref{strata:Y})
is defined by the slope $\mu$
determined by $(\ell,b)$, where $\ell=\delta_P(\chi)$ and $b$ is
the quadratic form on the component lattice induced by the chosen
Weyl-invariant quadratic form on $N(T)_{\mathbb Q}$.
The following result is proved in~\cite{PTlim}:

\begin{thm}\emph{(\cite[Theorem~6.13]{PTlim})}\label{thm:sym1}
Suppose that $Y$ is a symmetric $G$-representation. Then for 
the symmetric stack $\mathcal{Y}_i$ and a $\Theta$-stratum $\mathcal{S}_i \subset \mathcal{Y}_i$ as above, the 
statement of Claim~\ref{claim:1} holds. 
\end{thm}

We extend the above result in the following way. 
Let $U$ be an affine scheme with a $G$-action and a $G$-equivariant \'etale map 
$a \colon U \to Y$. It is called \textit{strongly \'etale} if the following 
diagram is Cartesian 
\begin{align*}
    \xymatrix{
    \mathcal{U}:=U/G \ar[r]^-{a} \ar[d] \diasquare & \Y=Y/G \ar[d] \\
    U\ssslash G \ar[r] & Y\ssslash G.
    }
\end{align*}
By pulling back the stratification (\ref{strata:Y}), we obtain 
the $\Theta$-stratification 
\begin{align*}
    \mathcal{U}=\mathcal{S}_1' \sqcup \cdots \sqcup \mathcal{S}_k' \sqcup \mathcal{U}^{\ell \text{-ss}}
\end{align*}
with respect to the pull-back of $(\ell, b)$. 
Each $\Theta$-stratum is 
\begin{align*}
    \mathcal{S}_i' \hookrightarrow \mathcal{U}_i:=a^{-1}(\mathcal{Y}_i).
\end{align*}
Note that each stack $\mathcal{U}_i$ is symmetric if $Y$ is a symmetric 
$G$-representation. 
\begin{lemma}\label{lem:setale}
If $Y$ is symmetric and $\delta \in \mathrm{CL}(\mathcal{U})_{\mathbb{Q}}^{\vee}$ is a pull-back from $\mathrm{CL}(\mathcal{Y})_{\mathbb{Q}}^{\vee}$, then 
    each $\Theta$-stratum $\mathcal{S}_i' \hookrightarrow \mathcal{U}_i$ satisfies the statement of Claim~\ref{claim:1}. 
\end{lemma}
\begin{proof}
    We show that the restriction functor 
    \begin{align}\label{rest:LU}
        \LL(\mathcal{U})_{\delta} \to \LL(\mathcal{U}^{\ell\text{-ss}})_{\delta}
    \end{align}
    is essentially surjective. Indeed this is the only input in the proof of~\cite[Theorem~6.13]{PTlim} from the assumption that the stack is a linear quotient; once this is proven, the rest of the arguments for $\mathcal{S}_i' \hookrightarrow \mathcal{U}_i$ are the same as those for $\mathcal{S}_i \hookrightarrow \mathcal{Y}_i$. 

By repeatedly using Theorem~\ref{thm:sym1}, 
there is a window subcategory 
\begin{align*}
    \mathcal{W}_{\mathcal{Y}}^{\ell\text{-ss}}
    \subset \Coh(\Y), \  \mathcal{W}_{\mathcal{Y}}^{\ell\text{-ss}}
    \stackrel{\sim}{\to} \Coh(\Y^{\ell\text{-ss}})
\end{align*}
such that the above equivalence restricts to the equivalence 
\begin{align}\label{equiv:WL2}
    \mathcal{W}_{\Y}^{\ell\text{-ss}} \cap \LL(\Y)_{\delta} \stackrel{\sim}{\to} \LL(\Y^{\ell\text{-ss}})_{\delta}. 
\end{align}

Now we have 
\begin{align}\label{thick}
    \mathrm{thick}\langle a^* \LL(\Y)_{\delta}\rangle
    =\LL(\mathcal{U})_{\delta}, \ 
       \mathrm{thick}\langle a^* \LL(\Y^{\ell\text{-ss}})_{\delta}\rangle
    =\LL(\mathcal{U}^{\ell\text{-ss}})_{\delta}.
\end{align}
Here $\mathrm{thick}\langle-\rangle$ is the subcategory split-generated by $(-)$, the first identity is mentioned in~\cite[Section~3.4.10]{BPT}
and the second one also follows similarly by noting that 
(cf. the proof of~\cite[Lemma~3.4.4]{BPT})
\begin{align*}\mathrm{thick}\langle a^* \Coh(\Y^{\ell\text{-ss}})\rangle=\Coh(\mathcal{U}^{\ell\text{-ss}})
\end{align*}
and the semiorthogonal decompositions of Theorem~\ref{thm:BPT} for $\Y^{\ell\text{-ss}}$, $\mathcal{U}^{\ell\text{-ss}}$ are compatible with $a^*$. 

By Theorem~\ref{thm:window}, there is also a window subcategory 
\begin{align*}
    \mathcal{W}_{\mathcal{U}}^{\ell\text{-ss}} \subset \Coh(\mathcal{U}), \  \mathcal{W}_{\mathcal{U}}^{\ell\text{-ss}}\stackrel{\sim}{\to} \Coh(\mathcal{U}^{\ell\text{-ss}})
\end{align*}
such that $a^* \mathcal{W}_{\mathcal{Y}}^{\ell\text{-ss}} \subset \mathcal{W}_{\mathcal{U}}^{\ell\text{-ss}}$, by choosing the same $m$ in Theorem~\ref{thm:window} at each $\Theta$-stratum. 
Then we have the fully faithful functor 
\begin{align*}
    \mathcal{W}_{\mathcal{U}}^{\ell\text{-ss}} \cap \LL(\mathcal{U})_{\delta} \hookrightarrow \LL(\mathcal{U}^{\ell\text{-ss}})_{\delta}
\end{align*}
which is indeed essentially surjective (hence an equivalence) 
by (\ref{thick}) and the equivalence (\ref{equiv:WL2}).
In particular, the functor (\ref{rest:LU}) is essentially surjective. 
\end{proof}

Let $Y$ be a (not necessarily symmetric) finite-dimensional 
$G$-representation, and $V$ be another $G$-representation. 
Suppose that the following equality holds: 
\begin{align}\label{UVdual}
    [Y]-[V]=[Y^{\vee}]-[V^{\vee}] \in K(BG).
\end{align}
Let $U$ be an affine scheme with $G$-action with a strongly 
\'etale map $a \colon U \to Y$, and a $G$-equivariant map 
$s \colon U \to V$. 
Then by (\ref{UVdual}), the quasi-smooth derived stack 
\begin{align*}
    \mathfrak{U}=s^{-1}(0)/G \hookrightarrow \mathcal{U}=U/G
\end{align*}
is a symmetric stack, where $s^{-1}(0)$ is the derived zero locus. For a cocharacter $\lambda \colon \mathbb{G}_m \to G$, 
we have the induced maps
\begin{align*}
    s^{\lambda \geq 0} \colon U^{\lambda \geq 0} \to V^{\lambda \geq 0}, \ 
     s^{\lambda} \colon U^{\lambda} \to V^{\lambda}.
\end{align*}
Consider the derived stacks 
\begin{align*}
    \mathfrak{U}^{\lambda \geq 0}:=(s^{\lambda \geq 0})^{-1}(0)/G^{\lambda \geq 0}, \ 
     \mathfrak{U}^{\lambda}:=(s^{\lambda})^{-1}(0)/G^{\lambda}
\end{align*}
and the induced diagram 
\begin{align}\label{diagram:fU}
    \mathfrak{U}^{\lambda} \leftarrow \mathfrak{U}^{\lambda \geq 0} \to \mathfrak{U}.
\end{align}
Suppose that 
\begin{align*}\cS_1=\mathcal{Y}^{\lambda \geq 0} \hookrightarrow \mathcal{Y}
\end{align*}
is a closed KN stratum. 
Then $\mathfrak{U}^{\lambda \geq 0} \hookrightarrow \mathfrak{U}$ is a closed immersion and it is a 
$\Theta$-stratum. Below we take $\delta \in \mathrm{CL}(\mathcal{Y})_{\mathbb{Q}}^{\vee}$, and 
use the same notation $\delta$ for its pull-back to $\mathfrak{U}$. 
\begin{prop}\label{prop:sym2}
For the symmetric stack $\mathfrak{U}$ and its $\Theta$-stratum 
$\mathfrak{U}^{\lambda \geq 0} \hookrightarrow \mathfrak{U}$ as above, 
suppose that the $\lambda$-weights of $V$ are non-positive. 
Then 
the statement of Claim~\ref{claim:1} holds for the $\Theta$-stratum 
$\mathfrak{U}^{\lambda \geq 0} \hookrightarrow \mathfrak{U}$. 
\end{prop}
\begin{proof}
The proof is essentially a generalization of the arguments  in~\cite[Proposition~8.23]{PTlim}. 

  Let $\mathcal{X}$ be the following smooth stack 
  \begin{align}\label{stack:X2}
      \mathcal{X}=(U\times  V^{\vee})/G.
  \end{align}
  The map $s$ defines the following function $f$ on $\X$
  \begin{align*}
      f\colon \mathcal{X} \to \mathbb{A}^1, \ 
      (u, v) \mapsto \langle s(u), v\rangle.
  \end{align*}
By (\ref{UVdual}), the $G$-representation $Y\oplus V^{\vee}$ is self-dual, which implies that $\mathcal{X}$ is a symmetric stack. 
There is a Koszul equivalence~\cite{I,MR2982435,Hirano,T} 
\begin{align}\label{Koszul}
\mathsf{K} \colon 
\Coh(\mathfrak{U}) \stackrel{\sim}{\to} \mathrm{MF}^{\mathbb{G}_m}(\mathcal{X}, f).
\end{align}
Here the right-hand side is the dg-category of $\mathbb{G}_m$-equivariant 
matrix factorizations of the function $f$, where $\mathbb{G}_m$ acts on the second factor $V^{\vee}$ of $\X$ by weight two; for details 
of matrix factorizations and Koszul equivalence, we also refer to~\cite[Section~4.2.1]{BPT}. 

The equivalence (\ref{Koszul}) restricts to an equivalence, see~\cite[Proposition~3.14]{PTlim}, \cite[Lemma~4.2.2]{BPT}
\begin{align}\label{equiv:MF}
    \mathsf{K} \colon \LL(\mathfrak{U})_{\delta} \stackrel{\sim}{\to}
    \mathrm{MF}^{\mathbb{G}_m}(\LL(\mathcal{X})_{\delta'}, f).
\end{align}
Here we have set (omitting the notation for pull-backs to $\X$)
\begin{align*}
    \delta'=\delta+\delta_P(\omega_{\mathcal{Y}}^{1/2}) \in \mathrm{CL}(\X)_{\mathbb{Q}}^{\vee}. 
\end{align*}
By applying Lemma~\ref{lem:setale} to $\X^{\lambda \geq 0} \hookrightarrow \X$, we obtain a semiorthogonal decomposition 
\begin{align}\label{sod:Ldelta'}
    \LL(\mathcal{X})_{\delta'}=\left\langle \Upsilon_{\mathcal{X}}(\LL(\mathcal{X}^{\lambda})_{\delta''}), \mathcal{W}_{\mathcal{X}}(m_{\delta'})\cap \LL(\mathcal{X})_{\delta'} \right\rangle
\end{align}
with an equivalence 
\begin{align}\label{equiv:WL}
\mathcal{W}_{\X}(m_{\delta'}) \cap \LL(\X)_{\delta'} \stackrel{\sim}{\to}
\LL(\X^{\circ})_{\delta'}.
\end{align}

Here we explain the notation in the semiorthogonal 
decomposition (\ref{sod:Ldelta'}). We have the diagram 
\begin{align*}
    \mathcal{X}^{\lambda} \stackrel{q_{\X}}{\leftarrow} \mathcal{X}^{\lambda\geq 0} \stackrel{p_{\mathcal{X}}}{\hookrightarrow} \mathcal{X}.
\end{align*}
Note that 
the closed substack $\X^{\lambda \geq 0} \hookrightarrow \X$
is a $\Theta$-stratum with respect to $(\ell, b)$ pulled back to $\X$. 
We set 
\begin{align*}
\delta''=(\delta'+\delta_P(\det \mathbb{L}_{p_{\mathcal{X}}}^{1/2}))|_{\mathcal{X}^{\lambda}} \in \mathrm{CL}(\X^{\lambda})_{\mathbb{Q}}^{\vee},
\end{align*}
and the functor $\Upsilon_{\mathcal X}$ is given by 
\begin{align*}
    \Upsilon_{\mathcal{X}}=p_{\mathcal{X}*}q_{\mathcal{X}}^* \colon \LL(\mathcal{X}^{\lambda})_{\delta''} \to \LL(\mathcal{X})_{\delta'}.
\end{align*}
The subcategory 
$\mathcal{W}_{\mathcal{X}}(m_{\delta'}) \subset \Coh(\mathcal{X})$
is the window-subcategory as in Theorem~\ref{thm:window} for $\mathcal{X}$,
where $m_{\delta'}\in \mathbb{Q}$ is the central weight of $\delta''$.
We have the equivalence 
\begin{align*}
    \mathcal{W}_{\mathcal{X}}(m_{\delta'}) \stackrel{\sim}{\to}\Coh(\mathcal{X}^{\circ}),
\end{align*}
where $\mathcal{X}^{\circ}:=\mathcal{X} \setminus \mathcal{X}^{\lambda \geq 0}$. 

By applying $\mathrm{MF}^{\mathbb{G}_m}(-, f)$ to the semiorthogonal 
decomposition (\ref{sod:Ldelta'}), see~\cite[Proposition~2.5]{PT0}, 
we have the semiorthogonal decomposition 
\begin{align}\label{sod:applyMF}
 &\mathrm{MF}^{\mathbb{G}_m}(\LL(\mathcal{X})_{\delta'}, f)= \\
 \notag &\left\langle \Upsilon_{\mathcal{X}} \mathrm{MF}^{\mathbb{G}_m}(\LL(\mathcal{X}^{\lambda})_{\delta''}, f), \mathrm{MF}^{\mathbb{G}_m}(\mathcal{W}_{\mathcal{X}}(m_{\delta'}) \cap \LL(\mathcal{X})_{\delta'}, f) \right\rangle.
\end{align}

For an object $\cE\in \LL(\mathfrak{U})_{\delta}$, consider the distinguished triangle (\ref{sod:l1}). 
By the compatibility of the functor $\Upsilon_{\mathcal{X}}$ with the Koszul equivalence in~\cite[Section~2.4]{T}, the distinguished triangle 
\begin{align*}
    \mathsf{K}(\cE_1) \to \mathsf{K}(\cE) \to \mathsf{K}(\cE_2)
\end{align*}
corresponds to the decomposition with respect to the semiorthogonal 
decomposition (\ref{sod:applyMF}). 
It follows that 
\begin{align*}
    \mathsf{K}(\cE_1) \in \mathrm{MF}^{\mathbb{G}_m}(\mathcal{W}_{\mathcal{X}}(m_{\delta'}) \cap \LL(\mathcal{X})_{\delta'}, f)
    \subset \mathrm{MF}^{\mathbb{G}_m}(\LL(\mathcal{X})_{\delta'}, f).
\end{align*}
By (\ref{equiv:MF}), we conclude that $\cE_1 \in \LL(\mathfrak{U})_{\delta}$. Therefore the first statement of Claim~\ref{claim:1} holds. In particular 
we have the semiorthogonal decomposition 
\begin{align*}
    \LL(\mathfrak{U})_{\delta}=\left\langle \Upsilon \LL(\mathfrak{U}^{\lambda})_{\delta+\delta_P(\det \mathbb{L}_p^{1/2})}, \mathcal{W}(m_{\delta})\cap \LL(\mathfrak{U})_{\delta} \right\rangle,
\end{align*}
which corresponds to (\ref{sod:applyMF}) under the Koszul equivalence $\mathsf{K}$. 

Let $\pi_{\mathcal{U}} \colon \mathcal{X} \to \mathcal{U}$ be the projection. 
Then the assumption that $V$ has non-positive $\lambda$-weights implies 
that 
\begin{align*}
    \X^{\lambda \geq 0}=(U\times V^{\vee})^{\lambda \geq 0}/G^{\lambda \geq 0}=(U^{\lambda \geq 0} \times V^{\vee})/G^{\lambda \geq 0}=\pi_{\mathcal{U}}^{-1}(\mathcal{U}^{\lambda \geq 0}).
\end{align*}
Therefore by setting $\mathcal{U}^{\circ}:=\mathcal{U}\setminus \mathcal{U}^{\lambda \geq 0}$, we have 
$\mathcal{X}^{\circ}=\pi_{\mathcal{U}}^{-1}(\mathcal{U}^{\circ})$.
We have the commutative diagram 
\begin{align*}
    \xymatrix{
\mathcal{W}(m_{\delta}) \cap \LL(\mathfrak{U})_{\delta}  \ar[r]^-{\sim}_-{\mathsf{K}} \ar[dd] & \mathrm{MF}^{\mathbb{G}_m}(\mathcal{W}_{\mathcal X}(m_{\delta'}) \cap \LL(\mathcal{X})_{\delta'}, f)\ar[d]_-{\sim} \\
& \mathrm{MF}^{\mathbb{G}_m}(\LL(\mathcal{X}^{\circ})_{\delta'}, f)
\ar@{=}[d] \\
\LL(\mathfrak{U}^{\circ})_{\delta} \ar[r]^-{\sim}_-{\mathsf{K}^{\circ}} & 
    \mathrm{MF}^{\mathbb{G}_m}(\LL(\pi_{\mathcal{U}}^{-1}(\mathcal{U}^{\circ}))_{\delta'}, f).
    }
\end{align*}
The right top equivalence follows from (\ref{equiv:WL}). 
The bottom equivalence follows as in (\ref{equiv:MF}), see~\cite[Proposition~3.14]{PTlim},~\cite[Lemma~4.2.2]{BPT}. The left vertical arrow is the pull-back with respect to the open immersion 
$\mathfrak{U}^{\circ} \hookrightarrow \mathfrak{U}$. 
By the above diagram, the left vertical arrow is an equivalence, and the second statement of Claim~\ref{claim:1} holds. 
\end{proof}

\subsection{Construction of local model}
We will reduce Claim~\ref{claim:1} in a global 
setting to the case of the local model in the previous subsection. 
We will use the following proposition: 
\begin{prop}\label{prop:loc}
Let $\mathfrak{M}$ be a derived stack as in Assumption~\ref{assum:dag}.
For each point $x\in \mathfrak{M}(k)$, let 
$G_x'\subset G_x=\mathrm{Aut}(x)$ be a reductive subgroup. 
Then 
there is a quasi-smooth derived affine scheme $\Spec A$ over $k$ with 
a $G_x'$-action and a smooth map 
$\Spec A/G_x' \to \mathfrak{M}$
such that we have the following Cartesian square 
\begin{align}\label{Cart:A}
\xymatrix{
    \Spec A/G_x' \ar[r] \diasquare & \mathfrak{M} \\
    BG_x' \ar[u] \ar[r] & BG_x. \ar[u]
    }
\end{align}

Moreover the image $o \in \Spec A/G_x'$
of the left vertical arrow is the unique closed point in the fiber 
of $\Spec A/G_x' \to \Spec A \ssslash G_x'$
at $\overline{o}\in \Spec A \ssslash G_x'$.
\end{prop}
\begin{proof}
    The result for the classical stack is due to~\cite[Theorems~1.1 and~8.1]{AHR2}.
  It can be extended to quasi-smooth stacks as in~\cite[Theorem~2.1.1]{HalpK32}.
\end{proof}
\begin{remark}\label{rmk:isomSpec}
    In Proposition~\ref{prop:loc},
    we have the isomorphism 
    \begin{align}\label{isom:tangent}
        \cH^0(\mathbb{T}_{\Spec A}|_{o}) \stackrel{\cong}{\to} \cH^0(\mathbb{T}_{\mathfrak{M}}|_{x}).
    \end{align}
Indeed for the map 
\begin{align*}
f\colon \Spec A \to \Spec A/G_x' \to \mathfrak{M},
\end{align*}
the above isomorphism (\ref{isom:tangent}) follows from the exact sequence 
\begin{align*}
0 \to \cH^{-1}(\mathbb{T}_{\mathfrak{M}}|_{x}) \to \cH^0(\mathbb{T}_f|_{o}) \to    \cH^0(\mathbb{T}_{\Spec A}|_{o}) \to \cH^0(\mathbb{T}_{\mathfrak{M}}|_{x}) \to 
 0
\end{align*}
and the map $\cH^{-1}(\mathbb{T}_{\mathfrak{M}}|_{x}) \to \cH^0(\mathbb{T}_f|_{o})$
is an isomorphism, as both are naturally isomorphic to the Lie algebra of $G_x$ by the Cartesian square (\ref{Cart:A}).
\end{remark}

Let $(\mathfrak{M}, \mathfrak{S}, \mathfrak{Z})$ be a derived stack with a $\Theta$-stratum satisfying Assumption~\ref{assum:dag}. 
For a closed point $z\in \mathfrak{Z}(k)$, its 
automorphism group $\mathrm{Aut}_{\mathfrak{Z}}(z)$ is a reductive 
group since $\mathfrak{Z}$ admits a good moduli space. 
Suppose that $z$ corresponds to $x\in \mathfrak{M}(k)$ with homomorphism \begin{align*}\lambda \colon \mathbb{G}_m \to G_x:=\mathrm{Aut}_{\mathfrak{M}}(x).
\end{align*}
Then because $\mathrm{Aut}_{\mathfrak{S}}(x)=\mathrm{Aut}_{\mathfrak{M}}(x)$, we have 
\begin{align*}
    \mathrm{Aut}_{\mathfrak{Z}}(z)=G_x^{\lambda} \subset G_x=G_x^{\lambda \geq 0}.
\end{align*}
Here $G_x^{\lambda}$ is 
    the subgroup which commutes with $\lambda$. 
It follows that 
\begin{align*}\mathfrak{n}_x:=\mathfrak{g}_x/\mathfrak{g}_x^{\lambda}
\end{align*}
    has positive $\lambda$-weights. 
    Here $\mathfrak{g}_x$ and $\mathfrak{g}_x^{\lambda}$ are the Lie algebras of $G_x$, $G_x^{\lambda}$, respectively. 
    Moreover $G_x^{\lambda}$ is a maximal reductive subgroup of $G_x$, 
    and we have the bijection 
    \begin{align}\label{closed:bij}
    \mathfrak{Z}(k)^{\mathrm{closed}} \stackrel{\cong}{\to} Z(k).
    \end{align}
    Here the left-hand side is the set of isomorphism classes of closed $k$-points of $\mathfrak{Z}$. 
 
We set 
\begin{align*}
    Y_x:=\cH^0(\mathbb{T}_{\mathfrak{M}}|_{x}), \ 
    V_x :=\cH^1(\mathbb{T}_{\mathfrak{M}}|_{x}).
\end{align*}
The condition (iii) in 
Assumption~\ref{assum:dag} implies that the $\lambda$-weights of $V_x$ are non-positive.

\begin{prop}\label{prop:localmodel}
There is a smooth affine scheme $U_x$ with 
$G_x^{\lambda}$-action, a strongly \'etale $G_x^{\lambda}$-equivariant map $U_x \to Y_x$
and a $G_x^{\lambda}$-equivariant map $\kappa \colon U_x \to V_x$
with a smooth map 
\begin{align*}
 g_x\colon   \mathfrak{U}_x:=\kappa^{-1}(0)/G_x^{\lambda} \to \mathfrak{M}
\end{align*}
such that we have the Cartesian square 
\begin{align*}
\xymatrix{
\mathfrak{U}_x \ar[r]^-{g_x} \diasquare & \mathfrak{M} \\
    BG_x^{\lambda} \ar[u] \ar[r] & BG_x. \ar[u]
    }
\end{align*}
\end{prop}
\begin{proof}
We apply Proposition~\ref{prop:loc} for $G_x'=G_x^{\lambda}$, 
and consider the quasi-smooth affine derived scheme $\Spec A$
with a smooth map 
\begin{align}\label{map:g}
    g\colon \Spec A/G_x^{\lambda} \to \mathfrak{M}
\end{align}
as in Proposition~\ref{prop:loc}. Since $\Spec A$ is quasi-smooth and $G_x^{\lambda}$ is reductive, by~\cite[Proposition~7.3.1]{T}, \cite[Lemma~4.2.3]{BPT}, there is a smooth affine 
scheme $U$ with a $G_x^{\lambda}$-action, 
$G_x^{\lambda}$-equivariant closed immersion $\Spec A \hookrightarrow U$, 
a $G_x^{\lambda}$-equivariant map $s \colon U \to V$ for a finite-dimensional $G_x^{\lambda}$-representation $V$
such that 
\begin{align}\label{Vx}
\Spec A/G_x^{\lambda} \simeq s^{-1}(0)/G_x^{\lambda}.
\end{align}

By shrinking $\Spec A$ if necessary, we can take the closed embedding 
$\Spec A \hookrightarrow U$ such that 
$\cH^0(\mathbb{T}_{\Spec A}|_{o}) \cong T_{U}|_{o}$.
Then by Remark~\ref{rmk:isomSpec}, $T_{U}|_{o}$ is isomorphic to $Y_x=\cH^0(\mathbb{T}_{\mathfrak{M}}|_{x})$ as a $G_x^{\lambda}$-representation, 
and 
$ds|_{o}=0$. It follows that $V=V_x=\cH^1(\mathbb{T}_{\mathfrak{M}}|_{x})$
and 
we can assume that there is a strongly \'etale $G_x^{\lambda}$-equivariant map $U \to Y_x$, 
see~\cite[Lemma~3.4.6]{BPT}.
By setting $U_x=U$, we obtain the proposition. 
\end{proof}

The derived stack $\mathfrak{U}_x$ is not necessarily symmetric, and we modify it to a symmetric stack. 
We define 
\begin{align*}
    \fU_x^{\dag}:=(\kappa^{-1}(0)\times \mathfrak{n}_x^{\vee})/G_x^{\lambda}.
\end{align*}

\begin{lemma}\label{lem:symU}
    The derived stack $\fU_x^{\dag}$ is symmetric. 
\end{lemma}
\begin{proof}
    Since $\mathfrak{M}$ is symmetric, the element 
    \begin{align*}
        [Y_x]-[V_x]-[\mathfrak{g}_x] \in K(BG_x)
    \end{align*}
    is self-dual. By restricting to $G_x^{\lambda}$, it is also self-dual 
    in $K(BG_x^{\lambda})$. 
    Since $\mathfrak{g}_x=\mathfrak{g}_x^{\lambda}\oplus \mathfrak{n}_x$ in $K(BG_x^{\lambda})$, 
    the element 
    \begin{align*}
        [Y_x]-[V_x]-[\mathfrak{n}_x]-[\mathfrak{g}_x^{\lambda}] \in K(BG_x^{\lambda})
    \end{align*}
    is self-dual. It follows that the element 
    \begin{align*}
        [Y_x]+[\mathfrak{n}_x^{\vee}]-[V_x]-[\mathfrak{g}_x^{\lambda}] \in K(BG_x^{\lambda})
    \end{align*}
    is self-dual. This implies that $\fU_x^{\dag}$ is symmetric, as 
    its tangent complex is given by 
    \begin{align*}
    \mathbb{T}_{\mathfrak{U}_x^{\dag}}=
        \mathrm{Tot}\left(\mathfrak{g}_x^{\lambda}\otimes \mathcal{O}_{\mathfrak{U}_x^{\dag}} \to (Y_x \oplus \mathfrak{n}_x^{\vee}) \otimes \mathcal{O}_{\mathfrak{U}_x^{\dag}} 
        \to V_x \otimes \mathcal{O}_{\mathfrak{U}_x^{\dag}} \right).
    \end{align*}
\end{proof}
Let $T_x\subset G_x^{\lambda}$ be a maximal torus.
Consider the maps
\begin{align}\label{u:proj}
 h_x\colon\fU_x^\dag\stackrel{\pi_x}{\rightarrow}
 \mathfrak U_x\rightarrow\mathfrak M,
 \quad r_x\colon\fU_x^\dag\rightarrow BG_x^{\lambda},
\end{align}
where $r_x$ is the natural projection of the quotient stack.

\begin{lemma}\label{lem:local-numdata}
Let $\ell_x$, $\delta_x^{\mathrm{res}}$ and $b_x$ be the restrictions of $\ell$, $\delta$ and $b$ along
$BG_x^{\lambda}\to\mathfrak M$. Then $\ell_x=\delta_P(\chi_x)$ for a rational
character $\chi_x\in\mathrm{Pic}(BG_x^{\lambda})_{\mathbb Q}$, and $b_x$ is a
positive definite rational quadratic form on $N(T_x)_{\mathbb Q}$
invariant under $W_{G_x^{\lambda}}$, where $W_{G_x^{\lambda}}$ is the Weyl group of $G_x^{\lambda}$. Moreover
by shrinking $U_x$ if necessary, we have 
\begin{align}\label{pull-back:id}
 h_x^*\ell=r_x^*\ell_x,\quad h_x^*b=r_x^*b_x, \quad h_x^* \delta=r_x^* \delta_x^{\mathrm{res}}.
\end{align}
\end{lemma}
\begin{proof}
By~\cite[Example~3.1.2]{BPT} and~\cite[Example~3.2.5]{BHIK},
we have 
\begin{align*}
 \mathrm{CL}(BG_x^{\lambda})_{\mathbb Q}
 \cong N(T_x)_{\mathbb Q}/W_{G_x^{\lambda}}.
\end{align*}
The $W_{G_x^{\lambda}}$-invariant rational linear forms on $N(T_x)$ correspond to rational characters on $G_x^{\lambda}$
and quadratic norms correspond to $W_{G_x^{\lambda}}$-invariant positive definite quadratic
forms on $N(T_x)_{\mathbb{Q}}$; see also~\cite[Section~3.1.9]{BPT}. This proves the assertions
about $\chi_x$ and $b_x$.
For the assertion (\ref{pull-back:id}), see the argument of~\cite[Section~3.4.12]{BPT}.
\end{proof}

We consider the diagram induced by the cocharacter $\lambda \colon \mathbb{G}_m \to G_x^{\lambda}$
as in (\ref{diagram:fU})
\begin{align}\label{dia:Udaglambda}
    (\fU^{\dag}_x)^{\lambda} \stackrel{q_x^{\dag}}{\leftarrow} (\fU^{\dag}_x)^{\lambda \geq 0} \stackrel{p_x^{\dag}}{\hookrightarrow}\fU^{\dag}_x.
\end{align}
The stack $(\fU^{\dag}_x)^{\lambda \geq 0}$ is a closed 
Kempf--Ness stratum with respect to $(\chi_x,b_x)$. 

\begin{lemma}\label{lem:wU}
The derived stack $\fU_x^{\dag}$ and its $\Theta$-stratum $(\fU_x^{\dag})^{\lambda \geq 0} \hookrightarrow \fU_x^{\dag}$
satisfy Assumption~\ref{assum:dag}. 
\end{lemma}
\begin{proof}
The symmetric condition follows from Lemma~\ref{lem:symU}, and Assumption~\ref{assum:dag} (iii) follows from the fact that the 
$\lambda$-weights of $V_x=\cH^1(\mathbb{T}_{\mathfrak{M}}|_{x})$
are non-positive by Assumption~\ref{assum:dag} (iii)
for $\mathfrak{S} \hookrightarrow \mathfrak{M}$. The rest of the properties follow immediately. 
\end{proof}

\subsection{Semiorthogonal decomposition: global case}
In this subsection, we prove Claim~\ref{claim:1} in a global 
setting under Assumption~\ref{assum:dag} by the reduction to 
the case of the local model. 

Let $(\mathfrak{M}, \mathfrak{S}, \mathfrak{Z})$ be a derived stack with a $\Theta$-stratum satisfying Assumption~\ref{assum:dag}. 
For each closed point $z\in \mathfrak{Z}(k)$ corresponding to $x\in \mathfrak{M}(k)$ 
and $\lambda \colon \mathbb{G}_m \to G_x$, we have 
the following diagram by the construction in the previous subsection 
\begin{align}\label{dia:igx}
\xymatrix{
    \fU_x^{\dag} & \linclusion_-{i_x} \mathfrak{U}_x \ar[r]^-{g_x} & \mathfrak{M}.
    }
\end{align}
Here $i_x$ is the closed immersion induced by $0\in \mathfrak{n}_x^{\vee}$. 
Since $g_x$ is smooth and $i_x$ is a closed immersion, we have the functor 
\begin{align}\label{funct:tilde}
    (i_x)_{*}g_x^* \colon \Coh(\mathfrak{M}) \to \Coh(\fU_x^{\dag}).
\end{align}

\begin{lemma}\label{lem:ig}
The functor (\ref{funct:tilde}) restricts to the functor 
\begin{align*}
    (i_x)_{*}g_x^* \colon \LL(\mathfrak{M})_{\delta} \to \LL(\fU_x^{\dag})_{\delta_x}.
\end{align*}
Here $\delta_x\in \mathrm{CL}(\fU_x^{\dag})_{\mathbb{Q}}^{\vee}$ is given by 
\begin{align}\label{deltax}
\delta_x:=h_x^* \delta -\delta_P(\pi_x^* \det \mathbb{L}_{i_x}^{1/2}).
\end{align}
\end{lemma}
\begin{proof}
    The lemma follows from Proposition~\ref{prop:functL}, as it satisfies the 
    condition (\ref{c1zero}) (by shrinking $U_x$ if necessary)
    \begin{align*}
        \mathbb{L}_{i_x}[-1]=\mathbb{L}_{g_x}^{\vee}=\mathfrak{n}_x \otimes \mathcal{O}_{\mathfrak{U}_x}.
    \end{align*}
\end{proof}

We next consider the compatibility of the functor (\ref{funct:tilde}) with the semiorthogonal 
decomposition of $\Coh(\fU_x^{\dag})$. 
\begin{lemma}\label{lem:CarU}
By shrinking $\fU_x$ if necessary, 
we have the following commutative diagram 
\begin{align}\label{dia:comU}
\xymatrix{
(\fU_x^{\dagger})^{\lambda}
  \ar@{}[dr]|{\square}
&
(\fU_x^{\dagger})^{\lambda \geq 0}
  \ar[l]_{q_x^{\dagger}}
  \ar@{^{(}->}[r]^-{p_x^{\dag}}
&
\fU_x^{\dagger}
\\
\fU_x^{\lambda}
  \ar[d]_-{g_x^{\lambda}}
   \ar[u]_-{\cong}^-{i_x^{\lambda}}
&
\fU_x^{\lambda \geq 0}
  \ar[l]_{q_x}
  \ar@{^{(}->}[r]^-{p_x}
  \ar[d]_-{g_x^{\lambda \geq 0}}
   \ar[u]_-{\cong}^-{i_x^{\lambda \geq 0}}
  \diasquare
&
\fU_x
  \ar[d]_-{g_x}
    \ar@{^{(}->}[u]^-{i_x}
\\
\mathfrak{Z}
&
\mathfrak{S}
\inclusion^-{p}
  \ar[l]_{q}
&
\mathfrak{M}
}
\end{align}
Here the upper-left and lower-right squares are Cartesian, 
and the map $g_x^{\lambda}$ is \'etale and induces isomorphisms on 
stabilizer groups. 
\end{lemma}
\begin{proof}
Since the $\lambda$-weights of $\mathfrak{n}_x^{\vee}$ are negative, the vertical maps $i_x^{\lambda}$ and $i_x^{\lambda\geq0}$ are isomorphisms. Hence the upper-left square is Cartesian. 

We show that the lower-right square is Cartesian. 
Since $\lambda$ is central in $G_x^{\lambda}$, the evaluation map $\fU_x^{\lambda \geq 0}\to\mathfrak U_x$ is a closed immersion.
By shrinking $\fU_x$ if necessary, we may assume that 
$\fU_x^{\lambda \geq 0}$ is connected, hence it is a connected 
component of $\mathrm{Filt}(\mathfrak{U}_x)$.
Then the induced map $\mathrm{Filt}(g_x)$ on 
the stacks of filtrations gives a map 
$g_x^{\lambda \geq 0} \colon \mathfrak{U}_x^{\lambda \geq 0}
\to \mathfrak{S}$ such that the diagram (\ref{dia:comU}) commutes. 
It is enough to show that the induced closed immersion 
\begin{align*}
\alpha \colon 
    \fU_x^{\lambda \geq 0}\hookrightarrow \fU_x \times_{\mathfrak{M}}\mathfrak{S}
\end{align*}
is an equivalence of derived stacks. 

The point $o\in \fU_x$ together with $\lambda \colon \mathbb{G}_m \to G_x^{\lambda}$ gives a point $f_o\in \fU_x^{\lambda \geq 0}$. 
By the formula of the tangent complexes for mapping stacks
\cite[Lemmas~1.2.2--1.2.3]{Halpinstab}, we have
\[
\mathbb T_\alpha|_{f_o}\simeq
\operatorname{fib}\bigl(\mathfrak n_x^{\lambda\geq0}\to\mathfrak n_x\bigr)=0.
\]
Thus $\alpha$ is étale at $f_o$, hence an equivalence over an open neighborhood of $o$.
By shrinking $\fU_x$ if necessary, the map $\alpha$ is an equivalence. 

Similarly $g_x^{\lambda}$ is induced by $\mathrm{Grad}(g_x)$, 
which is \'etale at $f_o$ since $\mathbb{L}_{g_x^{\lambda}}|_{f_o}=\mathfrak{n}_x^{\lambda}=0$. 
Moreover $g_x^{\lambda}$ induces an isomorphism on stabilizer 
groups at $f_o$, $\id \colon G_x^{\lambda} \to G_x^{\lambda}$.  
Therefore by shrinking $\fU_x$ if necessary, we obtain the claim for $g_x^{\lambda}$. 
\end{proof}
We also have the following lemma: 
\begin{lemma}\label{lem:central-weight-local}
The central weight $m_{\delta}\in \mathbb{Q}$ given by (\ref{cent:m}) equals the $\lambda$-weight of $\delta_x+\delta_P(\det \mathbb{L}_{p_x^{\dag}}^{1/2})$.
\end{lemma}
\begin{proof}
    The lemma follows from the right squares in (\ref{dia:comU}) and the definition of $\delta_x$ in (\ref{deltax}); 
    the central weight of $\det \mathbb{L}_{p}^{1/2}$ on $\mathfrak{Z}$
    equals the $\lambda$-weight of $\det \mathbb{L}_{p_x}^{1/2}$, 
    which equals the $\lambda$-weight of $\det \mathbb{L}_{p_x^{\dag}}^{1/2} \otimes \det \mathbb{L}_{i_x}^{-1/2}$. 
\end{proof}

Let $\Upsilon^{\dag}_x$ be the functor defined from the diagram (\ref{dia:Udaglambda})
\begin{align*}
    \Upsilon^{\dag}_x :=p^{\dag}_{x*}(q^{\dag}_x)^*  \colon \Coh((\fU_x^{\dag})^{\lambda}) \to \Coh(\fU_x^{\dag}).
\end{align*}
By Theorem~\ref{thm:window}, Lemma~\ref{lem:wU}, we have the semiorthogonal 
decomposition 
\begin{align}
    &\Coh(\fU_x^{\dag})=\\
  \notag  &\left\langle \ldots,\Upsilon^{\dag}_x(\Coh((\fU_x^{\dag})^{\lambda})_{\lfloor m_{\delta} \rfloor}), \mathcal{W}^{\dag}_x(m_{\delta}), \Upsilon^{\dag}_x(\Coh((\fU_x^{\dag})^{\lambda})_{\lfloor m_{\delta} \rfloor+1}), \ldots \right\rangle. 
\end{align}
Here $m_{\delta}\in \mathbb{Q}$ is given as in Lemma~\ref{lem:central-weight-local}, and 
$\lfloor m_{\delta} \rfloor\in  \mathbb{Z}$ is its round-down. Moreover 
we have the equivalence 
\begin{align}\label{sod:dag2}
   (j_x^{\dag})^* \colon  \mathcal{W}^{\dag}_x(m_{\delta}) \stackrel{\sim}{\to} \Coh((\fU_x^{\dag})^{\circ})
\end{align}
where $j_x^{\dag}$ is the open immersion 
\begin{align*}
    j_x^{\dag} \colon (\fU_x^{\dag})^{\circ}:=\fU_x^{\dag} \setminus 
    (\fU_x^{\dag})^{\lambda \geq 0}\hookrightarrow 
    \fU_x^{\dag}.
\end{align*}
\begin{lemma}\label{lem:compatisdo}
The functor (\ref{funct:tilde}) is compatible with the semiorthogonal decompositions (\ref{sod:dag}), i.e. 
\begin{align*}
   & (i_x)_{*}g_x^* \colon \Upsilon(\Coh(\mathfrak{Z})_{a}) \to 
    \Upsilon_x^{\dag}(\Coh((\fU_x^{\dag})^{\lambda})_{a}), \\
    &(i_x)_{*}g_x^* \colon \mathcal{W}(m_{\delta}) \to \mathcal{W}^{\dag}_x(m_{\delta}).
\end{align*}
\end{lemma}
\begin{proof}
    The functor (\ref{funct:tilde}) is compatible with $\Upsilon$ and $\Upsilon_x^{\dag}$ by 
    the base change of the diagram (\ref{dia:comU}). 
    It sends $\mathcal{W}(m_{\delta})$ to $\mathcal{W}^{\dag}_x(m_{\delta})$ by 
    Lemma~\ref{lem:nuineq} applied to $(g_x, i_x)$, together with Lemma~\ref{lem:wtmin} and 
    Lemma~\ref{lem:wtW}. Therefore the lemma holds.  
\end{proof}

By combining the results so far, we have the following proposition, which is the main result in this section: 
\begin{prop}\label{prop:global}
Let $\mathfrak{M}$ be a derived stack with a $\Theta$-stratum 
$\mathfrak{S} \hookrightarrow \mathfrak{M}$ satisfying Assumption~\ref{assum:dag}. Then the statements of Claim~\ref{claim:1} 
hold.
\end{prop}
\begin{proof}
For each closed point $z\in \mathfrak{Z}(k)$ corresponding to $x\in \mathfrak{M}(k)$ 
and $\lambda \colon \mathbb{G}_m \to G_x$, we take the diagram of 
derived stacks 
\begin{align*}\fU_x^{\dag} \stackrel{i_x}{\hookleftarrow} \fU_x \stackrel{g_x}{\to} \mathfrak{M}
\end{align*}
as in (\ref{dia:igx}). For $\cE \in \LL(\mathfrak{M})_{\delta}$, 
take the distinguished triangle (\ref{sod:l1}). 
By applying $(i_x)_{*}g_x^*$, we obtain the distinguished triangle 
in $\Coh(\fU_x^{\dag})$
\begin{align*}
    (i_x)_{*}g_x^* \cE_1 \to (i_x)_{*}g_x^{*}\cE \to (i_x)_{*}g_x^* \cE_2.
\end{align*}

By Lemma~\ref{lem:ig}, we have 
\begin{align*}
    (i_x)_{*}g_x^* \cE \in \LL(\fU_x^{\dag})_{\delta_x}.
\end{align*}
Also by Lemma~\ref{lem:compatisdo}, we have 
\begin{align*}
(i_x)_{*}g_x^* \cE_1 \in \mathcal{W}^{\dag}_x(m_{\delta}), \ 
(i_x)_{*}g_x^* \cE_2 \in \Upsilon_x^{\dag}(\Coh((\fU_x^{\dag})^{\lambda})_{m_{\delta}}).
\end{align*}
Then by the result of the local model in 
Proposition~\ref{prop:sym2}, we have 
\begin{align*}
    (i_x)_{*}g_x^* \cE_1 \in \mathcal{W}^{\dag}_x(m_{\delta}) \cap \LL(\fU_x^{\dag})_{\delta_x}.
\end{align*}
Since the above holds for any choice of closed point 
$z\in \mathfrak{Z}(k)$, by 
Proposition~\ref{prop:closed2}, Lemma~\ref{lem:lift} below 
and the property $j^*\cE_1 \cong j^* \cE\in \LL(\mathfrak{M}^{\circ})_{\delta}$, we conclude that $\cE_1 \in \LL(\mathfrak{M})_{\delta}$. 

It remains to show that 
the functor 
\begin{align*}
    j^* \colon \mathcal{W}(m_{\delta}) \cap \LL(\mathfrak{M})_{\delta} \to \LL(\mathfrak{M}^{\circ})_{\delta}
\end{align*}
is an equivalence. By Theorem~\ref{thm:window}, the above functor is fully faithful. It is enough to show that the above functor is essentially 
surjective. For any $A\in \LL(\mathfrak{M}^{\circ})_{\delta}$, let $A' \in \mathcal{W}(m_{\delta})$ 
be its lift which exists and is unique up to isomorphism by Theorem~\ref{thm:window}. It is enough to show that $A'\in \LL(\mathfrak{M})_{\delta}$. 

We have the Cartesian squares 
\begin{align*}
    \xymatrix{
(\fU_x^{\dag})^{\circ} \dinclusion_-{j_x^{\dag}} \diasquare & (\fU_x)^{\circ} \linclusion_-{i_x^{\circ}} \ar[r]^-{g_x^{\circ}} \dinclusion_-{j_x} \diasquare & \mathfrak{M}^{\circ} \dinclusion_-{j} \\
   \fU_x^{\dag} & \linclusion_-{i_x} \fU_x \ar[r]^-{g_x} & \mathfrak{M}.
    }
\end{align*}
By the above diagram, we have the following commutative diagram 
\begin{align}\label{dia:WCoh}
    \xymatrix{
\mathcal{W}(m_{\delta}) \ar[r]^-{(i_x)_{*}g_x^*} \ar[d]_-{j^*}^-{\sim} & \mathcal{W}^{\dag}_x(m_{\delta}) \ar[d]^-{(j_x^{\dag})^{*}}_-{\sim} \\
\Coh(\mathfrak{M}^{\circ}) \ar[r]^-{(i_x^{\circ})_{*}(g_x^{\circ})^*} & \Coh((\fU_x^{\dag})^{\circ}).
    }
\end{align}
Here the vertical arrows are equivalences by Theorem~\ref{thm:window}. 
The above diagram restricts to the commutative diagram 
\begin{align}\label{dia:WCohL}
    \xymatrix{
\mathcal{W}(m_{\delta}) \cap \LL(\mathfrak{M})_{\delta} \ar[r]^-{(i_x)_{*}g_x^*} \ar[d]_-{j^*} & \mathcal{W}^{\dag}_x(m_{\delta}) \cap \LL(\fU_x^{\dag})_{\delta_x} \ar[d]^-{(j_x^{\dag})^{*}}_-{\sim} \\
\LL(\mathfrak{M}^{\circ})_{\delta} \ar[r]^-{(i_x^{\circ})_{*}(g_x^{\circ})^*} & \LL((\fU_x^{\dag})^{\circ})_{\delta_x}.
    }
\end{align}
The right vertical arrow is an equivalence by Proposition~\ref{prop:sym2}. 

By the commutative diagrams (\ref{dia:WCoh}), (\ref{dia:WCohL}), 
the object 
\begin{align*}
    (i_x)_{*}g_x^* A' \in \mathcal{W}_x^{\dag}(m_{\delta})
\end{align*}
satisfies 
\begin{align*}
    (j_x^{\dag})^*   (i_x)_{*}g_x^* A' \cong (i_x^{\circ})_{*}(g_x^{\circ})^* A \in 
    \LL((\fU_x^{\dag})^{\circ})_{\delta_x}.
\end{align*}
Again using the diagrams (\ref{dia:WCoh}), (\ref{dia:WCohL}), we conclude that 
\begin{align*}(i_x)_{*}g_x^* A'\in \mathcal{W}^{\dag}_x(m_{\delta}) \cap \LL(\fU_x^{\dag})_{\delta_x}
\subset \LL(\fU_x^{\dag})_{\delta_x}.
\end{align*}
Since the above holds for any choice of closed point 
$z\in \mathfrak{Z}(k)$, by 
Proposition~\ref{prop:closed2}, Lemma~\ref{lem:lift} 
and the property $j^*A' \cong A \in \LL(\mathfrak{M}^{\circ})_{\delta}$, we conclude that $A' \in \LL(\mathfrak{M})_{\delta}$.

\end{proof}

We have used the following lemma: 
\begin{lemma}\label{lem:lift}
For any map
$\nu\colon\bgm\to\mathfrak S$, there is a closed point $z\in\mathfrak Z(k)$, with associated local model
$g_x\colon\mathfrak U_x\to\mathfrak M$, such that
$p\circ\nu$ lifts to $\mathfrak U_x$.
\end{lemma}

\begin{proof}
We set $\overline\nu=q\circ\nu \colon \bgm \to \mathfrak{Z}$, and let $z$ be the closed point
corresponding to its image in $Z$.
Let $\fU_x \to \mathfrak{M}$ be the associated local model, 
and consider the \'etale map $g_x^{\lambda} \colon \fU_x^{\lambda} \to \mathfrak{Z}$ in Lemma~\ref{lem:CarU}. 
The map $g_x^{\lambda}$ is \'etale, sends the closed point $o$ to $z$, and induces an isomorphism on their stabilizer groups. By Luna's fundamental lemma for good moduli spaces~\cite[Theorem~6.10]{AlperLocalQuotient}, after restricting to a saturated neighborhood of $o$ in $\fU_x^{\lambda}$, the map on classical truncations is the base change of an \'etale map $Z_x\to Z$. The composite $\bgm\xrightarrow{\overline\nu}\mathfrak Z\to Z$ factors through the point $[z]\in Z(k)$. Choosing the point of $Z_x$ represented by $o$ therefore gives a lift
\begin{align}\label{lift:nubar}
\widetilde{\overline\nu}\colon\bgm
\longrightarrow \mathfrak{U}_x^{\lambda}
\subset\operatorname{Grad}(\mathfrak U_x).
\end{align}

In order to distinguish the two grading parameters, 
we use $(-)_a, (-)_b$ for the two $\mathbb{G}_m$, and 
write
\[
\nu\colon B\mathbb G_{m,a}\to
\mathfrak S\subset
\operatorname{Map}(\Theta_b,\mathfrak M).
\]
Thus $\nu$ corresponds to a morphism
\[
F_\nu\colon B\mathbb G_{m,a}\times\Theta_b
\to\mathfrak M.
\]
Its restriction to
$B\mathbb G_{m,a}\times B\mathbb G_{m,b}$
is the morphism corresponding to
$\overline\nu=q\circ\nu$.

By exchanging the factors,
$F_\nu$ determines a morphism
\[
\alpha_\nu\colon\Theta_b
\to\operatorname{Grad}_a(\mathfrak M),
\quad
\operatorname{Grad}_a(\mathfrak M)
:=\operatorname{Map}(B\mathbb G_{m,a},\mathfrak M).
\]
Its associated graded is
\[
\beta_\nu\colon B\mathbb G_{m,b}
\longrightarrow\operatorname{Grad}_a(\mathfrak M).
\]
More explicitly, $\beta_\nu$ corresponds to the composite
\[
B\mathbb G_{m,b}\times B\mathbb G_{m,a}
\xrightarrow{\mathrm{swap}}
B\mathbb G_{m,a}\times B\mathbb G_{m,b}
\xrightarrow{F_\nu|_{B\mathbb G_{m,a}\times B\mathbb G_{m,b}}}
\mathfrak M.
\]

The lift $\widetilde{\overline\nu}$ constructed in (\ref{lift:nubar}) 
corresponds to a morphism
\[
B\mathbb G_{m,a}\times B\mathbb G_{m,b}
\longrightarrow\mathfrak U_x.
\]
By exchanging its two factors, we obtain a lift
\[
\widetilde\beta_\nu\colon B\mathbb G_{m,b}
\longrightarrow\operatorname{Grad}_a(\mathfrak U_x)
\]
of $\beta_\nu$. Hence we have a commutative diagram
\[
\begin{tikzcd}
B\mathbb G_{m,b}
  \arrow[r, "\widetilde\beta_\nu"]
  \arrow[d]
&
\operatorname{Grad}_a(\mathfrak U_x)
  \arrow[d, "\operatorname{Grad}_a(g_x)"]
\\
\Theta_b
  \arrow[r, "\alpha_\nu"]
&
\operatorname{Grad}_a(\mathfrak M).
\end{tikzcd}
\]

The morphism $\operatorname{Grad}_a(g_x)$ is smooth by
\cite[Proposition~1.3.1(4)]{Halpinstab}.
By applying the pointwise lifting argument in the proof of
\cite[Lemma~4.4.4]{Halpinstab}, we obtain a lift
\[
\widetilde\alpha_\nu\colon\Theta_b
\longrightarrow\operatorname{Grad}_a(\mathfrak U_x)
\]
of $\alpha_\nu$.
By evaluating at $1\in\Theta_b$, we obtain a morphism
\[
\widetilde\alpha_\nu(1)\colon
B\mathbb G_{m,a}\longrightarrow\mathfrak U_x
\]
whose composition with $g_x$ is isomorphic to $p\circ\nu$.
This is the required lift.
\end{proof}

So far we assumed that the center $\mathfrak{Z}$ is connected, see Remark~\ref{rmk:connected}.
If $\mathfrak{Z}$ is not connected, by 
applying Proposition~\ref{prop:global} for each connected component, we 
obtain the following: 
\begin{thm}\label{thm:sodlim}
Under Assumption~\ref{assum:dag},
there is a semiorthogonal decomposition 
\begin{align*}
    \LL(\mathfrak{M})_{\delta}=\left\langle \Upsilon (\LL(\mathfrak{Z})_{\delta'}), j_{!}\LL(\mathfrak{M}^{\circ})_{\delta} \right\rangle.
\end{align*}
    Here $\Upsilon$ is fully faithful, $\delta'=(\delta+\delta_P(\det \mathbb{L}_p^{1/2}))|_{\mathfrak{Z}}$, 
    and $j_!$ is the fully faithful left adjoint of 
    $j^* \colon \LL(\mathfrak{M})_{\delta} \to 
    \LL(\mathfrak{M}^{\circ})_{\delta}$.
\end{thm}

We have the following corollary: 
\begin{cor}\label{cor:jleft}
For an object $\cE \in \LL(\mathfrak{M})_{\delta}$, it lies in 
$j_{!}\LL(\mathfrak{M}^{\circ})_{\delta}$ if and only if 
for any $z \in \mathfrak{Z}(k)$ with corresponding map 
and regularization 
\begin{align*}
    \nu \colon \bgm \to \mathfrak{M}, \ 
    \nu_V^{\mathrm{reg}} \colon V[-1]/\mathbb{G}_m \to 
    \mathfrak{M},
\end{align*}
the inequality (\ref{wt:condL}) is strict on the left:
\begin{align}\label{wt:condL2.5}
   &\wt( \iota_{*}(\nu_V^{\mathrm{reg}})^*(\cE))
   \subset \\ 
 \notag  &\left(\frac{1}{2}c_1(\nu^* \mathbb{L}_{\mathfrak{M}}^{<0})
   +c_1((V^{\vee})^{<0}), \frac{1}{2}c_1(\nu^* \mathbb{L}_{\mathfrak{M}}^{>0})
   +c_1((V^{\vee})^{>0})
   \right]+\delta(\nu).
\end{align}
\end{cor}
\begin{proof}
If $\mathfrak{Z}$ is connected, 
by the construction of $j_{!}$ from the equivalence (\ref{equiv:WK}), we have 
\begin{align*}
    j_{!}\LL(\mathfrak{M}^{\circ})_{\delta}=\mathcal{W}(m_{\delta}) \cap \LL(\mathfrak{M})_{\delta}.
\end{align*}
Then the corollary follows from Lemma~\ref{lem:wtmin}, 
Lemma~\ref{lem:wtW} and Lemma~\ref{lem:choicem}.
In general (when $\mathfrak{Z}$ is not necessarily connected), 
the corollary follows by applying the above argument for each connected component. 
\end{proof}

\subsection{Semiorthogonal decompositions (in general)}
Let $\mathfrak{M}$ be a quasi-smooth and symmetric stack, which 
is not necessarily quasi-compact. Let
$\ell\in\mathrm{CL}(\mathfrak M)_{\mathbb Q}^{\vee}$ and let $b$ be a
positive definite rational quadratic norm on $\mathrm{CL}(\mathfrak M)$.
Suppose that $(\ell,b)$ defines a locally finite $\Theta$-stratification
in the sense of Definition~\ref{def:ndata} 
\begin{align*}
\mathfrak{M}_{\preceq \mu}=\bigsqcup_{\mu' \in (0, \mu] \cap I} \mathfrak{S}_{\mu'} \sqcup \mathfrak{M}^{\ell\text{-ss}}, \ \mathfrak{M}=\bigcup_{\mu\in I_{>0}}\mathfrak{M}_{\preceq \mu}.
\end{align*}
We have the following result: 
\begin{thm}\label{thm:cgen}
Suppose that each
$\Theta$-stratum $p_{\mu} \colon \mathfrak{S}_{\mu} \hookrightarrow \mathfrak{M}_{\preceq \mu}$ satisfies Assumption~\ref{assum:dag}.
Then $\IndL(\mathfrak{M})_{\delta}$ for $\delta \in \mathrm{CL}(\mathfrak{M})_{\mathbb{Q}}^{\vee}$ is compactly generated, 
and its subcategory of compact objects $\LL(\mathfrak{M})_{\delta}$ admits a semiorthogonal decomposition 
\begin{align*}
    \LL(\mathfrak{M})_{\delta}=\left\langle \Upsilon_{\mu} (\LL(\mathfrak{Z}_{\mu})_{\delta_{\mu}}), j_{!}\LL(\mathfrak{M}^{\ell\text{-ss}})_{\delta} : \mu \in I_{>0}\right\rangle.
\end{align*}
    Here $\delta_{\mu}=(\delta+\delta_P(\det \mathbb{L}_{p_{\mu}}^{1/2}))|_{\mathfrak{Z}_{\mu}}$, 
    and $\Upsilon_{\mu}$, $j_{!}$ are fully faithful functors 
    \begin{align*}
        \Upsilon_{\mu} \colon \LL(\mathfrak{Z}_{\mu})_{\delta_{\mu}} \hookrightarrow \LL(\mathfrak{M})_{\delta}, \ 
        j_{!} \colon \LL(\mathfrak{M}^{\ell\text{-ss}})_{\delta} \hookrightarrow \LL(\mathfrak{M})_{\delta}
    \end{align*}
    where $j_{!}$ is the left adjoint of $j^* \colon \LL(\mathfrak{M})_{\delta} \to \LL(\mathfrak{M}^{\ell\text{-ss}})_{\delta}$
    for the open immersion $j \colon \mathfrak{M}^{\ell\text{-ss}} \hookrightarrow \mathfrak{M}$. The order of the semiorthogonal 
    decomposition is given by the decreasing order of $\mu$. 
\end{thm}
\begin{proof}
Using Theorem~\ref{thm:sodlim}, the proof is similar to the argument in the proof of~\cite[Theorem~7.16]{PTlim}. We prove it for completeness. Below we use Proposition~\ref{prop:dgcat} for the limits/colimits of dg-categories. 

For $\mu_2 \leq \mu_1$, let $j_{12}$ be the open immersion 
\begin{align}\label{map:j12}
    j_{12} \colon \mathfrak{M}_{\preceq \mu_2} \hookrightarrow
    \mathfrak{M}_{\preceq \mu_1}. 
\end{align}
We have the pull-back functor
\begin{align}\label{pback2.0}
    j_{12}^{\ast} \colon \LL(\mathfrak{M}_{\preceq \mu_1})_{\delta}
    \to \LL(\mathfrak{M}_{\preceq \mu_2})_{\delta}
\end{align}
which induces the continuous functor by taking its ind-completion
(since $\mathfrak{M}_{\preceq \mu}$ is QCA)
\begin{align}\label{pback2}
    j_{12}^{\ast} \colon \IndL(\mathfrak{M}_{\preceq \mu_1})_{\delta}
    \to \IndL(\mathfrak{M}_{\preceq \mu_2})_{\delta}. 
\end{align}
By the open covering (\ref{open:Mu})
and the definition of 
$\IndL(\mathfrak{M})_{\delta}$, 
we have an equivalence 
\begin{align*}
\IndL(\mathfrak{M})_{\delta} \simeq \lim_{\mu \in I_{>0}}\IndL(\mathfrak{M}_{\preceq \mu})_{\delta}
\end{align*}
where the limit is taken with respect to the pull-backs (\ref{pback2}).

For the open immersion $j_{12}$ in (\ref{map:j12}), 
by applying Theorem~\ref{thm:sodlim} for each $\mu \in I \cap (\mu_2, \mu_1]$ (which is finite by the locally finite assumption on $I$), there is a fully faithful 
left adjoint of the functor $j_{12}^{\ast}$ in (\ref{pback2.0})
\begin{align}\label{j12!}
	j_{12!} \colon \LL(\mathfrak{M}_{\preceq \mu_2})_{\delta} \hookrightarrow 
	\LL(\mathfrak{M}_{\preceq \mu_1})_{\delta}. 
	\end{align}	
By taking its ind-completion, we obtain the continuous fully faithful functor
\begin{align}\label{ind:j2}
	j_{12!} \colon \IndL(\mathfrak{M}_{\preceq \mu_2})_{\delta} \hookrightarrow 
	\IndL(\mathfrak{M}_{\preceq \mu_1})_{\delta}. 
	\end{align}
    Since $j_{12}^*$ admits continuous left adjoint $j_{12!}$ as above, 
by Proposition~\ref{prop:dgcat} (i) the dg-category 
$\IndL(\mathfrak{M})_{\delta}$ is also equivalent to the colimit 
\begin{align}\label{IndL:0}
	\IndL(\mathfrak{M})_{\delta} \simeq \colim_{\mu \in I_{>0}} \IndL(\mathfrak{M}_{\preceq \mu})_{\delta}
	\end{align}
where the colimit is taken 
with respect to the functors (\ref{ind:j2}). 

In particular for the open immersion \begin{align*}j_{\mu} \colon \mathfrak{M}_{\preceq \mu} \hookrightarrow \mathfrak{M},\end{align*} 
the pull-back 
\begin{align*}
	j_{\mu}^{\ast} \colon \IndL(\mathfrak{M})_{\delta} \to \IndL(\mathfrak{M}_{\preceq \mu})_{\delta}
	\end{align*}
admits a left adjoint 
\begin{align}\label{jc!}
	j_{\mu!} \colon \IndL(\mathfrak{M}_{\preceq \mu})_{\delta} \hookrightarrow \IndL(\mathfrak{M})_{\delta}
	\end{align}
which is fully faithful since each (\ref{ind:j2}) is fully faithful by Proposition~\ref{prop:dgcat} (iv). 

For $\mu_2 < \mu_1$, by Proposition~\ref{prop:dgcat} (iv) the natural morphism 
$j_{12!} \to j_{\mu_1}^{\ast}j_{\mu_2!}$
is an isomorphism since (\ref{jc!}) is fully faithful. 
In particular, the functor (\ref{jc!}) restricts to the functor 
\begin{align*}
    j_{\mu!} \colon \LL(\mathfrak{M}_{\preceq \mu})_{\delta} \hookrightarrow \Coh(\mathfrak{M}). 
\end{align*}
By Proposition~\ref{prop:dgcat} (iii), an object $\cE$ in the left-hand side of (\ref{IndL:0}) is compact if and only 
if 
\begin{align}\label{jcE}\mathcal{E} \cong j_{\mu!}\cE_{\mu} \in \Coh(\mathfrak{M})
	\end{align}
for some $\mu \in I_{>0}$ and 
$\cE_{\mu} \in \LL(\mathfrak{M}_{\preceq \mu})_{\delta}$. 

From the above argument together with Proposition~\ref{prop:dgcat} (iii), 
    we conclude that the dg-category $\IndL(\mathfrak{M})_{\delta}$ is compactly generated 
with the subcategory of compact objects given by
\begin{align}\label{LL:colim}
\LL(\mathfrak{M})_{\delta} \simeq \colim_{\mu \in I_{>0}} \LL(\mathfrak{M}_{\preceq \mu})_{\delta}. 
\end{align}
Here the colimit is taken with respect to (\ref{j12!}). 
The desired semiorthogonal decomposition follows from (\ref{LL:colim}) and
Theorem~\ref{thm:sodlim}.
\end{proof}

Recall the subcategory 
\begin{align*}
    \widetilde{\LL}(\mathfrak{M})_{\delta} \subset \IndL(\mathfrak{M})_{\delta}
\end{align*}
in 
Definition~\ref{defn:lim2}.
We have the following generalization of~\cite[Corollary~8.27]{PTlim}:
\begin{cor}\label{cor:Ltilde}
We have 
\begin{align}\label{LsubLtilde}
\LL(\mathfrak{M})_{\delta} \subset \widetilde{\LL}(\mathfrak{M})_{\delta}
\end{align}
and an object $\mathcal{E} \in \widetilde{\LL}(\mathfrak{M})_{\delta}$
lies in $\LL(\mathfrak{M})_{\delta}$ if and only if 
    there is $\mu \in I_{>0}$ such that for any $A \in \widetilde{\LL}(\mathfrak{M})_{\delta}$
    satisfying 
    \begin{align}\label{Zc}
\mathrm{Supp}(A) \subset \mathfrak{S}_{>\mu}:=\mathfrak{M} \setminus \mathfrak{M}_{\preceq \mu},
    \end{align}
    we have $\Hom(\mathcal{E}, A)=0$. 
\end{cor}
\begin{proof}
    The inclusion (\ref{LsubLtilde}) holds 
    by the description of compact objects (\ref{jcE}) and using the fact that the functor (\ref{ind:j2})
    preserves compact objects (by its construction). 

    We prove the second statement. Consider an object $\mathcal{E}$ 
    in $\LL(\mathfrak{M})_{\delta}$. It may be written as $\mathcal{E}=j_{\mu!}\mathcal{E}_{\mu}$
    as in (\ref{jcE}). Then, for an object $A$ satisfying (\ref{Zc}), we have 
    \begin{align*}
        \Hom(j_{\mu!}\mathcal{E}_{\mu}, A)=\Hom(\mathcal{E}_{\mu}, j_{\mu}^{\ast}A)=0. 
    \end{align*}
    Conversely, for $\mathcal{E} \in \widetilde{\LL}(\mathfrak{M})_{\delta}$, 
    we have the distinguished triangle 
    \begin{align}\label{tri:jc}
        j_{\mu!}j_{\mu}^{\ast}\mathcal{E} \to \mathcal{E} \to A
    \end{align}
    where $j_{\mu!}j_{\mu}^{\ast}\mathcal{E} \in \LL(\mathfrak{M})_{\delta}$
    and $A \in \widetilde{\LL}(\mathfrak{M})_{\delta}$. Since $j_{\mu}^*A=0$, the object $A$ 
    satisfies (\ref{Zc}). 
    If $\Hom(\mathcal{E}, A)=0$, then 
    $\cE$ is a direct summand of $j_{\mu!} j_{\mu}^* \cE$, 
    therefore we conclude that $\mathcal{E} \in \LL(\mathfrak{M})_{\delta}$. 
\end{proof}

\subsection{Limits and colimits of dg-categories}
In the proof of Theorem~\ref{thm:cgen}, we used the following 
general facts on limits/colimits of dg-categories. 
Let 
\begin{align*}
    i \mapsto \mathcal{C}_i, \ (i\to j) \mapsto (\phi_{ji} \colon \mathcal{C}_j \to \mathcal{C}_i)
\end{align*}
be a diagram of dg-categories indexed by $I^{\mathrm{op}}$ for a filtered category $I$, 
where each $\mathcal{C}_i=\Ind(\mathcal{C}_i^{\mathrm{cp}})$ is compactly generated 
by the subcategory of compact objects $\mathcal{C}_i^{\mathrm{cp}} \subset \mathcal{C}_i$, 
and $\phi_{ji}$ is a continuous functor. 
We have the limit of dg-categories $\mathcal{C}_i$ together with
the evaluation functor 
\begin{align}\label{C:lim}
    \mathcal{C}^{\mathrm{lim}}=\lim_{i\in I^{\mathrm{op}}} \mathcal{C}_i \stackrel{\mathrm{ev}_i}{\to} \mathcal{C}_i. 
\end{align}

Suppose that each $\phi_{ji}$ admits a continuous left adjoint 
\begin{align}\label{funct:phij}
\psi_{ij} \colon \mathcal{C}_i \to \mathcal{C}_j.     
\end{align}
Note that $\psi_{ij}$ must satisfy $\psi_{ij}(\mathcal{C}_i^{\mathrm{cp}}) \subset \mathcal{C}_j^{\mathrm{cp}}$
since it admits a continuous right adjoint $\phi_{ji}$. In this 
situation, by~\cite[Proposition~1.7.5]{DGbun}
there is a colimit 
\begin{align*}
    \mathcal{C}^{\mathrm{colim}}:=\colim_{i\in I}\mathcal{C}_i
\end{align*}
with respect to functors (\ref{funct:phij}), and 
there are functors
\begin{align}\notag
\psi_i \colon \mathcal{C}_i \to \mathcal{C}^{\mathrm{colim}}
\end{align}
which satisfy
$\psi_j \circ \psi_{ij} \cong \psi_i$. 
The following proposition is given in~\cite{DGbun}, 
and summarized in~\cite[Proposition~2.2]{PTlim}:
\begin{prop}\emph{(\cite{DGbun}, \cite[Proposition~2.2]{PTlim})}\label{prop:dgcat}

 (i) There is a canonical equivalence $\mathcal{C}^{\mathrm{lim}} \simeq \mathcal{C}^{\mathrm{colim}}$. 

 (ii) Under the equivalence in (i), the functor $\psi_i$ is a left adjoint of $\mathrm{ev}_i$. 

 (iii) We have 
 $\mathcal{C}^{\mathrm{colim}}=\Ind(\colim_{i\in I}\mathcal{C}_i^{\mathrm{cp}})$. In particular 
 an object $\mathcal{E} \in \mathcal{C}^{\mathrm{colim}}$ is compact if and only 
 if it is isomorphic to $\psi_i(\mathcal{E}_i)$ for some $i\in I$ and 
 $\mathcal{E}_i \in \mathcal{C}_i^{\mathrm{cp}}$. 

 (iv) If each $\psi_{ij}$ is fully faithful, then $\psi_i$ is also fully faithful. 
 In this case, for $i \to j$ in $I$ 
 the natural morphism $\psi_{ij} \to \mathrm{ev}_j \psi_i$ is an isomorphism. \end{prop}

\subsection{K\"unneth formula for limit categories}
As an application of the semiorthogonal decomposition in Theorem~\ref{thm:cgen}, we prove the K\"unneth formula for limit categories. 

Let $\mathfrak{M}_i$ be quasi-smooth and symmetric derived stacks as in Theorem~\ref{thm:cgen}, with $\Theta$-stratifications satisfying Assumption~\ref{assum:dag}. Assume also that the stacks $\mathfrak{M}_{i, \preceq 0}$ admit good moduli spaces. 
We also take $\delta_i \in \mathrm{CL}(\mathfrak{M}_i)_{\mathbb{Q}}^{\vee}$.
\begin{thm}\label{thm:kunneth}
In the above setting, there is an equivalence 
\begin{align*}
    \IndL(\mathfrak{M}_1)_{\delta_1}\otimes 
    \IndL(\mathfrak{M}_2)_{\delta_2} \stackrel{\sim}{\to} 
    \IndL(\mathfrak{M}_1 \times \mathfrak{M}_2)_{\delta_1 \boxtimes \delta_2}.
\end{align*}
In particular, by taking the subcategory of compact objects, there is an equivalence 
\begin{align}\label{equiv:M12}
    \LL(\mathfrak{M}_1)_{\delta_1}\otimes 
    \LL(\mathfrak{M}_2)_{\delta_2} \stackrel{\sim}{\to} 
    \LL(\mathfrak{M}_1 \times \mathfrak{M}_2)_{\delta_1 \boxtimes \delta_2}.
\end{align}
\end{thm}
\begin{proof}
    We first prove the case where $\mathfrak{M}_i$ are QCA and admit good moduli spaces. 
    In this case, we have an equivalence, see~\cite{MR3037900}
    \begin{align*}
        \IndCoh(\mathfrak{M}_1)\otimes \IndCoh(\mathfrak{M}_2) \stackrel{\sim}{\to} \IndCoh(\mathfrak{M}_1\times \mathfrak{M}_2).
    \end{align*}
    By taking the subcategory of compact objects, 
    we have 
    \begin{align}\label{rest:Phi}
      \Phi \colon   \Coh(\mathfrak{M}_1)\otimes \Coh(\mathfrak{M}_2) \stackrel{\sim}{\to} \Coh(\mathfrak{M}_1\times \mathfrak{M}_2).
    \end{align}
    By taking quadratic norms $b_{i,\mathrm{aux}}$ on $\mathrm{CL}(\mathfrak{M}_i)_{\mathbb{Q}}$, we have the semiorthogonal decompositions by Theorem~\ref{thm:BPT}
    \begin{align*}
        \Coh(\mathfrak{M}_i)=\left\langle \ast_{\lambda_i} \LL(\mathfrak{M}_i^{\lambda_i})_{\delta_{i, \lambda_i}} : \lambda_i \in \lvert\mathrm{CL}(\mathfrak{M}_i)_{\mathbb{Q}}\rvert\right\rangle.
    \end{align*}
    By taking the tensor products, we obtain the semiorthogonal decomposition
    \begin{align*}
         \Coh(\mathfrak{M}_1)\otimes \Coh(\mathfrak{M}_2) =\left\langle 
(\ast_{\lambda_1}\LL(\mathfrak{M}_1^{\lambda_1})_{\delta_{1,\lambda_1}})\otimes (\ast_{\lambda_2}\LL(\mathfrak{M}_2^{\lambda_2})_{\delta_{2,\lambda_2}}): 
          \lambda_i \in \lvert\mathrm{CL}(\mathfrak{M}_i)_{\mathbb{Q}}\rvert
         \right\rangle.
    \end{align*}
    
    On the other hand, note that 
    \begin{align*}\mathrm{CL}(\mathfrak{M}_1\times \mathfrak{M}_2)_{\mathbb{Q}}= \mathrm{CL}(\mathfrak{M}_1)_{\mathbb{Q}}
    \times  \mathrm{CL}(\mathfrak{M}_2)_{\mathbb{Q}}.
    \end{align*}
    By applying Theorem~\ref{thm:BPT} for $\mathfrak{M}_1\times \mathfrak{M}_2$, $\delta_1\boxtimes \delta_2$ and 
    $b_{1,\mathrm{aux}}+b_{2,\mathrm{aux}}$, we have the semiorthogonal decomposition 
    \begin{align*}
        \Coh(\mathfrak{M}_1 \times \mathfrak{M}_2)=\left\langle  \ast_{\lambda}\LL(\mathfrak{M}_1^{\lambda_1} \times \mathfrak{M}_2^{\lambda_2})_{\delta_{1,\lambda_1}\boxtimes \delta_{2,\lambda_2}} : 
        \lambda=(\lambda_1, \lambda_2) \in
        \lvert\mathrm{CL}(\mathfrak{M}_1\times \mathfrak{M}_2)_{\mathbb{Q}}\rvert
        \right\rangle.
    \end{align*}
    Since $*_{(\lambda_1, \lambda_2)}=\ast_{\lambda_1}\boxtimes \ast_{\lambda_2}$, the equivalence (\ref{rest:Phi})
    restricts to the fully faithful functor 
    \begin{align}\label{equiv:Phi2}
        \Phi \colon (\ast_{\lambda_1}\LL(\mathfrak{M}_1^{\lambda_1})_{\delta_{1,\lambda_1}}) \otimes 
         (\ast_{\lambda_2}\LL(\mathfrak{M}_2^{\lambda_2})_{\delta_{2,\lambda_2}})
         \hookrightarrow \ast_{(\lambda_1, \lambda_2)}\LL(\mathfrak{M}_1^{\lambda_1} \times \mathfrak{M}_2^{\lambda_2})_{\delta_{1,\lambda_1}\boxtimes \delta_{2,\lambda_2}}.
    \end{align}
    It follows from the above semiorthogonal decompositions and the equivalence (\ref{rest:Phi}) that the functor (\ref{equiv:Phi2}) is an equivalence. 
    In particular, we have the equivalence (\ref{equiv:M12}) in this case. 

    In general ($\mathfrak{M}_i$ are locally QCA, not necessarily having good moduli spaces), since the categories $\IndL(\mathfrak{M}_i)_{\delta_i}$ are compactly generated, they are dualizable; hence tensoring with them commutes with limits. Therefore we may assume that $\mathfrak{M}_i$ are QCA. Then we have the equivalence (\ref{rest:Phi}) (which holds without good moduli spaces), and from the definition of $\LL(\mathfrak{M}_i)_{\delta_i}$ it restricts to the fully faithful functor 
    \begin{align}\label{ff:limit}
      \Psi \colon   \LL(\mathfrak{M}_1)_{\delta_1}\otimes  \LL(\mathfrak{M}_2)_{\delta_2}
    \hookrightarrow \LL(\mathfrak{M}_1\times \mathfrak{M}_2)_{\delta_1\boxtimes \delta_2}.
    \end{align}
    Let 
    \begin{align*}\mathfrak{Z}_{i, \mu} \leftarrow \mathfrak{S}_{i, \mu} \hookrightarrow \mathfrak{M}_{i, \preceq \mu}
    \end{align*}
    be 
    $\Theta$-strata for $\mathfrak{M}_i$, indexed by $\mu\in I_i\setminus\{0\}$, where $I_i\subset\mathbb R_{\geq0}$ is the index set of each $\Theta$-stratification.
    Set $\mathfrak Z_{i,0}=\mathfrak S_{i,0}=\mathfrak M_{i,\preceq0}$ and $\delta_{i,0}=\delta_i|_{\mathfrak M_{i,\preceq0}}$.
    By Theorem~\ref{thm:cgen}, we have the semiorthogonal decomposition
    \begin{align*}
        \LL(\mathfrak{M}_1)_{\delta_1}\otimes  \LL(\mathfrak{M}_2)_{\delta_2}
           =\left\langle \LL(\mathfrak{Z}_{1, \mu_1})_{\delta_{1, \mu_1}}
           \otimes\LL(\mathfrak{Z}_{2, \mu_2})_{\delta_{2, \mu_2}} : \mu_1\in I_1,\ \mu_2\in I_2
           \right\rangle.
    \end{align*}
    Here in order to simplify the notation, we omit $!$-push forward and $\Upsilon_{\mu}$. 
    The order of the semiorthogonal components are given by the decreasing order of $\sqrt{\mu_1^2+\mu_2^2}$, choosing any order among pairs with the same value.

    On the other hand, we have the stratification 
    \begin{align}\label{strata:product}
        \mathfrak{M}_1 \times \mathfrak{M}_2=\bigsqcup_{\mu \in I_{12}} \mathfrak{S}_{\mu}, \ 
        \mathfrak{S}_{\mu}=\bigsqcup_{\substack{\mu_1\in I_1,\ \mu_2\in I_2\\\sqrt{\mu_1^2+\mu_2^2}=\mu}}
        \mathfrak{S}_{1, \mu_1}\times \mathfrak{S}_{2, \mu_2}.
    \end{align}
Here $I_{12}=\{\sqrt{\mu_1^2+\mu_2^2}:\mu_1\in I_1,\ \mu_2\in I_2\}$. The above stratification is a $\Theta$-stratification with respect to $(\ell, b)$, 
\begin{align*}
    \ell=\mathrm{pr}_1^*\ell_1+\mathrm{pr}_2^*\ell_2
       \in\mathrm{CL}(\mathfrak M_1\times\mathfrak M_2)_{\mathbb Q}^{\vee},\quad
    b=\mathrm{pr}_1^*b_1+\mathrm{pr}_2^*b_2.
\end{align*}
Here $b$ is the sum of the pulled-back quadratic forms; it is
positive on every non-degenerate graded point of the product.
Indeed for maps $f_i\colon\Theta\to\mathfrak M_i$, let $\mu_i$ be the
slope with respect to $(\ell_i,b_i)$. 
Then for the map 
\begin{align*}f=(f_1, f_2) \colon \Theta \to \mathfrak{M}_1 \times \mathfrak{M}_2,
\end{align*}
its slope with respect to $(\ell, b)$ is 
\begin{align*}
    \mu(f)=\frac{\mu_1 \sqrt{\lVert(f_1)_0 \rVert}+\mu_2 \sqrt{\lVert(f_2)_0 \rVert}}{\sqrt{\lVert(f_1)_0 \rVert+\lVert(f_2)_0 \rVert}}
    \leq \sqrt{\mu_1^2+\mu_2^2}. 
\end{align*}
Here the equality holds if and only if 
$(\sqrt{\lVert(f_1)_0 \rVert}, \sqrt{\lVert(f_2)_0 \rVert})
\in \mathbb{R}_{>0}(\mu_1, \mu_2)$. 

    The derived stack $\mathfrak{M}_1 \times \mathfrak{M}_2$
    and each $\Theta$-stratum 
    $\mathfrak{S}_{\mu} \hookrightarrow (\mathfrak{M}_1 \times \mathfrak{M}_2)_{\preceq \mu}$ satisfy Assumption~\ref{assum:dag}. By 
    applying Theorem~\ref{thm:cgen} for $\mathfrak{M}_1\times \mathfrak{M}_2$, we also obtain the semiorthogonal decomposition 
    \begin{align*}
        \LL(\mathfrak{M}_1 \times \mathfrak{M}_2)_{\delta_1 \boxtimes \delta_2}=\left\langle \LL(\mathfrak{Z}_{1, \mu_1}\times \mathfrak{Z}_{2, \mu_2})_{\delta_{1, \mu_1}\boxtimes \delta_{2, \mu_2}} : \mu_1\in I_1,\ \mu_2\in I_2 \right\rangle.
    \end{align*}
   The functor (\ref{ff:limit}) is compatible with the above embeddings, and it induces the fully-faithful functor 
    \begin{align*}
        \Psi \colon \LL(\mathfrak{Z}_{1, \mu_1})_{\delta_{1, \mu_1}}
           \otimes \LL(\mathfrak{Z}_{2, \mu_2})_{\delta_{2, \mu_2}}
           \hookrightarrow \LL(\mathfrak{Z}_{1, \mu_1}\times \mathfrak{Z}_{2, \mu_2})_{\delta_{1, \mu_1}\boxtimes \delta_{2, \mu_2}}.
    \end{align*}
    The above functor is an equivalence by the equivalence (\ref{equiv:M12}) in the case where $\mathfrak{M}_i$ admit good moduli spaces, since $\mathfrak{Z}_{i, \mu_i}$ admit good moduli spaces by Assumption~\ref{assum:dag}.
    Therefore the functor (\ref{ff:limit}) is an equivalence. 
\end{proof}

\section{Examples}\label{sec:toy}\label{subsec:toy}
In this section, we first discuss the semiorthogonal decomposition 
of the 
two-dimensional toy example, and then generalize it to higher 
dimension. We observe that the functor $j_{!}$ from the semistable locus is related to the 
iteration of flop-equivalences, which we call the $\infty$-flop functor. 

\subsection{The two-dimensional toy model}\label{subsec:toy-two}
For each $i\in 1/2 +\mathbb{Z}$, we take a copy of 
the quotient stack 
\begin{align*}
    \X_i=\mathbb{A}^2/\mathbb{G}_m
\end{align*}
where $\mathbb{G}_m$ acts on $\mathbb{A}^2$ by weight $(1, -1)$. 
For the cocharacter $\lambda=\id \colon \mathbb{G}_m \to \mathbb{G}_m$, 
we set 
\begin{align*}
    L_i^{\pm } :=\X_i^{\lambda^{\pm 1} \geq 0}, \
    \mathcal{X}_i^{\pm }:=\mathcal{X}_i \setminus L_i^{\pm}.
\end{align*}
We have the stratifications 
\begin{align*}
    \X_i=L_i^+ \sqcup \X_i^{+} =\X_i^- \sqcup L_i^-
\end{align*}
where each $\X_i^{\pm} \subset \X_i$ is an open substack isomorphic to $\mathbb{A}^1$, by the isomorphisms 
\begin{align}\label{isom:A1}
\mathbb{A}^1 \stackrel{\cong}{\to} \X_i^{-}, \ t\mapsto (1, t), \ 
\mathbb{A}^1 \stackrel{\cong}{\to} \X_i^{+}, \ t\mapsto (t, 1).
\end{align}

We glue $\X_i$ and $\X_{i+1}$ for $i\in 1/2 +\mathbb{Z}$ 
by the isomorphisms (\ref{isom:A1})
\begin{align*}\X_i^{-} \stackrel{\cong}{\leftarrow}
\mathbb{A}^1 \stackrel{\cong}{\to} \X_{i+1}^{+}.
\end{align*}
We then obtain the non-quasi-compact 
stack $\X$
    \begin{align}\label{stack:X}
        \X=\bigcup_{i\in 1/2+\mathbb{Z}} \X_i.
    \end{align}
By the construction, we have the stratification 
\begin{align}\label{strata:X}
\X=\cdots \sqcup L_{-3/2}^{-} \sqcup  L_{-1/2}^{-} \sqcup \X_{1/2}^+ \sqcup L_{1/2}^+ \sqcup L_{3/2}^+ \sqcup \cdots.
    \end{align}
    The good moduli space maps $\X_i \to \mathbb{A}^1$ glue to 
    give a map 
    \begin{align}\label{piX}
        \pi_{\X} \colon \X \to \mathbb{A}^1.
    \end{align}
    We define the closed substack $L_i \hookrightarrow \X_i$ by 
\begin{align*}
    L_i:=L_i^+ \cup L_i^- =(\pi_{\X}|_{\X_i})^{-1}(0).
\end{align*}
\begin{remark}\label{rmk:X0}
The map $\pi_{\X}$ in (\ref{piX}) is an isomorphism away from $0\in \mathbb{A}^1$. 
Over $0\in \mathbb{A}^1$, the fiber $\X_0:=\pi_{\X}^{-1}(0)$
is a $\mathbb{G}_m$-quotient of an infinite chain of $\mathbb{P}^1$
\begin{align*}
    \X_0=\left(\bigcup_{n\in \mathbb{Z}}\mathbb{P}^1_n  \right)/\mathbb{G}_m.
\end{align*}
Here $\mathbb{P}_n^1$ is a copy of $\mathbb{P}^1$, and $\infty \in \mathbb{P}^1_n$ and $0\in \mathbb{P}_{n+1}^1$ are glued 
transversally. 
It is realized as a $\mathbb{G}_m$-rigidified moduli stack of 
rank-one torsion-free sheaves on the union of two copies of $\mathbb{P}^1$ 
intersecting transversally at one point; see \cite[Subsection~4.5]{TodaSurfaceLimits}.

The stack $\X_0$ has $\mathbb{Z}$-copies of open substacks isomorphic to 
$\mathrm{pt}$, each of which degenerates to two closed points with $\mathbb{G}_m$-stabilizer groups. It is depicted as follows: 
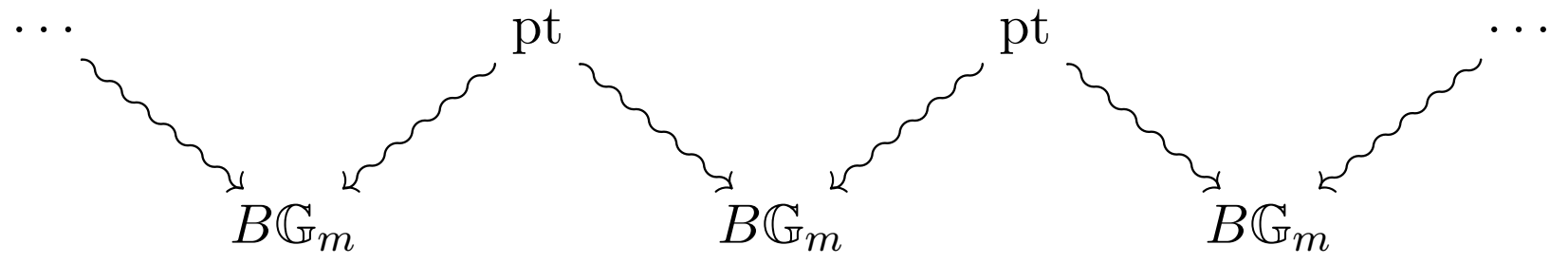
\begin{align*}
    \xymatrix{
\cdots \ar@{~>}[rd]& & \mathrm{pt} \ar@{~>}[ld] \ar@{~>}[rd] & & \mathrm{pt} \ar@{~>}[rd] \ar@{~>}[ld] & & \cdots \ar@{~>}[ld] \\
    & \bgm & & \bgm & & \bgm & & 
    }
\end{align*}
\end{remark}
    
    For $k\in \mathbb{Q}$, 
let $\mathcal{O}_{\X_i}(k) \in \mathrm{Pic}(\X_i)_{\mathbb{Q}}$
be the $\mathbb{Q}$-line bundle given by the pull-back of $\mathcal{O}_{\bgm}(k)$ by the map $\X_i \to \bgm$. We define $\ell\in\mathrm{CL}(\X)_{\mathbb Q}^{\vee}$ and
$\mathcal O_\X(k)\in\mathrm{Pic}(\X)_{\mathbb Q}$
by the condition
\begin{align*}
    \ell|_{\X_i}=\delta_P(\mathcal{O}_{\X_i}(-i)), \ \mathcal{O}_{\X}(k)|_{\X_i}=\mathcal{O}_{\X_i}(k).
\end{align*}
   Let $b$ be the quadratic norm induced by
$2\ch_2(\mathcal O_\X(1))\in H^4(\X,\mathbb Q)$.
With the minus sign in~(\ref{def:slope:lb}), the above choice gives
$\mu(L_i^+)=i$ and $\mu(L_i^-)=-i$. Thus 
   we have the $\Theta$-stratification of $\X$ 
   with respect to $(\ell, b)$ 
   by setting 
   \begin{align*}
       \X_{\preceq \mu}=L_{-\mu}^{-}\sqcup \cdots \sqcup L_{-1/2}^{-}
       \sqcup \X_{1/2}^{+} \sqcup L_{1/2}^+ \sqcup \cdots \sqcup L_{\mu}^+
   \end{align*}
   for each $\mu \in 1/2+\mathbb{Z}_{\geq 0}$, 
   such that $\X_{1/2}^+=\X^{\ell\text{-ss}}$
   and 
   \begin{align*}(L_{\mu}^+ \sqcup L_{-\mu}^-) \subset \X_{\preceq \mu}
   \end{align*}
   is a 
   $\Theta$-stratum. 

We use the following notation of open immersions 
\begin{align*}
    j_{(\mu-1, \mu]} \colon \X_{\prec \mu} \hookrightarrow \X_{\preceq \mu}, \ j_{\preceq \mu} \colon \X_{\preceq \mu} \hookrightarrow \X, \ 
    j \colon \X_{1/2}^+ \hookrightarrow \X.
\end{align*}
We take $\delta=\delta_P(\mathcal{O}_{\X}(-1/2))$. By Theorem~\ref{thm:sodlim}, 
we have the semiorthogonal decomposition 
\begin{align*}
    \LL(\mathcal{X}_{\preceq \mu})_{\delta}=\left\langle \mathcal{O}_{L_{\mu}^+}(-1), \mathcal{O}_{L_{-\mu}^-}, (j_{(\mu-1, \mu]})_{!}\LL(\X_{\prec \mu})_{\delta}\right\rangle. 
\end{align*}
  By Theorem~\ref{thm:cgen}, we have the semiorthogonal decomposition 
    \begin{align}\label{sod:LX}
        \LL(\X)_{\delta}=\left\langle 
        A_{\mu}, B_{\mu}, j_{!} \Coh(\X_{1/2}^+)
         : \mu \in 1/2+\mathbb{Z}_{\geq 0}
        \right\rangle. 
    \end{align}
Here $A_{\mu}$, $B_{\mu}$ are given by 
\begin{align*}
    A_{\mu}:=
        (j_{\preceq \mu})_{!}\mathcal{O}_{L_{\mu}^+}(-1), 
        \ B_{\mu} :=
        (j_{\preceq \mu})_{!} \mathcal{O}_{L_{-\mu}^-}.
\end{align*}
    One can calculate that 
    \begin{align*}
        A_{\mu}|_{\X_i}=\begin{cases}
        0, & i<\mu \\
        \mathcal{O}_{L_{\mu}^+}(-1), & i=\mu \\
        \mathcal{O}_{L_i}, & i>\mu
            \end{cases}, \quad 
     B_{\mu}|_{\X_i}=\begin{cases}
        0, & i>-\mu \\
        \mathcal{O}_{L_{-\mu}^-}, & i=-\mu \\
        \mathcal{O}_{L_i}(-1), & i<-\mu.
            \end{cases}
    \end{align*}
    Moreover $j_{!}\Coh(\X_{1/2}^+)$ is generated by $j_{!}\mathcal{O}_{\mathbb{A}^1}$ which is given by 
\begin{align*}
   (j_{!}\mathcal{O}_{\mathbb{A}^1})|_{\X_i}=\begin{cases}
       \mathcal{O}_{\X_i}, & i>0, \\
       \mathcal{O}_{\X_i}(-1), & i<0.
   \end{cases}
\end{align*}
By the above computations, one can also compute that 
 \begin{align}\label{calc2}
       & \Hom(A_{i}, j_{!}\mathcal{O}_{\mathbb{A}^1})
        = \Hom(B_i, j_{!}\mathcal{O}_{\mathbb{A}^1})
        =k[-1], \\
    \notag    &\Hom(A_i, A_j)
        =\Hom(B_i, B_j)=k\oplus k[-1], \ (i>j), \\
   \notag     &\Hom(A_i, B_j)
        =\Hom(B_j, A_i)=0.
    \end{align}

Let $\mathcal{Y}=\mathbb{A}^2/\mathbb{G}_m$ be the quotient 
stack as in Example~\ref{exam:1}. 
    By comparing with Example~\ref{exam:1}, 
    we have the following: 
    \begin{prop}\label{prop:equivX}
         There is an equivalence 
    \begin{align}\label{equiv:A2L}
    \Phi \colon 
    \Coh(\mathcal{Y}) \stackrel{\sim}{\to} \LL(\X)_{\delta}.
    \end{align}
    \end{prop}
    \begin{proof}
   Let
    \begin{align*}\cP \in \Coh^{\heartsuit}(\mathcal{Y} \times_{\mathbb{A}^1} \X)
    \end{align*}
be a rank one Cohen--Macaulay sheaf such that 
\begin{align*}
    \cP|_{\mathcal{Y} \times_{\mathbb{A}^1} \X_i}
    =\cP_i:=(x, u_i) \otimes \mathcal{O}(1/2+i, 0).
\end{align*}
Here $(x, u_i)$ is the ideal 
\begin{align*}
    (x, u_i) \subset \mathcal{O}_{\mathcal{Y} \times_{\mathbb{A}^1} \X_i}
    =k[x, y, u_i, v_i]/(xy-u_i v_i)
\end{align*}
and $\mathcal{O}(a, b)$ stands for the pull-back of 
$\mathcal{O}_{\bgm^2}(a, b)$ under the natural map 
$\Y\times_{\mathbb{A}^1}\mathcal{X} \to \bgm^2$.
Note that the sheaves $\cP_i$ glue to give $\cP$ by the isomorphism
\begin{align*}
    \cP_i|_{\mathcal{Y}\times_{\mathbb{A}^1}\X_i^{\pm}} \cong \mathcal{O}_{\mathcal{Y}}(\mp1/2+i) \in \Coh(\mathcal{Y}),
\end{align*}
hence isomorphisms
\begin{align*}
    \cP_i|_{\Y\times_{\mathbb{A}^1}\mathcal{X}_i^{-}} \cong
    \cP_{i+1}|_{\Y\times_{\mathbb{A}^1}\mathcal{X}_{i+1}^{+}}.
\end{align*}

Let $\Phi$ be the 
Fourier--Mukai functor 
\begin{align*}
    \Phi=p_{2*}(p_1^*(-)\otimes \cP) \colon \Coh(\mathcal{Y}) \to \Coh(\X)
\end{align*}
where $p_i$ are the projections. A direct calculation shows that 
\begin{align*}
    \Phi(\mathcal{O}_{\mathbb{A}^2/\mathbb{G}_m}) \cong j_{!}\mathcal{O}_{\mathbb{A}^1}, \ 
    \Phi(\mathcal{O}_{L^+}(-i))=A_{i-1/2}, \ \Phi(\mathcal{O}_{L^-}(j))=B_{j-1/2}.
\end{align*}
Recall in (\ref{calc1}), (\ref{calc2})
that the Hom spaces of these 
objects are isomorphic; 
indeed one can check that $\Phi$ induces isomorphisms of these 
Hom spaces. 
Together with the semiorthogonal decompositions (\ref{sod:A2gm}), (\ref{sod:LX}), 
the functor $\Phi$ gives an equivalence (\ref{equiv:A2L}). 
\end{proof}

\subsection{The higher-dimensional generalization: gluing of standard flops}\label{subsec:toy-higher}
\subsubsection*{Construction}
Fix $n\geq1$ and an $n$-dimensional vector space $V$.
Let $\mathbb{G}_m$ act on $V$ with weight $1$, and on its dual
$V^\vee$ with weight $-1$. 
In this and the next subsection, we suppress the dependence on $n$
and use the notation of Subsection~\ref{subsec:toy-two}. For
$i\in\frac12+\mathbb Z$, put
\begin{align*}
 V_i=\begin{cases}
 V,&i-\frac12\in2\mathbb Z,\\
 V^\vee,&i-\frac12\in2\mathbb Z+1.
 \end{cases}
\end{align*}
Note that we have
$V_{i+1}=V_i^\vee$ as $\mathbb{G}_m$-representations. Put
$\varepsilon_i=(-1)^{i-\frac12}$, the weight of $V_i$. Set
\begin{align*}
 \X_i=(V_i\oplus V_i^\vee)/\mathbb G_m,\quad
 t\cdot(x,y)=(t^{\varepsilon_i}x,t^{-\varepsilon_i}y),
\end{align*}
and write
\begin{align*}
 L_i^+=(V_i\oplus0)/\mathbb G_m,\quad
 L_i^-=(0\oplus V_i^\vee)/\mathbb G_m,\quad
 \X_i^\pm=\X_i\setminus L_i^\pm.
\end{align*}
The open substacks $\X_i^{\pm} \subset \mathcal{X}_i$ are varieties, and given by 
\begin{align*}
 \X_i^-&=\operatorname{Tot}_{\mathbb P(V_i)}
       \bigl(\mathcal O_{\mathbb P(V_i)}(-1)\otimes V_i^\vee\bigr),\\
 \X_i^+&=\operatorname{Tot}_{\mathbb P(V_i^\vee)}
       \bigl(\mathcal O_{\mathbb P(V_i^\vee)}(-1)\otimes V_i\bigr).
\end{align*}
There is a diagram of birational maps  
\begin{align*}
    \xymatrix{
    \X_i^+ \ar[dr] \ar@{.>}[rr] & & \X_i^- \ar[dl] \\
    &\X_i^0, &
    }
\end{align*}
where $\X_i^0=(V_i\oplus V_i^{\vee})\ssslash \mathbb{G}_m$. 
The above diagram is 
called a \textit{standard flop}. 

The exchange of the factors induces the isomorphism
\begin{align}\label{toy:glue-n}
 f_i\colon \X_i^-\xrightarrow{\cong}\X_{i+1}^+,
 \quad (x,y)\mapsto(y,x).
\end{align}
These isomorphisms give a gluing of $\mathcal{X}_i$ for $i\in 1/2 +\mathbb{Z}$ to define a smooth non-quasi-compact stack
\begin{align*}
 \X=\bigcup_{i\in\frac12+\mathbb Z}\X_i.
\end{align*}

The two open substacks $\mathcal{X}_i$ and $\mathcal{X}_j$ for $\lvert i-j \rvert \geq 2$
 meet only on the locus $(x, y)$ where $x\neq 0$ and $y\neq 0$. We consider the following open substack of $\mathcal{X}$
\begin{align}\label{toy:Yn}
 \X^{\mathrm{ss}}=\X_{1/2}^+=\X_{-1/2}^-
 =\operatorname{Tot}_{\mathbb P(V^\vee)}
       \bigl(\mathcal O_{\mathbb P(V^\vee)}(-1)\otimes V\bigr).
\end{align}
As we will see below, the above open substack is the semistable locus with respect to some $\ell \in \mathrm{CL}(\mathcal{X})_{\mathbb{Q}}^{\vee}$.
For $n=1$, it recovers the example in the previous subsection. 

\subsubsection*{Line bundles and numerical data on $\mathcal{X}$}
We write $\mathcal O_{\X_i}(k)$ for the line bundle associated with the
character $t\mapsto t^k$. These line bundles glue
to a line bundle 
\begin{align*}
    \mathcal{O}_{\X}(k) \in \mathrm{Pic}(\X).
\end{align*}

The points in $\mathcal{X}$ with nontrivial stabilizers are the origins
$o_i\in\X_i$.
For $a\in\mathbb Z$, let $\nu_{i, a} \colon \bgm \to \X$ be given by
\begin{align*}
    \mathrm{pt} \mapsto o_i, \ \mathbb{G}_m \ni t \mapsto t^a \in \mathbb{G}_m.
\end{align*}
We define $\ell, \delta \in \mathrm{CL}(\mathcal{X})_{\mathbb{Q}}^{\vee}$
and a quadratic norm $b$ on $\mathrm{CL}(\mathcal{X})_{\mathbb{Q}}$ by
\begin{align}\label{toy:data-n}
    \ell(\nu_{i, a})=-i\varepsilon_i a, \quad b(\nu_{i, a})=a^2, \quad \delta(\nu_{i, a})=-\frac{n}{2}\varepsilon_i a.
\end{align}

\begin{remark}\label{rmk:nonglobal}
The numerical data $\ell$ above do not come from a global $\mathbb{Q}$-line bundle on $\X$ when
$n\geq2$. Indeed for $n\geq 2$, 
we have the isomorphism 
\begin{align*}
    \mathrm{Pic}(\mathcal{X}) \stackrel{\cong}{\to} \mathrm{Pic}(\mathcal{X}^{\mathrm{ss}})=\mathbb{Z}
\end{align*}
and the generator $\mathcal O_{\X}(1)$ has $\nu_{i, a}$-weight $a$.
Since $-i\varepsilon_i$ is not constant in $i$, no rational multiple
of this generator realizes $\ell$.
\end{remark}

\subsubsection*{The $\Theta$-stratification}
We set $I_{>0}=\frac12+\mathbb Z_{\geq0}$ and
\begin{align*}
 \X_{\preceq0}=\X^{\mathrm{ss}},\quad
 \X_{\preceq\mu}=\bigcup_{-\mu\leq i\leq\mu}\X_i
       \quad(\mu\in I_{>0}),\quad
 \X_{\prec\mu}=\bigcup_{\mu'<\mu}\X_{\preceq\mu'},
\end{align*}
where $\mu'$ ranges over $\{0\}\cup I_{>0}$. 
We have the following lemma: 
\begin{lemma}\label{lem:theta}
The numerical data $(\ell, b)$ define the $\Theta$-stratification 
\begin{align}\label{toy:strata-n}
 \X=\bigcup_{\mu\in I_{>0}}\X_{\preceq\mu}, \  \mathcal{S}_{\mu}:=\X_{\preceq\mu}\setminus \X_{\prec\mu}=L_\mu^+\sqcup L_{-\mu}^-
\end{align}
with centers $B\mathbb G_m\sqcup B\mathbb G_m$ at each $\mu>0$, and 
$\X^{\ell\text{-ss}}=\X^{\text{ss}}$. Moreover each $\Theta$-stratum 
$\mathcal{S}_{\mu} \hookrightarrow \X_{\preceq \mu}$ satisfies Assumption~\ref{assum:dag}. 
\end{lemma}
\begin{proof}
For $i\in 1/2+\mathbb{Z}$ and $x\in L_i^{\pm}$, there is a map 
\begin{align*}
    f \colon \Theta \to L_i^{\pm}, \ f(1)=x
\end{align*}
such that $f(0) \colon \bgm \to L_i^{\pm}$ maps to 
the origin and with $\mathbb{G}_m \to \mathbb{G}_m$, 
$t\mapsto t^{\pm \varepsilon_i}$. 
Then we have 
\begin{align*}
    \ell(f(0))=-i \varepsilon_i \cdot (\pm \varepsilon_i)=\mp i, \ 
    \mu(f)=\pm i.
\end{align*}
It follows from the stratification 
\begin{align*}
    \X=\X^{\mathrm{ss}} \sqcup L_{\pm 1/2}^{\pm} \sqcup L_{\pm 3/2}^{\pm}\sqcup \cdots 
\end{align*}
that $\X^{\ell\text{-ss}}\subset\X^{\mathrm{ss}}$.

Conversely, let $f \colon \Theta \to \X$ be a map such that 
$f(1) \in \X^{\mathrm{ss}}$
and $f(0) \colon \bgm \to \X$ is non-degenerate. Then we have 
\begin{align*}f(1) \in L_{1/2}^- \mbox{ or }f(1) \in L_{-1/2}^+, \ 
\ell(f(0))>0.
\end{align*}
Therefore $\X^{\mathrm{ss}}\subset \X^{\ell\text{-ss}}$, and we conclude
$\X^{\ell\text{-ss}}=\X^{\mathrm{ss}}$.
The other statements are clear from the construction. 
\end{proof}

\begin{prop}\label{prop:toy-n}
For the data~(\ref{toy:data-n}), $\IndL(\X)_\delta$ is compactly
generated. Its subcategory of compact objects has a semiorthogonal decomposition
\begin{align}\label{toy:sod-n}
 \LL(\X)_\delta=
 \left\langle A_\mu,B_\mu,
          j_!\Coh(\X^{\mathrm{ss}}):\mu\in\tfrac12+\mathbb Z_{\geq0}\right\rangle,
\end{align}
where the order of the semiorthogonal decomposition is given by the decreasing order of $\mu$.

Here, writing
$j_\mu\colon \X_{\preceq\mu}\hookrightarrow\X$ and
$j\colon \X^{\mathrm{ss}}\hookrightarrow\X$, we have
\begin{align*}
 A_\mu=(j_\mu)_!\mathcal O_{L_\mu^+}(-n\varepsilon_\mu),\quad
 B_\mu=(j_\mu)_!\mathcal O_{L_{-\mu}^-}.
\end{align*}
The objects $A_\mu$ and $B_\mu$ are exceptional, and $A_\mu$ and
$B_\mu$ are orthogonal for the same $\mu$.
The functors $(j_\mu)_!$ and $j_!$ are the fully faithful left adjoints
of restriction on the corresponding limit categories.
\end{prop}
\begin{proof}
We apply Theorem~\ref{thm:cgen} to the $\Theta$-stratification in Lemma~\ref{lem:theta}.
For the closed immersions $p_i^\pm\colon L_i^\pm\hookrightarrow\X_i$,
the determinant of the relative cotangent complex has weight
$\mp n\varepsilon_i$ at $o_i$: the conormal bundle has weight
$\pm\varepsilon_i$, and occurs in
cohomological degree $-1$. 
Therefore the $\mathbb{G}_m$-weights of $\delta_{\mu}$ are 
\begin{align*}
 -\frac n2\varepsilon_\mu+\frac12(-n\varepsilon_\mu)
 =-n\varepsilon_\mu\quad\text{on }L_\mu^+,
 \quad
 -\frac n2\varepsilon_{-\mu}+\frac12(n\varepsilon_{-\mu})
 =0\quad\text{on }L_{-\mu}^-.
\end{align*}
Taking into account the cocharacters $t\mapsto t^{\varepsilon_\mu}$ and 
$t\mapsto t^{-\varepsilon_{-\mu}}$ determined by the $\Theta$-stratification, 
the central weights of $\delta_{\mu}$ are $-n$ and $0$,
respectively.

Then we have
$\LL(\mathfrak{Z}_{\pm\mu})_{\delta_\mu}
\simeq \mathrm{Vect}$,
and the functor $\Upsilon_{\mu}$
sends the generator to $\mathcal O_{L_\mu^+}(-n\varepsilon_\mu)$, 
$\mathcal O_{L_{-\mu}^-}$ respectively. 
Since $\X^{\mathrm{ss}}$ is a smooth scheme, we have $\LL(\X^{\mathrm{ss}})_{\delta}=\Coh(\X^{\mathrm{ss}})$.
In particular, for the open immersion $j_{\prec \mu, \preceq \mu}\colon \X_{\prec\mu}\hookrightarrow \X_{\preceq\mu}$,
Theorem~\ref{thm:sodlim} gives
\begin{align}\label{toy:finite-sod-n}
 \LL(\X_{\preceq\mu})_\delta=
 \left\langle\mathcal O_{L_\mu^+}(-n\varepsilon_\mu),
       \mathcal O_{L_{-\mu}^-},
       (j_{\prec \mu, \preceq \mu})_!\LL(\X_{\prec\mu})_\delta\right\rangle.
\end{align}
The first two components are orthogonal because their supports are
disjoint in $\X_{\preceq\mu}$. Theorem~\ref{thm:cgen} gives compact generation and~(\ref{toy:sod-n}).
Since $(j_\mu)_!$ is fully faithful, the assertion follows.
\end{proof}

\subsubsection*{The $\infty$-flop functor}
For $m\in\mathbb Z$, let
\begin{align*}
 Y^{(m)}:=\X_{m-1/2}^-=\X_{m+1/2}^+,
\end{align*}
where the identification is given by $f_{m-1/2}$. In particular,
$Y^{(0)}=\X^{\mathrm{ss}}$ is the semistable locus. For $n\geq2$, the two open
substacks $\X_i^+$ and $\X_i^-$ are related by a standard flop.
Therefore we have an infinite sequence of flops
\begin{align}\label{toy:infinite-flops}
 \begin{tikzcd}[column sep=large,ampersand replacement=\&]
 \cdots \& Y^{(-1)}\arrow[l,dashed,"\mathrm{flop}"']
 \& Y^{(0)}\arrow[l,dashed,"\mathrm{flop}"']
       \arrow[r,dashed,"\mathrm{flop}"]
 \& Y^{(1)}\arrow[r,dashed,"\mathrm{flop}"] \& \cdots .
 \end{tikzcd}
\end{align}
The flop between $Y^{(m)}$ and $Y^{(m+1)}$ is given by the two GIT
quotients of $\X_{m+1/2}$.

For $k\in \mathbb{Z}$, we have 
\begin{align*}
    \mathcal{O}_{\X_i}(k) \in \LL(\X_i)_{\delta} \mbox{ if and only if } -n \leq \varepsilon_i k \leq 0.
\end{align*}
Moreover, by Corollary~\ref{cor:jleft}, the image of $(j_{\prec \mu, \preceq \mu})_!$ is given by
the weight conditions
\begin{align*}
 -n<\varepsilon_\mu k\leq0\quad\text{at }o_\mu,\quad
 -n\leq\varepsilon_{-\mu}k<0\quad\text{at }o_{-\mu}.
\end{align*}
Based on the above observations, let 
$\mathcal W_i\subset\Coh(\X_i)$ be the thick subcategory
\begin{align*}
 \mathcal W_i=
 \begin{cases}
 \langle\mathcal O_{\X_i}(\varepsilon_i a): -n+1\leq a\leq0\rangle,&i>0,\\
 \langle\mathcal O_{\X_i}(\varepsilon_i a): -n\leq a\leq-1\rangle,&i<0.
 \end{cases}
\end{align*}
The restriction functors
\begin{align*}
 r_i^\pm\colon\mathcal W_i\stackrel{\sim}{\to}\Coh(\X_i^\pm)
\end{align*}
are equivalences by the window theorem, also see~\cite{DoSe}.
We denote the resulting flop equivalences by
\begin{align*}
 \Phi_i^{\mathrm R}&=r_i^-(r_i^+)^{-1}
       \colon\Coh(\X_i^+)\stackrel{\sim}{\to}\Coh(\X_i^-)
    \quad (i>0),\\
 \Phi_i^{\mathrm L}&=r_i^+(r_i^-)^{-1}
       \colon\Coh(\X_i^-)\stackrel{\sim}{\to}\Coh(\X_i^+)
    \quad (i<0).
\end{align*}

\begin{prop}\label{prop:toy-flop-extension}
For $F\in\Coh(\X^{\mathrm{ss}})$, the object $j_!F$ is obtained by gluing
$E_i\in\mathcal W_i$, where
\begin{align}\label{toy:flop-initial}
 E_{1/2}=(r_{1/2}^+)^{-1}F, \quad
 E_{-1/2}=(r_{-1/2}^-)^{-1}F,
\end{align}
and
\begin{align}\label{toy:flop-recursion}
 E_{i+1}&=(r_{i+1}^+)^{-1}(f_i)_*r_i^-E_i \quad (i>0),\\
 E_{i-1}&=(r_{i-1}^-)^{-1}f_{i-1}^*r_i^+E_i \quad (i<0).
 \notag
\end{align}
The gluing isomorphisms are given by 
\begin{align}\label{toy:flop-descent}
 f_i^*(E_{i+1}|_{\X_{i+1}^+})\cong E_i|_{\X_i^-}.
\end{align}
In particular, $(j_!F)|_{Y^{(m)}}$ is obtained from $F$ by $|m|$ flop
 equivalences along~(\ref{toy:infinite-flops}). 
\end{prop}
\begin{proof}
Let $j_{0,\mu}\colon \X^{\mathrm{ss}}\hookrightarrow \X_{\preceq\mu}$ be the inclusion.
For $\mu=1/2$, the description of the left adjoint $(j_{0, 1/2})_{!}$
of $(j_{0, 1/2})^*$
by the window subcategories shows that $(j_{0, 1/2})_{!}(F)$
is obtained as a gluing of $E_{1/2}$ and $E_{-1/2}$ in (\ref{toy:flop-initial}). 
The $!$-push-forward along the open immersion 
$\X_{\prec\mu} \hookrightarrow\X_{\preceq\mu}$
is similarly described using the windows; 
by induction, the objects $E_i$ for $i \in [-\mu, \mu]$
glue to give 
$(j_{0,\mu})_!F$ on $\X_{\preceq\mu}$.
By the construction of $j_{!}F$, we conclude that $j_{!}F$ is obtained as a gluing of $E_i$. 

The last assertion follows by restricting~(\ref{toy:flop-recursion})
to $Y^{(m)}$.
\end{proof}

We call the fully faithful functor
\begin{align*}
 j_!\colon\Coh(\X^{\mathrm{ss}})\hookrightarrow\LL(\X)_\delta
\end{align*}
the \emph{$\infty$-flop functor}. Namely, we apply the flop
 equivalences in~(\ref{toy:infinite-flops}) successively in both
directions, lift the resulting objects to $\mathcal W_i$, and glue
them on $\X$.

\begin{remark}\label{rmk:toy-BO}
The functors $\Phi_i^{\mathrm R}$ and $\Phi_i^{\mathrm L}$ can also
be described by the Bondal--Orlov equivalences~\cite[Section~3]{B-O2}.
For $n\geq2$, let
\begin{align*}
 \X_i^-\xleftarrow{q_i}\widetilde Y_i\xrightarrow{p_i}\X_i^+,
 \quad
 \widetilde Y_i=\operatorname{Tot}_{\mathbb P(V_i)\times\mathbb P(V_i^\vee)}
                    \mathcal O(-1,-1)
\end{align*}
be the common blow-up, and let $D_i \subset \widetilde{Y}_i$ be the exceptional divisor.
Set
\begin{align*}
 \mathrm{BO}_{i,t}=(p_i)_*\bigl(\mathcal O(tD_i)\otimes q_i^*(-)\bigr),\
 \mathrm{BO}'_{i,t}=(q_i)_*\bigl(\mathcal O(tD_i)\otimes p_i^*(-)\bigr).
\end{align*}
By~\cite[Proposition~3.1]{ADMflop}, one can check that 
\begin{align}\label{toy:BO-formulas}
 \Phi_i^{\mathrm R}&=\mathrm{BO}'_{i,n-1}
   \colon\Coh(\X_{i}^+)\stackrel{\sim}{\to}\Coh(\X_{i}^-)
   \quad (i>0),\\
 \Phi_i^{\mathrm L}&=\mathrm{BO}_{i,-1}
   \colon\Coh(\X_{i}^-)\stackrel{\sim}{\to}\Coh(\X_{i}^+)
   \quad (i<0).\notag
\end{align}
\end{remark}

\subsection{The gluing of Mukai flops}\label{subsec:toy-moment-zero}
We keep the notation of
Subsection~\ref{subsec:toy-higher}. The functions
$h_i(x,y)=\langle x,y\rangle$ on $\X_i$ are compatible with~(\ref{toy:glue-n}). They therefore define a map
\begin{align*}
 h\colon\X\to\mathbb A^1.
\end{align*}
We consider its derived zero locus 
\begin{align*}
 \mathfrak X_0:=\X\times_{\mathbb A^1}\{0\},
 \quad
 \mathfrak X_{i,0}:=\X_i\times_{\mathbb A^1}\{0\}.
\end{align*}
The latter derived stack is represented by the $\mathbb{G}_m$-equivariant cdga
\begin{align}\label{cdga}
 \bigl(k[V_i\oplus V_i^\vee]\otimes k[\epsilon],d\bigr),
 \quad \deg \epsilon=-1,\quad d\epsilon=\langle x,y\rangle,
\end{align}
where $\epsilon$ has trivial $\mathbb{G}_m$-weight. 

Write $\X_{i,0}^\pm=\X_i^\pm\times_{\mathbb A^1}\{0\}$.
By the Euler sequence 
\begin{align*}
0\to \Omega_{\mathbb{P}(V_i)} \to 
 \mathcal O_{\mathbb P(V_i)}(-1)\otimes V_i^\vee
 \longrightarrow\mathcal O_{\mathbb P(V_i)}\to 0
\end{align*}
we have the isomorphisms
\begin{align*}
 \X_{i,0}^-\cong \Omega_{\mathbb P(V_i)},\quad
 \X_{i,0}^+\cong \Omega_{\mathbb P(V_i^\vee)}.
\end{align*}
In particular we have 
\begin{align*}
 \mathfrak X_0^{\mathrm{ss}}:=\X^{\mathrm{ss}}\times_{\mathbb A^1}\{0\}
 \cong \Omega_{\mathbb P(V^\vee)}.
\end{align*}
We have the diagram of birational maps, called a \textit{Mukai flop} 
\begin{align*}
    \xymatrix{
    \X_{i, 0}^+ \ar[dr] \ar@{.>}[rr] & & \X_{i, 0}^- \ar[dl] \\
    &\X_{i, 0}^0, &
    }
\end{align*}
where $\mathfrak{X}_{i, 0}\to \mathcal{X}_{i, 0}^0$ is the good moduli space. 

Note that 
\begin{align*}
 \{0\}\times_{\mathbb A^1}\{0\}
   =\operatorname{Spec}k[\epsilon]=k[-1].
\end{align*}
The derived $\Theta$-strata on $\mathfrak X_{i, 0}$ and their centers are given by 
\begin{align*}
 \mathfrak S_{i,0}^\pm=L_i^\pm\times k[-1],\quad
 \mathfrak Z_{i,0}^\pm=B\mathbb G_m\times k[-1].
\end{align*}
These identifications are obtained by taking the attracting and
fixed loci of the $\mathbb{G}_m$-equivariant cdga (\ref{cdga}). 

We set
\begin{align*}\mathfrak X_{0,\preceq\mu}=\X_{\preceq\mu}\times_{\mathbb A^1}\{0\}, \ 
\mathfrak X_{0,\prec\mu}=\X_{\prec\mu}\times_{\mathbb A^1}\{0\}.
\end{align*}
The data~(\ref{toy:data-n}) define the same $\Theta$-stratification on $\mathfrak{X}_0$, 
with semistable locus $\mathfrak X_0^{\mathrm{ss}}$ and two derived $\Theta$-strata $\mathfrak S_{\mu,0}^+$ and $\mathfrak S_{-\mu,0}^-$
at every $\mu\in\frac12+\mathbb Z_{\geq0}$.
The following proposition follows as in Proposition~\ref{prop:toy-n}, and we omit details. 

\begin{prop}\label{prop:toy-moment-zero}
The dg-category $\IndL(\mathfrak X_0)_\delta$
is compactly generated, and its subcategory of compact objects admits a
semiorthogonal decomposition
\begin{align}\label{toy:sod-moment-zero}
 \LL(\mathfrak X_0)_\delta
 =\left\langle\mathcal A_{\mu,0},
       \mathcal B_{\mu,0},
       j_!\Coh\bigl(\Omega_{\mathbb P(V^\vee)}\bigr)
       :\mu\in\tfrac12+\mathbb Z_{\geq0}\right\rangle,
\end{align}
where the order of the semiorthogonal decomposition is given by the decreasing order of $\mu$.

Here $j\colon \mathfrak X_0^{\mathrm{ss}}\hookrightarrow
\mathfrak X_0$, and both $\mathcal A_{\mu,0}$ and
$\mathcal B_{\mu,0}$ are equivalent to $\Coh(k[-1])$.
The two components at the same $\mu$ are mutually orthogonal.
More explicitly, let $j_\mu\colon \mathfrak X_{0,\preceq\mu}\hookrightarrow
\mathfrak X_0$, and set
\begin{align*}
 A_{\mu,0}=(j_\mu)_!\mathcal O_{L_\mu^+}(-n\varepsilon_\mu),\quad
 B_{\mu,0}=(j_\mu)_!\mathcal O_{L_{-\mu}^-}.
\end{align*}
Then
$\mathcal{A}_{\mu, 0}$, $\mathcal{B}_{\mu, 0}$ are split-generated by 
$A_{\mu, 0}$, $B_{\mu, 0}$ respectively. 
\end{prop}

\subsubsection*{The $\infty$-Mukai-flop functor}
Assume $n\geq2$. In the notation of~\eqref{toy:infinite-flops}, put
\begin{align*}
 Y_0^{(m)}:=Y^{(m)}\times_{\mathbb A^1}\{0\}
 =\X_{m-1/2,0}^-=\X_{m+1/2,0}^+,
 \quad m\in\mathbb Z.
\end{align*}
Thus $Y_0^{(0)}=\mathfrak X_0^{\mathrm{ss}}=\Omega_{\mathbb P(V^\vee)}$, and 
we have an infinite sequence of Mukai flops
\begin{align}\label{toy:infinite-Mukai-flops}
 \begin{tikzcd}[column sep=large,ampersand replacement=\&]
 \cdots \& Y_0^{(-1)}\arrow[l,dashed,"\text{Mukai flop}"']
 \& Y_0^{(0)}\arrow[l,dashed,"\text{Mukai flop}"']
       \arrow[r,dashed,"\text{Mukai flop}"]
 \& Y_0^{(1)}\arrow[r,dashed,"\text{Mukai flop}"] \& \cdots .
 \end{tikzcd}
\end{align}
For $n=2$ the birational maps are isomorphisms, but the associated
derived equivalences can still be nontrivial.

Recall derived equivalences under Mukai flops
by Kawamata~\cite[Section~5]{MR1949787} and Namikawa~\cite{Na}
(also see~\cite{Hara} for a description of these equivalences in terms of non-commutative crepant
resolutions):
\begin{align*} \mathrm{KN}_{i,t}\colon\Coh(\X_{i,0}^-)\stackrel{\sim}{\to}\Coh(\X_{i,0}^+).
\end{align*}
We also denote by $\mathrm{KN}'_{i,t}$ the equivalence in the reverse direction, with the
normalization of~\cite[Definition~4.4 and Theorem~4.5]{ADMflop}.
Their Fourier--Mukai kernels are the restrictions of the kernels of
$\mathrm{BO}_{i,t}$ and $\mathrm{BO}'_{i,t}$ in
Remark~\ref{rmk:toy-BO}; see~\cite[proof of Proposition~4.6]{ADMflop}.
Consequently, the choices in~\eqref{toy:BO-formulas} induce
\begin{align}\label{toy:KN-formulas}
 \Psi_i^{\mathrm R}&=\mathrm{KN}'_{i,n-1}
   \colon\Coh(\X_{i,0}^+)\stackrel{\sim}{\to}\Coh(\X_{i,0}^-)
   \quad (i>0),\\
 \Psi_i^{\mathrm L}&=\mathrm{KN}_{i,-1}
   \colon\Coh(\X_{i,0}^-)\stackrel{\sim}{\to}\Coh(\X_{i,0}^+)
   \quad (i<0).\notag
\end{align}

The functor
\begin{align*}
 j_!\colon\Coh(\mathfrak X_0^{\mathrm{ss}})\hookrightarrow
          \LL(\mathfrak X_0)_\delta
\end{align*}
is therefore an \emph{$\infty$-Mukai-flop functor}.
More precisely, for $F\in\Coh(\mathfrak X_0^{\mathrm{ss}})$ and $m\geq1$, its
restrictions are
\begin{align}\label{toy:KN-iteration}
 (j_!F)|_{Y_0^{(m)}}
   &\simeq\Psi_{m-1/2}^{\mathrm R}\circ\cdots\circ
          \Psi_{1/2}^{\mathrm R}(F),\\
 (j_!F)|_{Y_0^{(-m)}}
   &\simeq\Psi_{-m+1/2}^{\mathrm L}\circ\cdots\circ
          \Psi_{-1/2}^{\mathrm L}(F),\notag
\end{align}
where the gluing maps $f_i$ are
understood between successive factors.

\subsubsection*{The case of $n=1$}
For $n=1$ in our notation, the semistable open $\mathfrak X_0^{\mathrm{ss}}$ is a point.
Using the notation of Subsection~\ref{subsec:toy-two}, set
\begin{align*}
 \mathcal Y_0:=\mathcal Y\times_{\mathbb A^1}\{0\}
   =\{xy=0\}/\mathbb G_m,
 \quad
 \X_0:=\X\times_{\mathbb A^1}\{0\}
   =\mathfrak X_0.
\end{align*}
Here $\mathcal Y_0$ is the nodal quotient stack, whereas
$\mathfrak X_0^{\mathrm{ss}}=\operatorname{Spec}k$ denotes the central open in $\X_0$.
Then similarly to Proposition~\ref{prop:equivX}, one can show that there is an equivalence 
\begin{align}\label{toy:equiv-zero-rank-one}
 \Phi_0\colon\Coh(\mathcal Y_0)
       \stackrel{\sim}{\longrightarrow}\LL(\X_0)_\delta.
\end{align}
Its kernel object is the restriction $\cP_0$ of the kernel $\cP$
in Proposition~\ref{prop:equivX} to
$\mathcal Y_0\times\X_0$.

\providecommand{\bysame}{\leavevmode\hbox to3em{\hrulefill}\thinspace}
\providecommand{\MR}{\relax\ifhmode\unskip\space\fi MR }
\providecommand{\MRhref}[2]{%
  \href{http://www.ams.org/mathscinet-getitem?mr=#1}{#2}
}
\providecommand{\href}[2]{#2}
\newcommand{\etalchar}[1]{$^{#1}$}
\providecommand{\bysame}{\leavevmode\hbox to3em{\hrulefill}\thinspace}
\providecommand{\MR}{\relax\ifhmode\unskip\space\fi MR }
\providecommand{\MRhref}[2]{%
  \href{http://www.ams.org/mathscinet-getitem?mr=#1}{#2}
}
\providecommand{\href}[2]{#2}


\vspace{5mm}

{\small\textsc{Yukinobu Toda: Kavli Institute for the Physics and Mathematics of the Universe (WPI), University of Tokyo, 5-1-5 Kashiwanoha, Kashiwa, 277-8583, Japan.}\par}
{\small\textsc{Inamori Research Institute for Science, 620 Suiginya-cho, Shimogyo-ku, Kyoto 600-8411, Japan.}\par}
{\small\textit{E-mail address:} \texttt{yukinobu.toda@ipmu.jp}\par}

\end{document}